
\catcode`\@=11

\magnification=1200 \parindent=0mm  \hfuzz=10pt
\hsize=13cm  \vsize=19cm   \hoffset=4mm		 \voffset=1cm
\pretolerance=500   \tolerance=1000   \brokenpenalty=5000

\catcode`\;=\active
\def;{\relax\ifhmode\ifdim\lastskip>\z@
\unskip\fi\kern.2em\fi\string;}

\catcode`\:=\active
\def:{\relax\ifhmode\ifdim\lastskip>\z@\unskip\fi
\penalty\@M\ \fi\string:}

\catcode`\!=\active
\def!{\relax\ifhmode\ifdim\lastskip>\z@
\unskip\fi\kern.2em\fi\string!}

\catcode`\?=\active
\def?{\relax\ifhmode\ifdim\lastskip>\z@
\unskip\fi\kern.2em\fi\string?}

\def\^#1{\if#1i{\accent"5E\i}\else{\accent"5E #1}\fi}
\def\"#1{\if#1i{\accent"7F\i}\else{\accent"7F #1}\fi}
\frenchspacing


\catcode`\@=11
\def\system#1{\left\{\null\,\vcenter{\openup1\jot\m@th
\ialign{\strut\hfil$##$&$##$\hfil&&\enspace$##$\enspace&
\hfil$##$&$##$\hfil\crcr#1\crcr}}\right.}
\catcode`\@=12

 \def\Fhd#1#2{\smash{\mathop{\hbox to 14mm{\rightarrowfill}}
\limits^{\scriptstyle#1}_{\scriptstyle#2}}}

\def\Fhg#1#2{\smash{\mathop{\hbox to 14mm{\leftarrowfill}}
\limits^{\scriptstyle#1}_{\scriptstyle#2}}}

 \def\fhd#1#2{\smash{\mathop{\hbox to 8mm{\rightarrowfill}}
\limits^{\scriptstyle#1}_{\scriptstyle#2}}}

\def\fhg#1#2{\smash{\mathop{\hbox to 8mm{\leftarrowfill}}
\limits^{\scriptstyle#1}_{\scriptstyle#2}}}

\def\Fvb#1#2{\llap{$\scriptstyle#1$}\left\downarrow
\vbox to 8mm{}\right.\rlap{$\scriptstyle#2$}}

\def\fvb#1#2{\llap{$\scriptstyle#1$}\left\downarrow
\vbox to 4mm{}\right.\rlap{$\scriptstyle#2$}}

\def\fvh#1#2{\llap{$\scriptstyle#1$}\left\uparrow
\vbox to 4mm{}\right.\rlap{$\scriptstyle#2$}}

\def\diagram#1{\def\normalbaselines{\baselineskip=0pt
\lineskip=10pt\lineskiplimit=1pt}   \matrix{#1}}


\centerline{\bf SUR DES NOTIONS DE $n$-CATEGORIE ET $n$-GROUPOIDE} \par
\centerline{\bf NON-STRICTES VIA DES ENSEMBLES MULTI-SIMPLICIAUX}

\par\vskip 1cm

\centerline{\bf Zouhair Tamsamani}\par
\centerline{	\sevenrm Laboratoire de Topologie et G\'eom\'etrie (URA-CNRS
1408)}\par
\centerline{ \sevenrm Universit\'e Paul Sabatier, 118, route de Narbonne,
31062 Toulouse, France.}\par
\vskip 1cm

\centerline{\bf INTRODUCTION}
\vskip 1cm\hskip 5mm
Dans son manuscrit [1], A.Grothendieck a souvent fait r\'ef\'erence \`a la
notion
de $n$-cat\'egorie et $n$-groupoide comme moyen de formalisation en alg\`ebre
homologique et homotopique. Selon lui on devrait pouvoir associer a un espace
topologique $X$ un $n$-groupoide fondamentale $\Pi _{_{n}}(X)$ qui
g\'en\'eralise
le groupoide de Poincar\'e de $X$ et qui contiendrait l'information des groupes
d'homotopie $\pi _{_{i}}(X)$ pour $i\leq n$. Aussi on pourrait d\'efinir
une notion de
r\'ealisation g\'eom\'etrique d'un $n$-groupoide qui permettrait
par la suite de voir que tout $n$-groupoide est $n$-\'equivalent a un $\Pi
_{_{n}}(X)$.
D'autre part les $n$-cat\'egories permettrait une bonne interpretation de la
cohomologie non ab\'elienne suivant les travaux de J. Giraud [2] et L.
Breen [3].
\par\vskip 2mm\hskip 5mm
A. Grothendieck a donn\'e quelques indications de la fa\c con dont on peut
d\'efinir une
$n$-cat\'egorie, en essayant d'expliciter les diff\'erentes compositions
des $i$-fl\`eches,
les morphismes de coh\'erences et les relations aux quelles elles satisfont.
Mais le probl\`eme qui oblige les usagers de ses notions
(J. Giraud [2] - L. Breen [3] - O. Leroy [4]) \`a se restreindre au cas $n$
= 1,2 et 3 c'est
qu'\`a partir de $n\geq 3$ les donn\'ees deviennent trop nombreuses et leur
comportement se complique.\par\vskip 2mm
\hskip 5mm L'objet de cette th\`ese consiste en un premier temps, de donner
une bonne
d\'efinition formelle d'une $n$-cat\'egorie et d'un $n$-goupoide non
strictes en termes
d'ensembles multi-simpliciaux avec les propri\'et\'es qui en d\'ecoulent,
et dans un
second temps de v\'erifier avec cette nouvelle notion certains des id\'ees de
Grothendieck. Comme le mot $n$-cat\'egorie veut dire stricte dans la
litt\'erature on
appellera notre construction de $n$-cat\'egorie large
un $n$-nerf. \par\vskip 2mm\hskip 5mm
On appelle {\it 1-nerf} un foncteur contravariant
$\Phi $, de la cat\'egorie simpliciale $\Delta $ vers celle des ensembles,
qui v\'erifie la
condition suivante : \par $(\star)$ Pour tout entier non nul  $m$,
l'application :
\par\centerline{$\diagram{\Phi (m) &\Fhd{{\delta }_{[m]}}{\sim}&
{\Phi (1){\times}_{_{\Phi (0)}},\dots,{\times}_{_{\Phi (0)}}\Phi (1)}\cr
x&\Fhd{}{}&
({\delta }^{'}_{01}(x),\dots,{\delta }^{'}_{m-1,m}(x))}$} est une
bijection, o\`u pour
$i,j\in\{0,\dots,m\}$, $\delta _{_{ij}} : [1]\fhd{}{}[m]$ est l'application
qui envoie $0,1$
respectivement sur $i,j$, et $\delta ^{'}_{_{ij}}=\Phi (\delta _{_{ij}})$.
L'id\'ee de base de
notre approche de $n$-cat\'egorie non stricte est due au fait qu'un 1-nerf
est une cat\'egorie au sens usnel et inversement. En effet si on
consid\'ere $\Phi _{_{0}}$
(respectivement $\Phi _{_{1}}$) comme l'ensemble d'objets (respectivement
des fl\`eches)
de $\Phi $, on d\'efinit la composition dans $\Phi $ par $gf:=\delta
_{_{02}}^{'}(\sigma )$
o\`u $(f,g)$ est un \'el\'ement de $\Phi _{_{1}}{\times}_{\Phi _0}\Phi
_{_{1}}$ et
$\sigma $ est l'unique \'el\'ement de $\Phi (2)$ tel que $\delta
_{_{[2]}}(\sigma )=(f,g)$.
Alors $\Phi $ muni de cette composition est une cat\'egorie au sens usuel,
et inversement
le nerf d'une cat\'egorie est un ensemble simplicial qui satisfait bien la
condition
$(\star)$. Soit maintenant $\Phi :\Delta \fhd{}{}Hom(\Delta ,Ens)$
un foncteur contravariant tel que pour tout objet $m$ de $\Delta $, le foncteur
$\Phi _{_{m}}$ est un 1-nerf et la morphisme
$\delta _{_{[m]}}$ soit une \'equivalence de cat\'egories, alors $\Phi $ se
comporte comme
une 2-cat\'egorie non stricte (Th\'eor\`eme (1.4.3)). Ce qui montre qu'une
approche
recurrente dans cette direction permet de d\'efinir une $n$-cat\'egorie non
stricte comme
un $n$-nerf, mais cela n\'ecessite une bonne notion de $n$-\'equivalence
que doivent
v\'erifier les transformations naturelles du type $\delta _{_{[m]}}$. Pour
r\'esoudre ce
probl\'eme on introduit la notion de troncation qui va repr\'esenter
l'aspect local des
$n$-\'equivalences. \par\vskip 2mm\hskip 5mm
On appelle {\it $n$-pr\'e-nerf} un foncteur contravariant de la cat\'egorie
$\Delta ^{n}$ vers celle des ensembles. Un morphisme entre deux
$n$-pr\'e-nerf est une
transformation naturelle entre les foncteurs qui les repr\'esentent. Soient
$\Phi $ un
$n$-pr\'e-nerf et $i$ un entier tel que $0\leq i\leq n$. On appelle
$i$-fl\`eche de
$\Phi $ un \'el\'ement de l'ensemble $C_{_{i}}:=\Phi (I_{_{i}},0_{_{n-i}})$
(o\`u $I_{_{i}}=(1,\dots,1)$ $i$-fois) et objet une 0-fl\`eche. lorsque $i$
est tel que
$1\leq i\leq n$, on obtient deux applications  $s, b :
C_{_{i}}\fhd{}{}C_{_{i-1}}$ qui sont
respectivement les images par $\Phi $ des fl\`eches $\delta ^{k}_{0}$ et
$\delta ^{k}_{1}$ de $\Delta ^{^{n}}$ telles que :\par
\centerline{$\matrix{\delta ^{k}_{i}=I_{[1]}\times ..\times
\delta _{i}\times ..\times I_{[0]}}$\hskip 5mm et \hskip 5mm
$\diagram{[0]&\Fhd{\delta _{i}}{}&[1]\cr 0&\Fhd{}{}&i}$}\par
Les applications $s$ , $b$ sont appell\'ees respectivement source et but des
$i$-fl\`eches de $\Phi $. On appelle ($n$-$i$)-pr\'e-nerf des $i$-fl\`eches de
$\Phi $ le  ($n$-$i$)-pr\'e-nerf ${\cal F}_{_{i}}(\Phi ):= \Phi
_{_{I_{_{i}}}}$. \par
\hskip 5mm Un $n$-pr\'e-nerf $\Phi $ est {\it 1-troncable} si et seulement
si pour
tout objet $M$ de $\Delta ^{^{n-1}}$ l'ensemble simpliciale $\Phi _{_{M}}$ qui
\`a un objet $m$ de $\Delta $ fait correspondre l'ensemble
$\Phi (M,m):=\Phi _{_{M}}(m)$ est un 1-nerf. On montre qu'alors
$\Phi $ induit un ($n$-1)-pr\'e-nerf $T\Phi $ faisant associer \`a chaque
$M$ l'ensemble des classes d'isomorphismes d'objets de la cat\'egorie
$\Phi _{_{M}}$ (deux objets de $\Phi _{_{M}}$ repr\'esentent la m\^eme classe
s'ils sont respectivement source et but d'une fl\`eche inversible de $\Phi
_{_{M}}$).
On dit qu'un $n$-pr\'e-nerf $\Phi $ est $k$-troncable, o\`u $2\leq k\leq
n$, si et
seulement si $\Phi $ est ($k$-1)-troncable et $T^{^{k-1}}\Phi $ est un
1-nerf, et
on pose $T^{^{k}}\Phi =T(T^{^{k-1}}\Phi )$. Avec cette d\'efinition de
troncation on voit
appara\^itre deux notions de $k$-\'equivalences :\par\vskip 2mm\hskip 5mm
{\it Equivalence int\'erieure :} Soit $\Phi $ un $n$-pr\'e-nerf $k$-troncable,
avec $1\leq k\leq n$ et consid\'erons les applications :
\par\centerline{$\diagram{C_{_{n-k}}
&\Fhd{s,b}{} &C_{_{n-k-1}}\cr\fvb{}{t^{k}}&&\cr  T^{k}\Phi (I_{n-k})&&}$}
Soient $u$ et $v$ deux ($n$-$k$)-fl\`eches de $\Phi $ telles que $s(u)=s(v)$ et
$b(u)=b(v)$. On dit que $u$ et $v$ sont $k$-\'equivalents si et seulement
si on a :
$t^{k}(u) = t^{k}(v)$ (c.\`a.d repr\'esentent la m\^eme classe dans
$T^{k}\Phi (I_{n-k})$).
On v\'erifie facilement que la relation "$k$-\'equivalent" est une relation
d'\'equivalence
sur l'ensemble des ($n$-$k$)-fl\`eches de $\Phi $, et on appelle ce genre
de relation
l'\'equivalence int\'erieure dans $\Phi $.
\par\vskip 2mm\hskip 5mm
{\it Equivalence ext\'erieure :} Soit $F: \Phi \fhd{}{}\Psi $ un morphisme
entre deux
$n$-pr\'e-nerfs $k$-troncables. On dit que $F$ est une $k$-\'equivalence
ext\'erieure
si et seulement si pour tout entier $h$ dans $\{ 0,\dots,k\}$ , toutes
($n$-$h$-1)-fl\`eches $u, v$ de $\Phi $ et toute ($n$-$h$)-fl\`eche $ w$ de
$\Psi $
tels que : $s(w)=F(I_{n-h-1},0_{h+1})(u)$ et
$b(w)=F(I_{n-h-1},0_{h+1})(v)$, il existe une ($n$-$h$)-fl\`eche $x$ de
$\Phi $ v\'erifiant :\par \vskip 2mm
{\bf (a)} $s(x)=u$ , $b(x)=v$ et $F(I_{n-h},0_{h})(x)$ , $w$ sont
$h$-\'equivalents dans $\Psi $\par
{\bf (b)} Si une ($n$-$h$)-fl\`eche $y$ de $\Phi $ satisfait {\bf (a)},
alors $x$ et $y$ sont $h$-\'equivalents dans $\Phi $.\par\vskip 2mm
On peut representer une $k$-\'equivalence ext\'erieure par la famille des
diagrammes du genre ci-dessous, avec $0\leq h	\leq k$.
\par\centerline{$\diagram{{\forall u, v \in
\Phi (I_{n-h-1},O_{h+1})}&\Fhd{F(I_{n-h-1},O_{h+1})}{}&\Psi
(I_{n-h-1},O_{h+1})\cr
\fvh{s,b}{}&&\fvh{}{s,b}\cr {\exists x\in\Phi
(I_{n-h},O_h)}&\Fhd{}{F(I_{n-h},O_h)}&
{\forall w\in\Psi (I_{n-h},O_h)}}$}\par\vskip 2mm	\hskip 5mm
On appelle {\it $n$-nerf} ou {\it $n$-cat\'egorie large}
un $n$-pr\'e-nerf $n$-troncable
$\Phi $ v\'erifiant pour tout entier $s$ tel que $1\leq s \leq n-1$ et tout
objet $(M,m)$ de
$\Delta ^{n-s-1}\times\Delta $ les deux axiomes suivants :\par
{\bf (C1)} Le foncteur $\Phi _{_{M,0}} : {\Delta}^{s}\fhd{}{}Ens$ est
constant.\par
{\bf (C2)} La morphisme $\delta ^{^{[m]}}_{_{M}} :
\Phi _{_{M,m}}\fhd{}{}\Phi _{_{M,1}}{\times}_{_{\Phi _{_{M,0}}}}\dots.
{\times}_{_{\Phi _{_{M,0}}}}\Phi _{_{M,1}}$ est une $s$-\'equivalence
ext\'erieure.\par
\vskip 2mm\hskip 5mm On en d\'eduit que pour tout entier $s$ tel que $1\leq
s \leq n-1$,
$\Phi _{_{M}}$ poss\'ede la structure d'un ($n$-$s$)-nerf. Donc pour tout
$i$ tel que
$1\leq i \leq n$ le ($n$-$i$)-pr\'e-nerf ${\cal F}_{i}(\Phi )$ des
$i$-fl\`eches est un
($n$-$i$)-nerf. L'op\'erateur troncation est un foncteur covariant
de $n$-Cat vers ($n$-1)-Cat, o\`u $n$-Cat est la cat\'egorie dont les objets
sont les $n$-nerfs et les fl\`eches sont les morphismes.\par
\vskip 2mm\hskip 2mm {\it Conjecture :} Soit $HO$-$n$-Cat le localis\'e
 ([5] P. Gabriel-M. Zisman) de la cat\'egorie $n$-Cat par rapport aux
$n$-\'equivalences
ext\'erieures. On peut d\'efinir un ($n$+1)-nerf $n$-CAT telle que
$HO$-$n$-Cat = $T^{n}$($n$-CAT).\par\vskip 2mm\hskip 5mm
Soit ${\it C}$ une 2-cat\'egorie large usuelle, on lui associe l'ensemble
bisimplicial
${\cal N}$ d\'efini par le 2-pr\'e-nerf suivant :
\par\centerline{$\diagram{\Delta  ^{^{2}}&\Fhd{}{}&Ens\cr (m,n)&\Fhd{}{}&
{\cal N}^{^{n}}_{_{m}}}$}o\`u ${\cal N}^{^{n}}_{_{m}}$ est l'ensemble des
familles de
quadruplets $(x_{_{i}}, f_{_{ij}}^{^{\alpha }}, \lambda _{_{ij}}^{^{\alpha
\beta }},
\varepsilon_{_{ijk}}^{^{\alpha }})$ \'el\'ement du produit \par\vskip 2mm
${{\cal C}_{_{0}}\times {\cal C}_{_{1}}\times {\cal C}_{_{2}}\times {\cal
C}_{_{2}}}$,
\hskip 3mm avec \hskip 3mm$0\leq i < j < k \leq m$\hskip 3mm ,\hskip 3mm
$ 0\leq \alpha  < \beta < \gamma \leq n$  \hskip 3mm et \par\vskip 3mm
$$x_{_{i}}\Fhd{f_{_{ij}}^{^{\alpha }}}{}x_{_{j}}\hskip 1cm
f_{_{ij}}^{^{\alpha }}
\Fhd{\lambda _{_{ij}}^{^{\alpha \beta }}}{}f_{_{ij}}^{^{\beta }}\hskip 1cm
f_{_{jk}}^{^{\alpha }}f_{_{ij}}^{^{\alpha }}\Fhd
{\varepsilon_{_{ijk}}^{^{\alpha }}}{\sim}f_{_{ik}}^{^{\alpha }}$$\par\vskip 3mm
tels qu'on a les relations de coh\`erences suivantes :\par\vskip 4mm
\centerline{$\matrix{\lambda _{_{ik}}^{^{\alpha \beta }}\cdot
\varepsilon_{_{ijk}}^{^{\alpha }}&=&\varepsilon_{_{ijk}}^{^{\beta }}\cdot
(\lambda _{_{ik}}^{^{\alpha \beta }}\star\lambda _{_{ij}}^{^{\alpha \beta }})
\cr\lambda _{_{ij}}^{^{\beta \gamma }}\cdot\lambda _{_{ij}}^{^{\alpha \beta
}}&=&
\lambda _{_{ij}}^{^{\alpha \gamma }}}$}\par\vskip 2mm\hskip 2mm
${\cal N}$ est un 2-nerf qu'on appelle {\it nerf double de $C$}.
Inversement si $\Phi $ un 2-nerf $\Phi $ les ensembles
$\Phi _{_{0,0}}$ , $\Phi _{_{1,0}}$ et $\Phi _{_{1,1}}$ se comportent
respectivement
comme les ensembles d'objets, de fl\`eches et de 2-fl\`eches d'une
2-cat\'egorie large,
o\`u les compositions ainsi que le reste de la structure sont
d\'etermin\'ees par
l'\'equivalence ext\'erieure $\delta _{_{[2]}}$ et la fonctorialit\'e de
$\Phi $
(Th\'eor\`eme (1.4.2)).\par\vskip 2mm\hskip 2mm
Soit ${\cal C}$ une 2-cat\'egorie large, et consid\'erons
la 2-cat\'egorie large ${\cal C}^{'}$, qui correspond \`a son nerf double.
Les 2-cat\'egorie
${\cal C}$ et ${\cal C}$ ont les m\^emes $i$-fl\`eches, la composition
verticale
des 2-fl\`eches est la m\^eme dans ${\cal C}$ et ${\cal C}^{'}$, mais la
composition
des fl\`eches dans ${\cal C}$ n'est pas forcement la m\^eme que celle dans
${\cal C}^{'}$, ils sont isomorphes. Cela se
traduit par le fait qu'il y a un 2-foncteur large ([4] O. Leroy) entre
${\cal C}$ et
${\cal C}^{'}$ qui admet un inverse stricte (ou exacte).\par\vskip
2mm\hskip 2mm
Soit maintenant $\Phi $ un 2-nerf, et consid\'erons le nerf double
$\Psi $ de la 2-cat\'egorie large ${\cal C}$ qui correspond \`a $\Phi $.
Alors ils existent
un 2-nerf strict ${\cal S}$ et deux 2-\'equivalences ext\'erieures
$\alpha $ , $\beta $ comme suite :\par
\centerline{$\diagram{\Phi &\Fhd{\alpha }{\sim}&{\cal S}&\Fhg{\beta
}{\sim}&\Psi }$}
\par (Remarques (1.4.5)) \par\vskip 2mm\hskip 2mm
On v\'erifie de m\^eme que le nerf multiple d'ordre $n$ d'une
$n$-cat\'egorie stricte est un
$n$-nerf stricte (c.\`a.d les transformations naturelle
$\delta _{_{[m]}}$ sont des isomorphismes) (Proposition (1.3.8).
En ayant ces resultats nous proposons la notion de $n$-nerf comme reponse
\`a la
recherche d'une notion de $n$-cat\'egorie large entreprise par A.
Grothendieck [1] et
d'autres.\par\vskip 5mm
\par\vskip 2mm\hskip 5mm La notion de $n$-nerf fait penser de fa\c con
naturelle \`a la notion de $n$-groupoide o\`u on peut inverser les
$i$-fl\`eches \`a
($n$-$i$)-\'equivalence pr\`es. On appelle {\it $n$-groupoide} la donn\'ee d'un
$n$-nerf $\Phi $ telle que ; pour tout $i$ tel que  $1\leq i\leq n$ la
cat\'egorie
${\cal C}_{_{i}}(\Phi ) = T^{^{n-i}}\Phi _{_{N}}$ avec $N=I_{_{(i-1)}}$ est
un groupoide.
Soit $F:\Phi \fhd{}{}\Psi $ un morphisme entre $n$-nerfs, alors
pour tout $i\in\{1,\dots,n\}$ , $F$  induit de fa\c con naturelle un morphisme
:\par \vskip 2mm\centerline{${\cal C}_{_{i}}(F) : {\cal C}_{_{i}}(\Phi )
\Fhd{}{}{\cal C}_{_{i}}(\Psi )$}\par\vskip 2mm  On appelle
{\it $i$-\`eme groupe d'homotopie} de $\Phi $ de base un objet $f$ de
${\cal C}_{_{i}}(\Phi )$, le groupe $Aut_{_{{\cal C}_{i}(\Phi )}}(f)$ qu'on
notera par
${\pi }_{_{i}}(\Phi ,f)$. Pour tout objet $f$ de ${\cal C}_{_{j}}(\Phi )$, o\`u
$1\leq j \leq i\leq n$, on d\'esignera par ${\pi }_{_{i}}(\Phi ,f)$ le groupe
${\pi }_{_{i}}(\Phi ,I^{^{i-j}}_{_{f}})$. On pose $\pi _{_{0}}(\Phi ) =
T^{n}\Phi $
l'ensemble des classes de $n$-\'equivalence d'objets de
$\Phi $. Ces groupes d'homotopies ainsi que l'ensemble $\pi _{_{0}}(\Phi )$
caract\'erisent le $n$-groupoide $\Phi $ \`a $n$-\'equivalence pr\`es. Lorsque
$2\leq i\leq n$ on peut consid\'erer les \'el\'ements du groupe $\pi
_{_{i}}(\Phi ,f)$
comme les 2-fl\`eches d'un 2-nerf, et en appliquant la relation de Godement
dans
cette 2-cat\'egorie large correspondante on montre que $\pi _{_{i}}(\Phi
,f)$ est ab\'elien
(Th\'eor\`eme (2.2.1)). \par\vskip 2mm\hskip 5mm
D'autre part \`a chaque espace topologique $X$ on fait associer
fonctorielement un
$n$-groupoide $\Pi _{_{n}}(X)$ qui g\'en\'eralise le groupoide fondamental
de Poincar\'e
pour $n = 1$ (Th\'eor\`eme (2.3.6)), il est construit de la fa\c con
suivante :\par
Pour tout entier positif $m$ on d\'esigne par le $m$-simplexe fondamental
l'ensemble :
\par $R^{m}=\{ (t_0,\dots,t_n)$ tel que $0\leq t_i\leq 1$ et $\sum t_i
=1\}$. Pour tout
$i\in\{0,\dots,m\}$, on d\'efinit les applications :\par\vskip 4mm
\centerline{$\matrix{R^{m-1}\fhd{d^{"}_{i}}{}R^{m}\cr (t_0,..,t_{m-1})\fhd{}{}
(t_0,..,0,..,t_{_{m-1}})}$ \hskip 1cm $\matrix{R^{m+1}\fhd{\varepsilon
^{"}_{i}}{}
R^{m}\cr
(t_{0},..,t_{_{m+1}})\fhd{}{}(t_{0},..,t_{i}+t_{_{i+1}},..,t_{_{m+1}})}$}\par
\vskip 3mm
On vient de construire un foncteur covariant $\cal R$ : $\Delta \fhd{}{} Ens$,
qui envoie les applications $\delta _{i}$ et $\delta _{ij}$ vers les
op\'erateurs c\^ot\'es
et sommets d'un $m$-simplexe suivants :\par\vskip 4mm
\centerline{$\matrix{\{1\}\fhd{\delta _{i}^{"}}{}R^{m}\cr
1\fhd{}{}(0,..,0,1,0,..,0)}$
\hskip 2cm $\matrix{R\fhd{\delta _{ij}^{"}}{}R^{m}\cr
(t,s)\fhd{}{}(0,..,0,t,0,..,0,s,0,..,0)}$}
\par\vskip 3mm\hskip 5mm On appelle {\it multi-complexe singulier}
(o\`u {\it $\infty$-complexe simgulier}) d'un espace topologique $X$, la
famille
${\cal X} =(X^{n})_{n\geq 1}$ de $n$-pr\'e-nerfs $X^{n}$
d\'efinie de fa\c con r\'ecurrente pour tout $n\geq 1$ et tout objet $(M,m)$ de
$\Delta ^{n}\times \Delta $ par : \par\vskip 3mm
\centerline{$\matrix{X^{1}(m)=Hom(R^{m},X),\hskip 5mm X^{1}(0)=X
\cr\cr X^{n+1}(M,m) =\{f\in Hom(R^{m},X^{n}(M))\mid
\forall x\in R^{m} , \forall i\in \{0,\dots,m_n\}\hskip 4mm
\delta ^{'}_i (f(x)) = f_i\}}$}\par	\vskip 3mm
o\`u  $\delta ^{'}_i : X^{n+1}_{_{M,m}}\fhd{}{}X^{n+1}_{_{M,0}}=X^{n}_{_{M}}$
est l' images de l'application $\delta _i$ par le foncteur $X^{n+1}_{_{M}}$
\par\vskip 2mm
et $f_i$ est un \'el\'ement de $X^{n}_{_{M}}$ ind\'ependant de $x$. Pour
simplifier
les notations on posera dans la suite  $X_{_{M}}=X^{n}_{_{M}}$.
Soient $f$, $g$ deux \'el\'ements de $X_{_M}$ avec $M = (m_1,\dots,m_n)$,
on dit que $f$ et $g$ sont homotopes et note $\overline f =\overline g$ si
et seulement si
il existe $\gamma \in X_{_{M,1}}$ tel que $\delta _{0}^{'}(\gamma) = f$ et
$\delta _{1}^{'}(\gamma) = g$. L'homotopie dans $X_{_M}$ est une relation
d'\'equivalence , et on d\'esignera par $\overline{(X_{_M})}$ l'ensemble de ses
classes d'\'equivalences.\par \vskip 2mm\hskip 2mm
Finalement pour un espace topologique $X$ son $n$-groupoide de Poincar\'e
associ\'e est
le foncteur $\Pi _{_{n}}(X)$ : $\Delta ^{n}\fhd{}{}Ens$ d\'efini pour tout
objet $M$ de
$\Delta ^{n}$ par : \par\vskip 3mm\centerline{$\Pi _{_{n}}(X)(0_{n})=X$,
\hskip 5mm
$\Pi _{_{n}}(X)(M) = \overline{(X_{_{M}})},\hskip 5mm\Pi _{_{n}}(X)(d^{k}_i) =
\overline{X^{n}(d^{k}_i)}$\hskip 3mm et \hskip 3mm$\Pi
_{_{n}}(X)(\varepsilon ^{k}_i) =
\overline{X^{n}(\varepsilon ^{k}_i)}$}\par\vskip 2mm\hskip 2mm
 On montre que les groupes d'homotopies sup\'erieures $\pi _{_{i}}(X ,x)$
de $X$ pour $x$
dans $X$ et $1\leq i\leq n$ sont alors isomorphes aux groupes $\pi
_{_{i}}(\Pi _{_{n}}(X) ,x)$
(Th\'eor\`eme (2.4.4)). Il y a aussi un foncteur dans l'autre sens du
pr\'ec\`edent qui
\`a un $n$-groupoide $\Phi $ fait correspondre un espace topologique
(sa r\'ealisation g\'eom\'etrique) $\mid \Phi \mid$ avec une application
naturelle :\par
\vskip 2mm\centerline{${\cal F}\ :\ Hom\bigl(\mid \Phi \mid,X\bigr)\ \Fhd{}{}\
Hom\bigl(\Phi ,\Pi _{_{n}}(X)\bigr)$}\par\vskip 2mm
et comme cons\'equence de ${\cal F}$, pour tout $n$-groupoide $\Phi $ il
existe un
morphisme de $\Phi $ vers $\Pi _{_{n}}(\mid \Phi \mid)$.
Nous conjecturons que le foncteur $\mid\ \mid$ est un inverse \`a
\'equivalence pr\`es
du foncteur $\Pi _{_{n}}(\ )$ de la cat\'egorie des espaces topologiques
$n$-tronqu\'es vers la cat\'egorie des $n$-groupoides, et nous montrons
quelques parties
de cela. \par
\vskip 5mm\hskip 5mm
J'aimerais exprimer ici ma profonde gratitude \`a mon directeur de th\`ese
Carlos
SIMPSON qui m'a initi\'e \`a la recherche et aupr\`es duquel j'ai beaucoup
appris.
Ses conseils avis\'es, sa patience et sa compr\'ehention m'ont apport\'e
une aide capitale pour l'\'elaboration de cette th\`ese.\par
\vfill\eject


{\bf CHAPITRE 1 :}\par\vskip 1cm
\centerline{\bf NOTION DE $n$-NERF}\vskip 2cm\hskip 5mm
On se propose de donner une approche simpliciale \`a la notion de
$n$-cat\'egorie non
stricte, par un foncteur contravariant $\Phi : \Delta ^{n}\fhd{}{}Ens$
de la cat\'egorie produit $n$-\`eme de la cat\'egorie simpliciale $\Delta $
vers celle des
ensembles qui satisfait certains axiomes. En suite on v\'erifie que cette
approche
est en accord avec les d\'efinitions usuelles bien connues pour $n=1$
et $2$. Cette d\'efinition permet de surmonter la difficult\'e d'expliciter
les contraintes
de coh\'erences des diverses compositions qui deviennent de plus en plus
nombreuses
lorsque $n$ augmente. On vera au cas $n$=2 que notre approche bien qu'elle
est simple
d'aspect englobe toute la complexit\'e d'une 2-cat\'egorie non stricte.
Comme le mot $n$-cat\'egorie veut dire stricte dans la litt\'erature on
appellera
notre construction de $n$-cat\'egorie large un $n$-nerf.\par\vskip 5mm


{\bf (1.1).--- Notations et d\'efinitions}\par\vskip 5mm\hskip 5mm
Soit $\Delta $ la cat\'egorie dont les objets sont les ensembles finies de
la forme :\par
$[n]=\{0,\dots,n\}$, o\`u $n$ est un entier naturel, et les fl\`eches sont les
applications $\sigma : [n]\fhd{}{}[m]$ telle que $0\leq \sigma (i)\leq
\sigma (j)\leq m $  si
$0\leq i < j\leq n$. Parmis les fl\`eches de $\Delta $ on distingue la
famille d'applications
\'el\'ementaires d\'efinies pour tout $m$ entier naturel et tout
$i\in\{0,\dots,m\}$ par :\par
\centerline{$\diagram{[m-1]&\fhd{d_{i}}{}&[m]\cr j < i&\fhd{}{}&j\cr j\geq i&
\fhd{}{}&j+1}$ \hskip 1cm$\diagram{[m+1]&\fhd{\varepsilon _{i}}{}&[m]\cr j
\leq i&
\fhd{}{}&j\cr j > i&\fhd{}{}&j-1}$}\par\hskip 5mm Soit $\sigma :
[n]\fhd{}{}[m]$ une
fl\`eche de $\Delta$, alors $\sigma $ s'ecrit de fa\c con unique sous la
forme :
$\sigma = d_{i_{1}}...d_{i_{s}}.\varepsilon _{j_{1}}...\varepsilon _{j_{t}}$ ;
o\`u $i_{1}...i_{s}$ dans l'ordre d\'ecroissant sont les \'el\'ements de
$[m]$ qui ne sont pas
atteints par $\sigma $, et $j_{1}...j_{t}$ dans l'ordre croissant les
\'el\'ements de $[n]$
tels que $\sigma (j)=\sigma (j+1)$, avec $n-t+s=m$ ([8] P. May). Un
foncteur de $\Delta $
vers $Ens$ est compl\`etement d\'etermin\'e par les images des $d_{i}$ et
$\varepsilon _{i}$.\par\vskip 5mm \hskip 5mm
On d\'esigne par $\Delta ^{n}$ la cat\'egorie n-produit de $\Delta $, ses
objets sont les
produits cartesiens $[m_{_{1}}]\times ..\times [m_{_{n}}]$, qu'on notera
dans la suite
par : $M=(m_{_{1}},\dots,m_{_{n}})$. Une fl\`eche $\mu : M\fhd{}{}M^{'}$
est de la forme
$\mu =(\mu _{_{1}},\dots,\mu _{_{n}})$, avec $\mu _{i}\in Mor_{\Delta
}(m_{i},m^{'}_{i})$.
Comme pour $\Delta $ les fl\`eches de $\Delta ^{n}$ sont engendr\'ees par
les fl\`eches
\'el\'ementaires $d^{k}_{i}$ et $\varepsilon ^{k}_{i}$ d\'efinies pour tout
$M=(m_{_{1}},\dots,m_{_{n}}) \in Ob(\Delta ^{n})$, $1\leq k\leq n$ et
$0\leq i\leq m_{_{k}}$
par : \vskip 3mm
\centerline{$\matrix{d^{k}_{i}=I_{[m_{1}]}\times ..\times d_{i}\times
..\times I_{[m_{n}]}
\cr\varepsilon ^{k}_{i}=I_{[m_{1}]}\times ..\times \varepsilon _{i}
\times ..\times I_{[m_{n}]}}$}\par\vskip 3mm\hskip 5mm
{\bf D\'efinition (1.1.1) :} Un {\it $n$-pr\'e-nerf} est un foncteur
contravariant $\Phi $ de la cat\'egorie $\Delta ^{n}$ vers celle des
ensembles. Pour
d\'eterminer un tel foncteur il suffit de conna\^itre les images des
applications
$d^{k}_{i}$ et $\varepsilon ^{k}_{i}$. Une transformation naturelle entre deux
$n$-pr\'e-nerf sera appell\'ee un {\it morphisme}.\par\hskip 5mm
Soient $\Phi $ un ($n$+$p$)-pr\'e-nerf et $(M,N)$ un objet de
$\Delta ^{n}\times\Delta ^{p}$. On d\'esigne par $\Phi _{_{M}}$ le
$p$-p\'er-nerf d\'efini par : $\Phi _{_{M}}(N)=\Phi (M,N)$.\par\hskip 5mm
Soit $\Phi $ un ($n$+1)-pr\'e-nerf. On d\'esigne par $b^{^{M}}$, $s^{^{M}}$
les images respectives par $\Phi $ des applications
$\delta ^{^{M}}_{_{0}}$ et  $\delta ^{^{M}}_{_{1}}$ d\'efinies par :\par
\centerline{$\diagram{[0]\times M&\Fhd{\delta ^{^{M}}_{_{i}}}{}&[1]\times M
\cr (0,.\alpha _{_{1}},..,\alpha _{_{n}})&\Fhd{}{}&(i,.\alpha
_{_{1}},..,\alpha _{_{n}})}$
\hskip 1cm
o\`u
\hskip 1cm
$i = 0,1$}
\par\hskip 5mm
{\bf D\'efinition (1.1.2) :} Soit $\Phi $ un ($n$+1)-pr\'e-nerf.
On appelle {\it 2-produit fibr\'e} de $\Phi _{1}$ par rapport \`a
$\Phi _{0}$ le $n$-pr\'e-nerf $\Phi _{1}{\times}_{_{\Phi _{0}}}\Phi _{1}$,
qui \`a un
objet $M$ de $\Delta ^{^{n}}$ fait correspondre le produit fibr\'e
ensembliste des
applications : \par\vskip 4mm\centerline
{$\Phi _{_{1,M}}\ \Fhd{b^{^{M}}}{}\ \Phi _{_{0,M}}\ \Fhg{s^{^{M}}}{}\ \Phi
_{_{1,M}}$}\par
\vskip 2mm\hskip 5mm Le foncteur $\Phi _{1}{\times}_{_{\Phi _{0}}}\Phi
_{1}$ est bien
d\'efini car si  $\sigma : M\fhd{}{}N$ est une fl\`eche de $\Delta
^{^{n}}$, l'application
$\Phi _{1}(\sigma )\times\Phi _{1}(\sigma )$ est \`a valeur dans l'ensemble
$\Phi _{_{1,M}}{\times}_{_{\Phi _{_{0,M}}}}\Phi _{_{1,M}}$. En effet, soit
$(f,g)$ un
\'el\'ement de $\Phi _{_{1,N}}{\times}_{_{\Phi _{_{0,N}}}}\Phi _{_{1,N}}$
et posons
$(f^{'},g^{'}) = (\Phi _{1}(\sigma )(f),\Phi _{1}(\sigma )(g))$.
Comme la fonctorialit\'e de $\Phi $ nous permet de transformer un diagramme
commutative en un autre, on a :\par
\centerline{$\diagram{[0]\times M&\fhd{\delta ^{^{M}}_{_{i}}}{}&[1]\times M\cr
\fvb{I_{_{[0]}}\times\sigma }{}&&\fvb{}{I_{_{[1]}}\times\sigma }\cr [0]\times
N&
\fhd{}{\delta ^{^{N}}_{_{i}}}&[1]\times N}$\hskip1cm $\Rightarrow$\hskip 1cm
$\diagram{\Phi _{_{0,M}}&\fhd{\Phi (\delta ^{^{M}}_{_{i}})}{}&\Phi _{_{1,M}}\cr
\fvh{\Phi _{0}(\sigma )}{}&&\fvh{}{\Phi _{1}(\sigma )}\cr \Phi _{_{0,N}}&
\fhd{}{\Phi (\delta ^{^{N}}_{_{i}})}&\Phi _{_{1,N}}}$}
La commutativit\'e du deuxi\`eme diagramme nous permet de dire que
$(f^{'},g^{'})$ est
bien un \'el\'ement de $\Phi _{_{1,M}}{\times}_{_{\Phi _{_{0,M}}}}\Phi
_{_{1,M}}$ :\par
\vskip 3mm
\centerline{$\matrix{b^{^{M}}(f^{'})&=&\Phi ({\delta }^{^{M}}_1)[{\Phi
}_{1}(\sigma )(f)]
\cr &=&{\Phi }_{1}(\sigma )[\Phi ({\delta }^{^{N}}_1)(f)]\cr\cr &=&
{\Phi }_{0}(\sigma )(b^{^{N}})(f)\cr\cr &=&\Phi _{0}(\sigma
)(s^{^{N}})(g)\cr\cr &=&
\Phi _{1}(\sigma )[\Phi (\delta ^{^{N}}_0)(g)]\cr\cr &=&
\Phi (\delta ^{^{M}}_0)[\Phi _{1}(\sigma )(g)]\cr\cr
&=&s^{^{M}}(g^{'}).}$}\hfill  \par

\vskip 3mm\hskip 5mm
On peux d\'efinir de fa\c con similaire le $m$-produit fibr\'e de $\Phi _1$
par rapport \`a $\Phi _0$, c'est le $n$-pr\'e-nerf  $\Phi
_1{\times}_{_{\Phi _0}}
\dots{\times}_{_{\Phi _0}}\Phi _1$ qui associe \`a un objet $M$ de $\Delta
^{n}$  le
produit fibr\'e des applications : \par \centerline{$\diagram{\Phi _{_{1,M}}
{\times}_{_{\Phi _{_{0,M}}}}\dots{\times}_{_{\Phi _{_{,M}}}}\Phi _{_{1,M}}&
\Fhd{b^{^{M}}P_{_{m-1}}}{}&\Phi _{_{0,M}}&\Fhg{s^{^{M}}}{}&\Phi
_{_{1,M}}}$}\par
o\`u le premier terme est un produit de ($m$-1)-facteurs et, $P_{_{m-1}}$
est la
projection sur le ($m$-1)-\`eme facteur. Ce $m$-produit fibr\'e est aussi
\'egale au
produit fibr\'e des applications suivantes :\par
\centerline{$\diagram{\Phi _{_{1,M}}&\Fhd{b^{^{M}}}{}&\Phi _{_{0,M}}&
\Fhg{s^{^{M}}P_{_{1}}}{}&\Phi _{_{1,M}}{\times}_{_{\Phi _{_{0,M}}}}
\dots{\times}_{_{\Phi _{_{,M}}}}\Phi _{_{1,M}}}$}\par
On en d\'eduit alors que le 2-produit fibr\'e est associative.
\par\vskip 5mm\hskip 5mm{\bf D\'efinition (1.1.3) :}
Un pr\'e-nerf $\Phi  :{\Delta} \fhd{}{} Ens$ est un {\it 1-nerf} si et
seulement si pour
tout entier $m\geq 2$, l'application : \par \hskip 3cm $\diagram{\Phi _m &
\Fhd{{\delta }_{[m]}}{\sim} &{\Phi _1{\times}_{\Phi _0}\Phi _1,\dots,
{\times}_{\Phi _0}\Phi _1}\cr x&\Fhd{}{}&({\delta
}^{'}_{01}(x),\dots,{\delta }^{'}_{m-1,m}
(x))}$\par est une bijection, avec pour $i,j\in\{0,\dots,m\}$, $\delta
_{ij} : [1]\fhd{}{}[m]$
est l'application qui envoie $0,1$ respectivement sur $i,j$, et ${\delta
}^{'}_{ij}=\Phi
(\delta _{ij})$.


\vskip 5mm\hskip 5mm{\bf Proposition (1.1.4) :} {\it Un 1-nerf est une
cat\'egorie au
sens usuel et inversement.}\par\vskip 5mm
\hskip 5mm{\bf Preuve :} {\bf (1)} Soit $\Phi $ un 1-nerf.
Les fl\`eches de $\Delta $ suivantes : \par\vskip 4mm
\centerline{$[0]\ \Fhd{\delta _{_{0}},\delta _{_{1}}}{}\ [1]\ \Fhd
{\delta _{_{ij}}}{}\ [m]$\hskip 1cm et \hskip 1cm
$[1]\ \Fhd{\delta _{_{00}}}{}\ [0]$}\par\vskip 2mm
se transforment par $\Phi $ en donnant les applications :\par\vskip 4mm
\centerline{$\Phi _{_{0}}\ \Fhg{\delta ^{'}_{_{0}},\delta ^{'}_{_{1}}}{}\
\Phi _{_{1}}\
\Fhg{\delta ^{'}_{_{ij}}}{}\ \Phi _{_{m}}$\hskip 1cm et \hskip 1cm
$\Phi _{_{1}}\ \Fhg{I}{}\ \Phi _{_{0}}$}\par\vskip 2mm\hskip 5mm
Soit $L_{_{2}}$ l'inverse $\delta _{_{[2]}}$. Si on consid\'ere $\Phi _{_{0}}$
(respectivement $\Phi _{_{1}}$) comme ensemble d'objets (respectivement de
fl\`eches)
de $\Phi $, on d\'efinit le compos\'e d'un \'el\'ement $(f,g)$
de $\Phi _{_{1}}{\times}_{\Phi _0}\Phi _{_{1}}$  par $gf:=\delta
_{_{02}}^{'}(\sigma )$
o\`u $\sigma =L_{_{2}}(f,g)$. \par \vskip 2mm\hskip 5mm
{\bf (a)} Montrons que cette composition est associative. Soit $(f,g,h)$ un
\'el\'ement de
$\Phi _{_{1}}{\times}_{\Phi _0}\Phi _{_{1}}{\times}_{\Phi _0}\Phi _{_{1}}$
et posons :	\par
\vskip 2mm
\centerline{$\matrix{\sigma _{_{1}} = L_{_{2}}(f,g) \hskip 5mm \tau _{_{1}}
=L_{_{2}}(g,h)
\hskip 5mm T =L_{_{3}}(f,g,h)\cr\cr
\sigma _{_{2}} = L_{_{2}}(gf,h) \hskip 5mm \tau _{_{2}}
=L_{_{2}}(f,hg)}$}\par\vskip 2mm
Pour tout $\scriptstyle i,j,k$ tels que $\scriptstyle 0\leq i<j<k\leq 3$,
on d\'esigne par
$T_{_{ijk}}= \delta ^{'}_{_{ijk}}(T)$ la face de c\^ot\'es $\scriptstyle
i,j,k$ du 3-simplexe $T$.
Soient les fl\`eches de $\Delta $ suivantes :\par
\centerline{$\diagram{[1]&\Fhd{\delta _{_{01}},\delta _{_{12}},\delta
_{_{02}}}{}&[2]&
\Fhd{\delta _{_{ijk}}}{}&[3]&\Fhg{\delta _{_{ij}},\delta _{_{jk}},\delta
_{_{ik}}}{}&[1]}$}
\par On a alors les relations :\par\vskip 2mm
\centerline{$\left\{\matrix{\delta _{ijk}\delta _{01}&=&\delta
_{ij}\cr{\delta }_{ijk}
\delta _{12}&=&\delta _{jk}\cr{\delta }_{ijk}{\delta }_{02}&=&\delta
_{ik}}\right.$
\hskip 5mm en appliquant $\Phi $ on obtient\hskip 5mm
$\left\{\matrix{{\delta }^{'}_{01}{\delta }^{'}_{ijk}&=&{\delta }^{'} _{ij}\cr
{\delta }^{'}_{12}{\delta }^{'}_{ijk}&=&{\delta }^{'}_{jk}\cr
{\delta }^{'}_{02}{\delta }^{'}_{ijk}&=&{\delta
}^{'}_{ik}}\right.$}\par\vskip 2mm
On en d\'eduit les relations :\par\vskip 2mm
\centerline{$\matrix{\delta _{_{[2]}}(T_{_{012}})=(f,g) =\delta
_{_{[2]}}(\sigma _{_{1}})
\hskip 5mm \delta _{_{[2]}}(T_{_{123}})=(g,h) =\delta _{_{[2]}}(\tau
_{_{1}}) }$}\par\vskip 2mm
ce qui montre que $\sigma _{_{1}}=T_{_{012}}$ et $\tau _{_{1}}=T_{_{123}}$.
Et par suite
on obtient :\par\vskip 2mm
\centerline{$\matrix{\delta _{_{[2]}}(\sigma _{_{2}})=(gf,h)=\delta
_{_{[2]}}(T_{_{023}})
\hskip 5mm \delta _{_{[2]}}(\tau _{_{2}})=(f,hg)=\delta
_{_{[2]}}(T_{_{013}})\cr\cr
\sigma _{_{2}}=T_{_{023}}\hskip 5mm \tau _{_{2}}=T_{_{013}}
\cr\cr (hg)f=\delta _{_{02}}^{'}(T_{_{013}})=\delta _{_{03}}^{'}(T)=
\delta _{_{02}}^{'}(T_{_{023}})=h(gf)}$}\hfill \par\vskip 2mm
\hskip 5mm{\bf (b)} Pour tout \'el\'ement $x$ de $\Phi (0)$, montrons que
$I_{_{x}} = I(x)$
la fl\`eche identit\'e de $x$ represente l'identit\'e pour la composition
des fl\`eches de
$\Phi $. Soient $f$ une fl\`eche de $\Phi $ , $\sigma
=L_{_{2}}(I_{_{s(f)}},f)$ et posons
$\tau =\delta ^{'}_{_{001}}(f)$. Donc on a $\delta _{_{[2]}}(\sigma
)=\delta _{_{[2]}}(\tau )$,
et par cons\'equent :\par $fI_{_{s(f)}}=\delta^{'} _{_{02}}(\tau )=\delta
^{'}_{_{02}}
\delta ^{'}_{_{001}}(f)=\delta ^{'}_{_{01}}(f)=f$. On v\'erifie de m\^eme que :
$I_{_{b(f)}}f=f$.\par\vskip 2mm\hskip 5mm
D'autre part pour $(f,g)$ dans $\Phi _{_{1}}{\times}_{\Phi _0}\Phi _{_{1}}$
et $\sigma =L_{_{2}}(f,g)$ on a :\par\vskip 3mm
\centerline{$\matrix{s(gf)=\delta _{0}^{'}\delta _{02}^{'}(\sigma )=
{\delta }_{0}^{'}(\sigma )=\delta _{0}^{'}\delta _{01}^{'}(\sigma )=
\delta _{0}^{'}(I_{_{s(f)}})=s(f)\cr\cr b(gf)=\delta _{1}^{'}\delta
_{02}^{'}(\sigma )=
{\delta }_{2}^{'}(\sigma )=\delta _{1}^{'}\delta _{12}^{'}(\sigma )=
\delta _{1}^{'}(I_{_{b(g)}})=b(g)}$}\hskip 5mm\vskip 3mm
{\bf (2)} Inversement soit $\cal C$ une cat\'egorie. On d\'efinit le {\it
nerf} de $\cal C$
comme le foncteur $\Phi $ qui \`a un objet $m$ de $\Delta $ fait correspondre
l'ensemble $\Phi (m)$ d\'efini par : \par\vskip 3mm
\centerline{$\Phi (m) = \{(x_{i},f_{_{ij}})_{_{0\leq i <j\leq m}}
\mid x_{i} = s(f_{_{ij}}), x_{j} = b(f_{_{ij}})
\ et\  f_{_{jk}}f_{_{ij}} = f_{_{ik}}\}$}\par\vskip 3mm
Et pour tout $k\in\{0,..,m\}$, les images de $d_{k}$  et  $\varepsilon
_{k}$ par $\Phi $ sont
donn\'ees par :\par\centerline{$\diagram{\Phi (m)&\Fhd{\Phi (d_{k})}{}&\Phi
(m-1)\cr
(x_{i},f_{_{ij}})&\Fhd{}{}&(x_{i},f_{_{ij}})_{_{i,j\not = k}}}$\hskip 2cm
$\diagram{\Phi (m)&\Fhd{\Phi (\varepsilon _{k})}{}&\Phi (m+1)\cr
(x_{i},f_{_{ij}})&\Fhd{}{}&(x_{i},f_{_{ij}},I_{_{x_{_{k}}}})}$}
La d\'ecomposition de $\delta _{_{k,k+1}}$ en fonction des $d_{i}$ s'\'ecrit :
$\delta _{_{k,k+1}} = d_{_{m}}...d_{_{k+1}}d_{_{k-1}}...d_{_{0}}$. Donc son
image par $\Phi $
est l'application $\delta _{_{k,k+1}}^{'}$ qui envoie l'\'el\'ement
$(x_{i},f_{_{ij}})$ vers la
fl\`eche $f_{_{k,k+1}}$. Ce qui montre que l'application :\par
\centerline{$\diagram{\Phi (m)&\Fhd{\delta _{_{[m]}}}{}&
\Phi (1){\times}_{_{\Phi (0)}}...{\times}_{_{\Phi (0)}}\Phi (1)\cr
(x_{i},f_{_{ij}})&\Fhd{}{}&(f_{_{i,i+1}})}$}
est bijective. Par cons\'equent $\Phi $ est un 1-nerf.\par\vskip 5mm
\hskip 5mm{\bf Remarque et d\'efinitions :} Un morphisme entre deux 1-nerfs
est une
transformation naturelle qui correspond \`a un vrai foncteur entre les
cat\'egories
usuelle correspondantes et inversement. \par\hskip 5mm
Soient $\Phi $ un $n$-pr\'e-nerf et $i$ un entier tel que $0\leq i\leq n$.
On appelle
$i$-fl\`eche de $\Phi $ un \'el\'ement de l'ensemble $C_{_{i}}:=\Phi
(I_{_{i}},0_{_{n-i}})$
(o\`u $I_{_{i}}=(1,\dots,1)$ $i$-fois) et objet une 0-fl\`eche.
lorsque $i$ est tel que $1\leq i\leq n$, on obtient deux applications
$s, b : C_{_{i}}\fhd{}{}C_{_{i-1}}$ qui sont respectivement les images par
$\Phi $ des fl\`eches $\delta ^{k}_{0}$ et $\delta ^{k}_{1}$ de $\Delta
^{^{n}}$
telles que :\par\centerline{$\matrix{\delta ^{k}_{i}=I_{[1]}\times ..\times
\delta _{i}\times ..\times I_{[0]}}$\hskip 5mm et \hskip 5mm
$\diagram{[0]&\Fhd{\delta _{i}}{}&[1]\cr 0&\Fhd{}{}&i}$}\par
Les applications $s$ , $b$ sont appell\'ees respectivement source et but des
$i$-fl\`eches de $\Phi $. \par	\vskip 2mm\hskip 5mm
On appelle ($n$-$i$)-pr\'e-nerf des $i$-fl\`eches de $\Phi $ le
($n$-$i$)-pr\'e-nerf
${\cal F}_{_{i}}(\Phi ):= \Phi _{_{I_{_{i}}}}$.
\vskip 5mm


{\bf (1.2).--- Troncation d'un $n$-pr\'e-nerf :}\par\vskip 5mm \hskip 5mm
{\bf D\'efinition (1.2.1) :} Soit $\Phi  : \Delta \fhd{}{} Ens$  une
cat\'egorie.
On appelle {\it 1-troncation} de $\Phi $ l'ensemble des classes
d'isomorphismes d'objets de
la cat\'egorie associ\'ee. On le note par $T\Phi $ = $[\Phi (0)]^{\sim}$.
Deux \'el\'ements
$x$ et $y$ de $\Phi (0)$ repr\'esentent la m\^eme classe si et seulement
si, ils existent
deux fl\`eches $f$  et  $g$ de $\Phi $ avec  $s(f)=b(g)=x$  et
$b(f)=s(g)=y$ tels que :
$gf=I_{_{x}}$  et  $fg=I_{_{y}}$. \par\vskip 5mm\hskip 5mm
{\bf D\'efinition (1.2.2) :} Un $n$-pr\'e-nerf $\Phi  : {\Delta }^{n}
\fhd{}{} Ens$ est
dite {\it $1$-troncable} si et seulement si pour tout $N$ objet de $\Delta
^{n-1}$,
le pr\'e-nerf $\Phi _{N}$ est une cat\'egorie. \par\vskip 5mm
\hskip 5mm{\bf Proposition (1.2.3) :} {\it Un $n$-pr\'e-nerf $1$-troncable
$\Phi $
induit un foncteur contravariant $T\Phi  :\Delta ^{n-1}\fhd{}{}Ens$ qui \`a
chaque objet
$N$ de $\Delta ^{n-1}$ fait correspondre l'ensemble $T(\Phi _{N})$. Le
($n$-1)-pr\'e-nerf $T\Phi $ est appell\'e la 1-troncation de $\Phi $.}
\par\vskip 5mm\hskip 5mm{\bf Preuve :} Soit $\sigma : M\fhd{}{}N$
une fl\`eche de $\Delta ^{n-1}$ et posons $\sigma _{0}=\sigma \times
I_{[0]}$. Alors il
existe une unique application $T\Phi (\sigma )$ rendant commutatif le
diagramme suivant :
\centerline{$\diagram{\Phi (N,0)&\Fhd{\Phi (\sigma _{0})}{}&\Phi (M,0)\cr
\fvb{t}{}&&\fvb{}{t}\cr T\Phi (N)&\Fhd{}{T\Phi (\sigma )}&T\Phi (M)}$
\hskip 1cm donc
\hskip 1cm $T\Phi (\sigma )(t(x))=t[\Phi (\sigma _{0})(x)]$}
Soient $\sigma :M\fhd{}{}N$ et $\tau :N\fhd{}{}S$ deux fl\`eches de $\Delta
^{n-1}$.
Donc pour tout $x$ dans $\Phi (S,0)$, on a :
\par\centerline{$\matrix{(\tau \sigma )_{0}=\tau _{0}\sigma _{0}&et
&\Phi (\tau _{0}\sigma _{0})=\Phi (\sigma _{0})\Phi (\tau _{0})\cr\cr
T\Phi (\tau \sigma )\Bigr(t(x)\Bigr)&=&t\Bigr[\Phi \Bigr((\tau \sigma
)_{0}\Bigr)(x)\Bigr]
\cr\cr &=&t\Bigr[\Phi (\sigma _{0})\Bigr(\Phi (\tau _{0}(x)\Bigr)\Bigr]\cr\cr
&=&T\Phi (\sigma )\Bigr[t\Bigr(\Phi (\tau _{0})(x)\Bigr)\Bigr]\cr\cr
&=&T\Phi \Bigr[T\Phi (\sigma )\Bigr(t(x)\Bigr)\Bigr]\cr	\cr &=&
\Bigr[T\Phi (\sigma )T\Phi (\sigma )\Bigr]\Bigr(t(x)\Bigr)\cr\cr }$}\par
D'o\`u $T\Phi (\tau \sigma )=T\Phi (\sigma )T\Phi (\tau )$.
Par cons\'equent $T\Phi $ est un ($n$-1)-pr\'e-nerf.\hfill
\vskip 5mm\hskip 5mm {\bf D\'efinition (1.2.4) :}
Soient $\Phi  : {\Delta }^{n} \fhd{}{} Ens$ un $n$-pr\'e-nerf et $k$ un
entier tel que
$1\leq k\leq n$. On dit que $\Phi $ est {\it $k$-troncable} si et seulement
si $\Phi $ est
($k$-$1$)-troncable et le ($n$-$k$+$1$)-pr\'e-nerf $T^{k-1}\Phi  :
\Delta ^{n-k+1}\fhd{}{} Ens$ est $1$-troncable. On pose alors : $T^{k}\Phi
:=T(T^{k-1}\Phi )$.
\par\vskip 5mm\hskip 5mm
Soit $\Phi  : \Delta ^{n}\fhd{}{}Ens$ un $n$-pr\'e-nerf $k$-troncable.
Alors pour tout $h$ tel que $1\leq h\leq k<n$ et tout $M$ objet de $\Delta
^{n-h}$,
on a une application naturelle : \par
\vskip 3mm\centerline{$t^{h}(\Phi )(M)\ :\ \Phi (M,0_{h})\ \Fhd{}{}\
T^{h}\Phi (M)$}\par
\vskip 3mm qui est le compos\'e $h$ fois de l'application canonique :
$t(\Phi )(M) : \Phi (M,0)\ \fhd{}{}\  T\Phi (M)$ qui \`a chaque \'el\'ement
de $\Phi (M,0)$ fait
correspondre sa classe dans $T\Phi (M)$. Pour tout  $h$ tel que
$1\leq h < k\leq n$ on a : $T^{h}(T\Phi )\ =\ T(T^{h}\Phi )$ qui se
d\'eduit de fa\c con
r\'ecurrente de $T^{h+1}(\Phi )\ :=\ T(T^{h}\Phi )$.
Dans la suite on notera $t^{h}(\Phi )$ par $t^{h}$ s'il n'y a pas
d'ambiguit\'e.
\par\vskip 5mm\hskip 5mm
Pour tout $n$ entier naturel non nul, on d\'esigne par $n$-{\cal PN} la
cat\'egorie o\`u les
objets sont les $n$-pr\'e-nerfs $n$-troncables et les fl\`eches sont les
morphismes. Si $F : \Phi \fhd{}{} \Psi $ est une fl\`eche de $n$-{\cal PN},
alors $F$ induit une fl\`eche $TF : T\Phi  \fhd{}{} T\Psi $ de
($n$-$1$)-{\cal PN}
d\'efinie de fa\c con \`a ce que le diagramme suivant soit commutative :\par
\centerline{$\diagram{\Phi (M,0)&\Fhd{F(M,0)}{}&\Psi (M,0)\cr
\fvb{t(\Phi )(M)}{}&&\fvb{}{t(\Psi )(M)}\cr T\Phi (M)&\Fhd{}{TF(M)}&T\Psi
(M)}$}
En effet $TF$ est d\'efinie par $TF(M)(t(x)) = t(F(M,0)(x))$.\par\vskip 5mm
\hskip 5mm


{\bf Proposition (1.2.5) :} {\it Pour tout $h\in \{1,\dots,n-1\}$ , $T^{h}$
est un foncteur
covariant de $n$-{\cal PN}.}
\par\vskip 5mm\hskip 5mm{\bf Preuve :} Soient $\Phi \fhd{F}{}\Psi
\fhd{G}{}\Theta $ deux
fl\`eches de $n$-{\cal PN}, et $M$ dans $Ob(\Delta ^{n-1})$. On obtient
deux diagrammes
commutatifs :\par
\centerline{$\diagram{\Phi (M,0)&\Fhd{F(M,0)}{}&\Psi (M,0)&\Fhd{G(M,0)}{}&
\Theta (M,0)\cr \fvb{t(\Phi )}{}&&\fvb{t(\Psi )}{}&&\fvb{t(\Theta )}{}\cr
T\Phi (M)&
\Fhd{}{TF(M)}&T\Psi (M)&\Fhd{}{TG(M)}&T\Theta (M)}$}
ce qui entraine que le grand diagramme est aussi commutatif, donc \par
$T(GF) = (TG)(TF)$. En plus la commutativit\'e du diagramme : \par
\centerline{$\diagram{\Phi (M,0)&\Fhd{I_{\Phi }(M,0)}{}&\Phi (M,0)\cr
\fvb{t}{}&&\fvb{t}{}\cr T\Phi (M)&\Fhd{}{T(I_{\Phi })(M)}&T\Phi (M)}$}\par
montre que  $T(I_{\Phi }) = I_{T\Phi }$ , et par cons\'equent l'opperation
de troncation est
bien un foncteur covariant. En effectuant cette opp\'eration $h$ fois on montre
la proposition.\hfill


\vskip 5mm\hskip 5mm  {\bf Proposition (1.2.6) :} {\it Soit $\Phi $ un
$n$-pr\'e-nerf $k$-troncable, avec $2 \leq k \leq n$. Alors pour tout $h$
tel que $1\leq h < k $, on a :\par
(a) $\Phi $ est $h$-troncable.\par
(b) $T^{h}\Phi $ est un ($n$-$h$)-pr\'e-nerf ($k$-$h$)-troncable.\par
(c) $T(T^{h}\Phi ) \ =\ T^{h}(T\Phi )$.}\par
\vskip 5mm \hskip 5mm {\bf Preuve :} (a) et (b)  $\Phi $ est $k$-troncable,
donc par
d\'efinition elle est ($k$-$1$)-troncable et $T^{k-1}\Phi $ est
$1$-troncable. De m\^eme
on a, $\Phi $ est ($k$-$2$)-troncable et $T^{k-2}\Phi $ est $1$-troncable.
On en d\'eduit
que $T^{k-2}\Phi $ est $2$-troncable. Ainsi de suite, on arrive \`a $\Phi $
$1$-troncable
et $T\Phi $ est ($k$-$1$)-troncable.\par
(c) Soit h tel que $1\leq h < k $ alors :\par\vskip 3mm
\centerline{$\matrix{T(T^{h}\Phi )&=&T[T(T^{h-1}\Phi )]\cr\cr
&=&T^{2}(T^{h-1})\cr
\hbox to 2cm{\dotfill}&&\hbox to 2cm{\dotfill}\cr &=& T^{h}(T\Phi
).}$}\hfill \par


\vskip 5mm \hskip 5mm  {\bf Lemme (1.2.7) :} {\it (a) Si $\Phi $ est un
$n$-pr\'e-nerf 1-troncable, avec $1 < n$ , alors pour tout entier $m$ ,
${\Phi }_m$ est un ($n$-$1$)-pr\'e-nerf 1-troncable et $T({\Phi }_m)=(T\Phi
)_m$.
\par (b) Si $\Phi $ est un $n$-pr\'e-nerf $k$-troncable, avec $1<k < n$ ,
alors pour tout entier $m$ , ${\Phi }_m$ est un ($n$-$1$)-pr\'e-nerf
$k$-troncable
et $T^{k}({\Phi }_m)=(T^{k}\Phi )_m$.}\vskip 5mm \hskip 5mm
{\bf Preuve :} (a) $\Phi $ est $1$-troncable, si et seulement si
$\Phi _{m,m_2,\dots,m_{n-1}}$ est $1$-troncable, d'o\`u ${\Phi }_m$ est
$1$-troncable.
Soient $m$ un entier et $M\in {\Delta }^{n-2}$, alors\par\vskip 4mm
\centerline{$\matrix{T({\Phi }_m)(M)&=&[{{\Phi }_m}(M,0)]^{\sim}\cr\cr &=
&[{\Phi }(m,M,0)]^{\sim}\cr\cr &=&{T\Phi }(m,M)\cr\cr &=&{{(T\Phi )}_m}(M)}$}.
\par\vskip 4mm Donc  $T({\Phi }_m)\ =\  {(T\Phi )}_m$.\par\vskip 2mm
(b) D'apr\`es la proposition (1.2.2) $\Phi $ est un ($n$-1)-pr\'e-nerf
($k$-1)-troncable, et en lui appliqant (a) on montre que $T(\Phi _{m})$ est
1-troncable et
$(T^{2}\Phi )_{m}=T^{2}(\Phi _{m})$. On en d\'eduit que $\Phi _m$ est
$2$-troncable.
Ainsi de suite on montre (b).\hfill \par


\vskip 5mm \hskip 5mm  {\bf Proposition (1.2.8) :} {\it Soit $\Phi $ un
$n$-pr\'e-nerf $k$-troncable avec $2 \leq k \leq n$. Soient $s$ un entier et
$M$ un objet de $\Delta ^{s}$.\par
(a) Si  $1\leq s\leq n-k$  alors $\Phi _{_{M}}$ est un ($n$-$s$)-pr\'e-nerf
$k$-troncable.\par
(b) Si  $n-k\leq s< n$  alors $\Phi _{_{M}}$ est un ($n$-$s$)-pr\'e-nerf
($n$-$s$)-troncable.\par
(c) Soit $h$ entier tel que, $1\leq h\leq k$. Si $\Phi _{_{M}}$ est
$h$-troncable alors
$T^{^{h}}(\Phi _{_{M}}) = (T^{^{h}}\Phi )_{_{M}}.$} \par
\vskip 5mm \hskip 5mm {\bf Preuve :} Pour $n$ et $k$ fix\'es, v\'erifiant
$2 \leq k
\leq n$, faisons une r\'ecurrence sur $s$ :\par
(a) Le cas $s = 1$ est justifi\'e par le Lemme pr\'ec\'edent. Supposons que
(a) est vraie
pour $s\in\{1,\dots,n-k-1\}$. Soit $ M^{'}=(m,m_{_{2}},\dots,m_{_{s+1}})$
un ($s$+$1$)-uplet
d'entiers, et posons $M=(m,m_2,\dots,m_s)$, alors $\Phi _{_{M}}$ est
$k$-troncable.
En lui appliquant le Lemme (1.2.7), on aura $\Phi _{_{M^{'}}}\ =\ (\Phi
_{_{M}})_{m_{s+1}}$
est un ($n$-$s$-$1$)-pr\'e-nerf $k$-troncable.\par (b) Le cas $s=n-k$ est
justifi\'e par (a).
Un raisonement analogue au pr\'ec\'edent permet de justifi\'e l'autre
partie de (b).
\par (c) Soient $ M=(m,m_{_{2}},\dots,m_{_{s}})\in{\Delta }^{n-s}$ et $h$
un entier
tel que, $1\leq h\leq k$. Si $\Phi _{_{M}}$ est $h$-troncable alors
:\par\vskip 3mm
\centerline{$\matrix{T^{^{h}}(\Phi _{_{M}})&=&T^{^{h-1}}[T(\Phi )_{_{M}}]
\cr\cr
&=&T^{^{h-1}}(T\Phi )_{_{M}}\cr \hbox to 2cm{\dotfill}&&\hbox to
2cm{\dotfill}\cr\cr &=&
T(T^{^{h-1}}\Phi )_{_{M}}\cr\cr &=&(T^{^{h}}\Phi )_{_{M}}.}$}\hfill
\par\vskip 5mm\hskip 5mm


{\bf Remarque (1.2.9) :} Consid\`erons le diagramme d'applications suivant :
\par\centerline{$\diagram{A&\fhd{s_{1},b_{1}}{}&B\cr
\fvb{s_{2},b_{2}}{}&&\fvb{}
{s_{3},b_{3}}\cr C&\fhd{}{s_{4},b_{4}}&D}$\hskip 1cm tels que \hskip
1cm$\left\{\matrix
{s_{3}s_{1}&=&s_{4}s_{2}\cr s_{3}b_{1}&=&s_{4}b_{2} \cr
b_{3}s_{1}&=&b_{4}s_{2}\cr
b_{3}b_{1}&=&b_{4}b_{2}}\right.$}\par Pour tout entiers $m, m^{'} > 0$ , on
a une
bijection :\par\centerline{$\diagram{\displaystyle\prod_{_{B{\times}_{_{D}}
\dots{\times}_{_{D}}B}}^{m fois}(A{\times}_{_{C}}\dots{\times}_{_{C}}A)&
\Fhd{\Gamma }{\sim}&\displaystyle\prod_{_{C{\times}_{_{D}}
\dots{\times}_{_{D}}C}}^{{m}^{'} fois}(A{\times}_{_{B}}
\dots{\times}_{_{B}}A)\cr \Bigr((x_{1}^{i},\dots,x_{{m}^{'}}^{i})
\Bigr)_{_{1\leq i\leq m}}&\Fhd{}{}&\Bigr((x_{j}^{1},\dots,x_{j}^{m})
\Bigr)_{_{1\leq j\leq {m}^{'}}}}$}\par\hskip 5mm


{\bf Proposition (1.2.10) :} {\it Soit $\Phi $ un $n$-pr\'e-nerf $1$-troncable,
avec $n > 1$ . Alors pour tout entier $m\geq 1$ , le $m$-produit fibr\'e
$\Psi = {\Phi }_{1}
{\times}_{{\Phi }_0}\dots{\times}_{{\Phi }_0}{\Phi }_{1}$ est un
(n-1)-pr\'e-nerf
1-troncable, et $T\Psi = {(T\Phi )}_{1}{\times}_{{(T\Phi
)}_0}\dots{\times}_{{(T\Phi )}_0}
{(T\Phi )}_{1}$.}\par\hskip 5mm
{\bf Preuve :} $\Psi $ est $1$-troncable si et seulement si pour tout
$(M,m^{'})$ objet de
$\Delta ^{n-2}\times\Delta $, l'application : \par
\centerline{$\diagram{\Psi (M,m^{'})&\Fhd{\delta ^{M}_{[m^{'}]}}{}&
\Psi (M,1){\times}_{_{\Psi (M,0)}}\dots{\times}_{_{\Psi (M,0)}}\Psi (M,1)}$}
\par est une bijective. Or, on sait que :\par\vskip 4mm
\centerline{$\matrix{\Psi (M,m^{'})& =& \Psi (1,M,m^{'}){\times}_{_{\Psi
(0,M,m^{'})}}
\dots{\times}_{_{\Psi (0,M,m^{'})}}\Psi (1,M,m^{'}) \hskip 1cm (m fois)\cr\cr
\delta ^{^{(0,M)}}_{_{[m^{'}]}}&:&\Phi (1,M,m^{'})\fhd{}{}\Phi (1,M,1)
{\times}_{_{\Phi (1,M,0)}}\dots{\times}_{_{\Phi (1,M,0)}}\Phi (1,M,1)\hskip
1cm (m^{'} fois)
\cr\cr\delta ^{^{(0,M)}}_{_{[m^{'}]}}& :&\Phi (0,M,m^{'})\fhd{}{}\Phi (0,1,M)
{\times}_{_{\Phi (0,M,0)}}\dots{\times}_{_{\Phi (0,M,0)}}\Phi (0,M,1)
\hskip 1cm (m^{'} fois)}$}\par\vskip 4mm
Si on pose : $A = \Phi (1,M,1)$ , $B = \Phi (0,M,1)$ , $C = \Phi (1,M,0)$
et $D = \Phi (0,M,0)$ on aura, d'apr\`es la remarque (1.2.9), un diagramme
commutative :\par\centerline{$\diagram{\Psi (M,m^{'})&\Fhd{\delta
^{M}_{[m^{'}]}}{}&
\Psi (M,1){\times}_{_{\Psi (M,0)}}\dots{\times}_{_{\Psi (M,0)}}\Psi (M,1)\cr
\fvb{\delta ^{(1,M)}_{[m^{'}]}\times\dots\times\delta
^{(1,M)}_{[m^{'}]}}{\sim}&&
\Vert\cr\displaystyle\prod_{_{B{\times}_{_{D}}
\dots{\times}_{_{D}}B}}^{m fois}(A{\times}_{_{C}}\dots{\times}_{_{C}}A)&
\Fhd{\Gamma }{\sim}&\displaystyle\prod_{_{C{\times}_{_{D}}
\dots{\times}_{_{D}}C}}^{{m}^{'}
fois}(A{\times}_{_{B}}\dots{\times}_{_{B}}A)}$}\par
En effet, si $\sigma =(\sigma _{1},\dots,\sigma _{m})$ est un \'el\'ement de
$\Psi (M,m^{'})$ alors :\par\vskip 4mm
\centerline{$\matrix{\delta ^{(1,M)}_{[m^{'}]}\times\dots\times
\delta ^{(1,M)}_{[m^{'}]}(\sigma )&=&\Bigr((\delta _{_{01}}(\sigma _{i}),
\dots,\delta _{_{m^{'}-1,m^{'}}}(\sigma _{i}))\Bigr)_{_{1\leq i\leq
m}}\cr\cr\cr
\delta ^{M}_{[m^{'}]}(\sigma )&=&\Bigr((\delta _{_{j,j+1}}(\sigma _{1}),
\dots,\delta _{_{j,j+1}}(\sigma _{m^{'}}))\Bigr)_{_{0\leq j\leq m^{'}-1}}
\cr\cr\cr &=&\Gamma \Bigr[\Bigr((\delta _{_{01}}(\sigma _{i}),
\dots,\delta _{_{m^{'}-1,m^{'}}}(\sigma _{i}))\Bigr)_{_{1\leq i\leq
m}}\Bigr]\cr\cr\cr&=&
\Gamma \Bigr[\delta ^{(1,M)}_{[m^{'}]}\times\dots\times
\delta ^{(1,M)}_{[m^{'}]}(\sigma )\Bigr]}$}\par\vskip 4mm
On en d\'eduit que $\delta ^{M}_{[m^{'}]}$ est une bijection, d'o\`u $\Psi
$ est
$1$-troncable. Soit $M$un objet de $\Delta ^{n-2}$ alors :\par\vskip 4mm
\centerline{$\matrix{(T\Psi )(M)&=&[\Psi (M,0)]^{\sim}\cr\cr &=&[\Phi _{1}(M,0)
{\times}_{_{\Phi _{0}(M,0)}}\dots{\times}_{_{\Phi _{0}(M,0)}}\Phi
_{1}(M,0)]^{\sim}\cr
\cr &=&[\Phi _{1}(M,0)]^{\sim}{\times}_{_{[\Phi _{0}(M,0)]^{\sim}}}
\dots{\times}_{_{[\Phi _{0}(M,0)]^{\sim}}}[\Phi _{1}(M,0)]^{\sim}\cr\cr &=&
(T{\Phi }_{1})(M){\times}_{_{(T{\Phi }_{0})(M)}}\dots{\times}_{_{(T{\Phi
}_{0})(M)}}
(T{\Phi }_{1})(M)\cr\cr &=& {(T\Phi )}_{1}(M){\times}_{{(T\Phi )}_{0}(M)}
\dots{\times}_{{(T\Phi )}_{0}(M)}{(T\Phi )}_{1}(M)\cr \cr&=&
\Bigr({(T\Phi )}_{1}{\times}_{{(T\Phi )}_0}\dots{\times}_{{(T\Phi
)}_0}{(T\Phi )}_{1}\Bigr)(M)}$}
\par\vskip 4mm D'o\`u \hskip 1cm $T\Psi  =  {(T\Phi )}_{1}{\times}_{{(T\Phi
)}_0}
\dots{\times}_{{(T\Phi )}_0}{(T\Phi )}_{1}$.\hfill \vskip 5mm


{\bf (1.3).--- $n$-Equivalence et $n$-nerf. } \par\vskip 5mm\hskip 5mm
{\it Equivalence int\'erieure :} Soit $\Phi $ un $n$-pr\'e-nerf $k$-troncable,
avec $1\leq k\leq n$ et consid\'erons les applications :
\par\centerline{$\diagram{C_{_{n-k}}
&\Fhd{s,b}{} &C_{_{n-k-1}}\cr\fvb{}{t^{k}}&&\cr  T^{k}\Phi (I_{n-k})&&}$}
Soient $u$ et $v$ deux ($n$-$k$)-fl\`eches de $\Phi $ telles que $s(u)=s(v)$ et
$b(u)=b(v)$. On dit que $u$ et $v$ sont $k$-\'equivalents si et seulement
si on a :
$t^{k}(u) = t^{k}(v)$ (c.\`a.d repr\'esentent la m\^eme classe dans
$T^{k}\Phi (I_{n-k})$).
On v\'erifie facilement que la relation "$k$-\'equivalent" est une relation
d'\'equivalence
sur l'ensemble des ($n$-$k$)-fl\`eches de $\Phi $, et on appelle ce genre
de relation
l'\'equivalence int\'erieure dans $\Phi $.
\par\vskip 3mm\hskip 5mm
{\it Cas limite :} $(k=0 , n > 0)$  Deux $n$-fl\`eches de $\Phi $ qui ont
m\^eme source et m\^eme but sont 0-\'equivalents si et seulement si
ils sont \'egales.\par\vskip 5mm\hskip 5mm
{\it Equivalence ext\'erieure :} Soit $F: \Phi \fhd{}{}\Psi $ un morphisme
entre deux
$n$-pr\'e-nerfs $k$-troncables. On dit que $F$ est une $k$-\'equivalence
ext\'erieure
si et seulement si pour tout entier $h$ dans $\{ 0,\dots,k\}$ , toutes
($n$-$h$-1)-fl\`eches $u, v$ de $\Phi $ et toute ($n$-$h$)-fl\`eche $ w$ de
$\Psi $
tels que : $s(w)=F(I_{n-h-1},0_{h+1})(u)$ et
$b(w)=F(I_{n-h-1},0_{h+1})(v)$, il existe une ($n$-$h$)-fl\`eche $x$ de
$\Phi $ v\'erifiant :\par \vskip 2mm
{\bf (a)} $s(x)=u$ , $b(x)=v$ et $F(I_{n-h},0_{h})(x)$ , $w$ sont
$h$-\'equivalents dans $\Psi $\par
{\bf (b)} Si une ($n$-$h$)-fl\`eche $y$ de $\Phi $ satisfait {\bf (a)},
alors $x$ et $y$ sont $h$-\'equivalents dans $\Phi $.\par\vskip 2mm
On peut representer une $k$-\'equivalence ext\'erieure par la famille des
diagrammes du genre ci-dessous, avec $0\leq h	\leq k$.
\par\centerline{$\diagram{{\forall u, v \in
\Phi (I_{n-h-1},O_{h+1})}&\Fhd{F(I_{n-h-1},O_{h+1})}{}&\Psi
(I_{n-h-1},O_{h+1})\cr
\fvh{s,b}{}&&\fvh{}{s,b}\cr {\exists x\in\Phi
(I_{n-h},O_h)}&\Fhd{}{F(I_{n-h},O_h)}&
{\forall w\in\Psi (I_{n-h},O_h)}}$}\par
On dira que le diagramme ci-dessus poss\'ede la propriet\'e de
$h$-\'equivalence.
\par\vskip 5mm{\bf Cas limites :}\par\hskip 5mm
{\bf (h=k=0)} $F$ est une $0$-\'equivalence si pour tout $u, v \in \Phi
(I_{n-1},0)$ et tout
$ w\in\Psi (I_n) $ tels que : $F(I_{n-1})(u)=s(w)$ et $F(I_{n-1})(v)=b(w) $
il existe
$ x\in \Phi (I_n)$ v\'erifiant : \par
(a) $s(x) = u$ , $b(x) = v$ et $F(I_{n})(x) = w$.\par
(b) $x$ est l'unique \'el\'ement de $\Phi (I_{n})$ v\'erifiant (a).
\par\vskip 3mm \hskip 5mm
{\bf (h=k=n)} Pour tout $ w\in \Psi (0_n) $, il existe $ x\in \Phi (0_n)$
v\'erifiant : \par
(a) $F(I_{n})(x)$ , $w$ sont n-\'equivalents.\par
(b) Si $y\in\Phi (0_n)$ v\'erifie (a), alors $y$ et $x$ sont
$n$-\'equivalents dans $\Phi $.
\par\vskip 5mm\hskip 5mm


Soit $F : \Phi\fhd{}{}\Psi $ un morphisme entre deux $n$-pr\'e-nerfs
$n$-troncables. Pour tout $i$ dans $\{1,\dots,n\}$ on
d\'esigne par ${\cal C}_{_{i}}(\Phi )$ la cat\'egorie $T^{^{n-i}}\Phi
_{_{N}}$, o\`u
$N = I_{_{i-1}}$ (C'est la cat\'egorie o\`u les objets sont les
($i$-1)-fl\`eches de $\Phi $
et les morphismes sont les classes de ($n$-$i$)-\'equivalences des
$i$-fl\`eches de $\Phi $).
Soient $M=I_{_{n-i-1}}$, donc $F$ induit une transformation naturelle :
\par\vskip 2mm\centerline{$T^{i}\Phi _{_{M}}:{\cal C}_{_{n-i}}(\Phi )
\Fhd{}{}{\cal C}_{_{n-i}}(\Psi  )$}\par\vskip 2mm qui induit par la suite
une application :
\par\centerline{$\diagram{Hom_{_{{\cal C}_{_{n-i}}(\Phi
)}}(u,v)&\Fhd{G^{u,v}_{_{n-i}}}{}
&Hom_{_{{\cal C}_{_{n-i}}(\Psi )}}(u{'},v^{'})\cr t^{i}(x)&\Fhd{}{}&
T^{i}F_{_{M,1}}(t^{h}(x))=t^{i}(F_{_{M,1}}(0_{i})(x))}$}
pour tout $u,v$ objets de ${\cal C}_{_{n-i}}(\Phi )$, o\`u
$u^{'}=F_{_{M}}(0_{i+1})(u)$ et $v^{'}=F_{_{M}}(0_{i+1})(v)$.\par\vskip
3mm\hskip 5mm
{\bf Proposition (1.3.1) :} {\it $F$ est une $n$-\'equivalence ext\'erieure
si et
seulement si pour tout $i$ dans $\{1,\dots,n\}$, et tout $u, v $ dans $\Phi
(I_{i-1},0_{n-i+1})$
les applications $G^{u,v}_{_{n-i}}$ et $T^{n}F :T^{n}\Phi \fhd{}{}T^{n}\Psi
$ sont bijectives.}
\par\vskip 5mm\hskip 5mm{\bf Preuve :} Il est claire qu'on a
les \'equivalences suivantes :\par\vskip 2mm
\hskip 1cm (a) est vrai si et seulement si $G^{u,v}_{_{n-i}}$ est
surjective.\par
\hskip 1cm (b) est vrai si et seulement si $G^{u,v}_{_{n-i}}$ est
injective.\par\vskip 2mm
Lorsque $h=n$ (a) et (b) est \'equivalent \`a $T^{n}F$ est bijective.
Ce qui montre la proposition (1.3.1).\hfill \par\vskip 5mm\hskip 5mm
{\bf D\'efinition (1.3.2) :} On appelle {\it $n$-nerf} ou {\it
$n$-cat\'egorie large}
un $n$-pr\'e-nerf $n$-troncable
$\Phi $ v\'erifiant pour tout entier $s$ tel que $1\leq s \leq n-1$ et tout
objet $(M,m)$ de
$\Delta ^{n-s-1}\times\Delta $ les deux axiomes suivants :\par
{\bf (C1)} Le foncteur $\Phi _{_{M,0}} : {\Delta}^{s}\fhd{}{}Ens$ est
constant.\par
{\bf (C2)} Le morphisme $\delta ^{^{[m]}}_{_{M}} :
\Phi _{_{M,m}}\fhd{}{}\Phi _{_{M,1}}{\times}_{_{\Phi _{_{M,0}}}}\dots.
{\times}_{_{\Phi _{_{M,0}}}}\Phi _{_{M,1}}$ est une $s$-\'equivalence
ext\'erieure.\par
\vskip 5mm\hskip 5mm{\bf Remarque (1.3.3) :} D'apr\`es la d\'efinition (1.3.2),
pour tout entier $s$ dans $\{0,\dots,n-1\}$, $\Phi _{_{M}}$ poss\`ede la
structure d'un
($n$-$s$)-nerf. Donc pour tout $i$ tel que $1\leq i \leq n$ le
($n$-$i$)-pr\'e-nerf
${\cal F}_{i}(\Phi )$ des $i$-fl\`eches est un ($n$-$i$)-nerf.\par\vskip 5mm

($i$-$1$)-fl\`eches

{\bf Compositions des $i$-fl\`eches par rapport aux ($i$-$1$)-fl\`eches
:}\vskip 5mm
\hskip 5mm Soit $\Phi $ un $n$-nerf avec $n\geq 1$ et, pour tout
$i\in\{1,\dots,n\}$
posons $N = I_{_{i-1}}$. On a une ($n$-$i$)-\'equivalence ext\'erieure :\par
\centerline{$\diagram{\Phi _{_{N,2}}&\Fhd{\delta _{_{N}}^{^{[2]}}}{\sim}&
\Phi _{_{N,1}}{\times}_{_{\Phi _{_{N,0}}}}\Phi _{_{N,1}}}$}\par Or,
l'ensemble $C_{_{i}}$ des
$i$-fl\`eches et $C_{_{i-1}}$ des ($i$-$1$)-fl\`eches sont tels que : \par
$C_{_{i}} = \Phi _{_{N,1}}(0_{_{n-i}})$ et $C_{_{i-1}} = \Phi
_{_{N}}(0_{_{n-i+1}})$.
Alors $\delta _{_{N}}^{^{[2]}}$ induit une application :\par\vskip 2mm
\centerline{$\diagram{\Phi _{_{N,2}}(0_{_{n-i}})&\Fhd{\delta
_{_{N}}^{^{[2]}}(0_{_{n-i}})}{}&
C_{_{i}}{\times}_{_{C_{_{i-1}}}}C_{_{i}}}$}\par v\'erifiant :\par
{\bf (i)} $\forall (f,g)\in C_{_{i}}{\times}_{_{C_{_{i-1}}}}C_{_{i}}$\hskip
3mm il existe
\hskip 3mm$\lambda \in \Phi _{_{N,2}}(0_{_{n-i}})$ \hskip 3mm tel que :\par
\centerline{$t^{^{n-i}}\bigl(\delta _{_{N}}^{^{[2]}}(\lambda )\bigr)=
\bigl(t^{^{n-i}}(f),t^{^{n-i}}(g)\bigr)$}
\par {\bf (ii)} L'unicit\'e de $\lambda $ est \`a  ($n$-$i$)-\'equivalence
pr\`es c.\`a.d si
${\lambda }^{'}$ satisfait (i) alors \par\vskip 2mm
	\centerline{$t^{^{n-i}}(\lambda ) = t^{^{n-i}}({\lambda }^{'})$.}
\par D'apr\`es l'Axiome de choix il existe une application :
${\cal L}_{_{i}} : C_{_{i}}{\times}_{_{C_{_{i-1}}}}C_{_{i}}\fhd{}{}\Phi
_{_{N,2}}(0_{_{n-i}})$\par\vskip 2mm
qui rend commutatif le diagramme suivant :\par
\centerline{$\diagram{C_{_{i}}{\times}_{_{C_{_{i-1}}}}C_{_{i}}&\Fhd{{\cal
L}_{_{i}}}{}&
\Phi _{_{N,2}}(0_{_{n-i}})\cr \fvb{t^{^{n-i}}}{}&&\fvb{}{\delta
_{_{N}}^{^{[2]}}}\cr
T^{^{n-i}}\Phi _{_{N,1}}{\times}_{_{T^{^{n-i}}\Phi
_{_{N,0}}}}T^{^{n-i}}\Phi _{_{N,1}}&
\Fhg{}{t^{^{n-i}}}&C_{_{i}}{\times}_{_{C_{_{i-1}}}}C_{_{i}}}$}\par
L'application ${\cal L}_{_{i}}$ est au fait le choix d'un syst\`eme de
repr\'esentants des
classes de \par ($n$-$i$)-\'equivalences des \'el\'ement de $\Phi
_{_{N,2}}(0_{_{n-i}})$.
Finalement on obtient une composition des $i$-fl\`eches par rapport aux
($i$-$1$)-fl\`eches :\par\vskip 2mm
\centerline{$\diagram{C_{_{i}}{\times}_{_{C_{_{i-1}}}}C_{_{i}}&\Fhd{\delta
^{'}_{_{02}}
\circ{{\cal L}_{_{i}}}}{}&C_{_{i}}\cr (f,g)&\Fhd{}{}&g{\bullet}_{_{i}}f :=
\delta ^{'}_{_{02}}
\bigl({\cal L}_{_{i}}(f,g)\bigr)}$}\par\hskip 5mm
{\bf Composition dans la cat\'egorie ${\cal C}_{_{i}}(\Phi )$ :}
\vskip 5mm\hskip 5mm La composition des fl\`eches de la cat\'egorie
${\cal C}_{_{i}}(\Phi )$ est donn\'ee par la bijection :\par\vskip 2mm
\centerline{$T^{^{n-i}}\delta _{_{N}}^{^{[2]}} :\ T^{^{n-i}}\Phi _{_{N,2}}\
\Fhd{}{}\
T^{^{n-i}}\Phi _{_{N,1}}{\times}_{_{T^{^{n-i}}\Phi
_{_{N,0}}}}T^{^{n-i}}\Phi _{_{N,1}}$}\par
\vskip 2mm de la mani\`ere suivante, en posant \hskip 3mm
$F = T^{^{n-i}}\delta _{_{N}}^{^{[2]}}$ , $C_{_{i}}^{'}=T^{^{n-i}}\Phi
_{_{N,1}}$ :
\par\centerline{$\diagram{C_{_{i}}^{'}
{\times}_{_{{\cal C}_{_{i-1}}}}C_{_{i}}^{'}&\Fhd{T^{^{n-i}}\delta
^{'}_{_{02}}\circ
F^{^{-1}}}{}&C_{_{i}}^{'}\cr (\alpha ,\beta )&\Fhd{}{}&\alpha
{\bar{\bullet}_{_{i}}}\beta
 := T^{^{n-i}}\delta ^{'}_{_{02}}\bigl(F^{^{-1}}(\alpha ,\beta )\bigr)}$}


\par\vskip 5mm\hskip 5mm
{\bf Lemme (1.3.4) :} {\it La composition dans la cat\'egorie ${\cal
C}_{_{i}}(\Phi )$ est
compatible avec celle des $i$-fl\`eches par rapport aux ($i$-1)-fl\`eches
dans $\Phi $
c.\`a.d que pour tout $(f,g)$ dans
$C_{_{i}}{\times}_{_{C_{_{i-1}}}}C_{_{i}}$ on a :
${\bar g}{\bar{\bullet}_{_{i}}}{\bar f} = \overline{g{\bullet}_{_{i}}f}$.}\par
\vskip 5mm\hskip 5mm{\bf Preuve :} Consid\'erons le diagramme :\par
\centerline{$\diagram{C_{_{i}}{\times}_{_{C_{_{i-1}}}}C_{_{i}}&\Fhd{{\cal
L}_{_{i}}}{}&
\Phi _{_{N,2}}(0_{_{n-i}})&\Fhd{\delta
_{_{N}}^{^{[2]}}}{}&C_{_{i}}{\times}_{_{C_{_{i-1}}}}
C_{_{i}}\cr \fvb{t^{^{n-i}}}{}&(1)&\fvb{t^{^{n-i}}}{}&(2)&\fvb{t^{^{n-i}}}{}\cr
C_{_{i}}^{'}{\times}_{_{{\cal C}_{_{i-1}}}}C_{_{i}}^{'}&
\Fhg{}{F^{^{-1}}}&T^{^{n-i}}\Phi _{_{N,2}}&\Fhd{}{F}&
C_{_{i}}^{'}{\times}_{_{{\cal C}_{_{i-1}}}}C_{_{i}}^{'}}$}\par
Il est claire que le diagramme (2) et le diagramme ext\'erieur form\'e par
(1) et (2) sont
commutatif. Or, l'application $F$ est bijective ce qui entraine que le
diagramme (1)
est aussi commutatif. On en d\'eduit que pour tout $(f,g)$ dans
$C_{_{i}}{\times}_{_{C_{_{i-1}}}}C_{_{i}}$ , la relation :	\par\vskip 2mm
\centerline{$t^{^{n-i}}\bigl[{\cal L}_{_{i}}(f,g)\bigr] = F^{^{-1}}(\bar
f,\bar g)$}\vskip 2mm
ce qui entraine en appliquant \hskip 3mm$T^{^{n-i}}\delta
_{_{02}}^{^{'}}$\hskip 3mm
la relation :\par\vskip 2mm
\centerline{$T^{^{n-i}}\delta _{_{02}}^{^{'}}\bigl[t^{^{n-i}}({\cal
L}_{_{i}}(f,g))\bigr] =
\bar g{\bar{\bullet}_{_{i}}}\bar f$}\par	\vskip 2mm
D'autre par, d'apr\`es (7-1) on a la relation v\'erifi\'e par ${\cal
L}_{_{i}}$ :\par\vskip 2mm
\centerline{$T^{^{n-i}}\delta _{_{02}}^{^{'}}\bigl[t^{^{n-i}}({\cal
L}_{_{i}}(f,g))\bigr] =
t^{^{n-i}}\bigl[\delta _{_{02}}^{^{'}}({\cal L}_{_{i}}(f,g))\bigr] =
\overline{g{\bullet}_{_{i}}f}$}\hfill \par


\vskip 5mm \hskip 5mm{\bf Proposition (1.3.5) :} {\it (a) Soit $F : \Phi
\fhd{}{} \Psi $
une $k$-\'equivalence entre deux $n$-nerfs $k$-troncables avec $1\leq k\leq n$,
alors $TF : T\Phi \fhd{}{} T\Psi $ est une ($k$-1)-\'equivalence \'exterieure.
\par (b) L'op\'erateur troncation est un foncteur covariant
de $n$-Cat vers ($n$-1)-Cat, o\`u $n$-Cat est la cat\'egorie dont les objets
sont les $n$-nerfs et les fl\`eches sont morphismes.}\vskip 5mm
\hskip 5mm{\bf Preuve :}\par \hskip 5mm
(a) Soient $h\in \{0,\dots,k-1\}$ et $N = I_{_{k-h-1}}$, on veut montrer
que le diagramme
suivant poss\'ede la propriet\'e de $h$-\'equivalence :\par
\centerline{$\diagram{T\Phi (N,0_{_{h+1}})&\Fhd{TF(N,0_{_{h+1}})}{}&T\Psi
(N,0_{_{h+1}})
\cr \fvb{s,b}{}&(1)&\fvb{}{s,b}\cr T\Phi
(N,1,0_{_{h}})&\Fhd{}{TF(N,1,0_{_{h}})}&
T\Psi (N,1,0_{_{h}})}$}\par
Pour $h > 1$ on a :  $T\Phi (N,0_{_{h}}) = \Phi (N,1,0_{_{h+1}})$  et
$T\Psi (N,0_{_{h}}) = \Psi (N,0_{_{h+1}})$
\par donc le diagramme (1) correspond \`a celui form\'e par $\Phi $, $\Psi
$, $F$ et qui
v\'erifie d\'ej\`a  la propriet\'e de $h$-\'equivalence.\par
Pour $h = 0$ le diagramme (1) s'ecrit :\par
\centerline{$\diagram{\Phi (N,0_{_{2}})&\Fhd{TF(N,0)}{}&\Psi (N,0_{_{2}})
\cr \fvb{s,b}{}&&\fvb{}{s,b}\cr T\Phi (N,1)&
\Fhd{}{TF(N,1)}&T\Psi (N,1)}$}\par
Soient  $w\in\Psi (N,1,0)$  et  $u,v\in\Phi (N,0_{_{2}})$  tels que
$s(w) = F(N,0_{_{2}})(u)$  et  $b(w) = F(N,0_{_{2}})(v)$.
La commutativit\'e du diagramme suivant :
\par\centerline{$\diagram{\Phi (N,0_{_{2}})&\Fhd{TF(N)}{}&\Psi (N,0_{_{2}})
\cr \fvb{s,b}{}&(2)&\fvb{}{s,b}\cr \Phi (N,1,0)&\Fhd{}{F(N,1)}&\Psi
(N,0)}$}\par
montre qu'il existe un \'el\'ement $x$ de $\Phi (N,1,0)$  unique \`a
1-\'equivalence
int\'erieure pr\`es tel que $s(x) = u$ , $b(x) = v$  et  $F(N,1)(x)$
1-\'equivalent \`a $w$.
Comme $t(w) = t\bigl[F(N,1)(x)\bigr] = TF(N,1)(t(x))$  et que la classe
$t(x)$ est unique,
on d\'eduit que le diagramme (2) poss\'ede la propriet\'e de 0-\'equivalence.
Par cons\'equent $TF$ est une ($k$-1)-\'equivalence ext\'erieure.
\par (b) On sait d\'ej\`a d'apr\`es la proposition (1.2.5) que $T$ est un
foncteur de
$n$-{\cal PN} vers ($n$-1)-{\cal PN} et en plus d'apr\`es (a) ce foncteur
respecte
l'\'equivalence ext\'erieure, donc il suffit de montrer que si $\Phi $ est un
$n$-nerf alors $T\Phi $ est un ($n$-1)-nerf. \par\hskip 5mm Soit
$(M,m,N)$ un objet de $\Delta ^{^{n-s-1}}\times\Delta \times\Delta ^{s}$ avec
$1\leq s\leq n-1$\par\vskip 2mm{\bf (i)} Le foncteur
$(T\Phi )_{_{M,0}} = T(\Phi _{_{M,0}})$ est constant, en effet :\par\vskip 2mm
$(T\Phi )_{_{M,0}}(N) = t\bigl [\Phi _{_{M,0}}(N,0)\bigr] = \Phi
(M,0_{_{s+1}})$.
\par\vskip 2mm{\bf (ii)} D'apr\`es {\bf (a)} $(T\delta )^{^{[m]}}_{_{M}} :
(T\Phi )_{_{M,m}}
\fhd{}{}(T\Phi )_{_{M,1}}{\times}_{(T\Phi
)_{_{M,0}}}\dots.{\times}_{_{(T\Phi )_{_{M,0}}}}
(T\Phi )_{_{M,1}}$ est une ($s$-1)-\'equivalence ext\'erieure, car
$\delta ^{^{[m]}}_{_{M}} : \Phi _{_{M,m}}\fhd{}{}\Phi
_{_{M,1}}{\times}_{\Phi _{_{M,0}}}\dots.
{\times}_{\Phi _{_{M,0}}}\Phi _{_{M,1}}$ est une $s$-\'equivalence
ext\'erieure et
$T(\Phi _{_{M,m}}) = (T\Phi )_{_{M,m}}$. On en d\'eduit que $T\Phi $ est un
($n$-1)-nerf et par cons\'equent $T$ est bien un foncteur de $n$-Cat vers
($n$-1)-Cat.\hfill


\par\vskip 5mm\hskip 5mm{\bf Lemme (1.3.6) :} {\it Le compos\'e de deux
$n$-\'equivalences ext\'erieures est une $n$-\'equivalence ext\'erieure.}
\par\vskip 5mm\hskip 5mm{\bf Preuve :} Soient $\Phi \fhd{F}{}\Psi
\fhd{G}{}\Theta$
deux $n$-\'equivalences ext\'erieures, on sait que le compos\'e $G.F$ est le
morphisme d\'efinie pour tout objet $M$ de $\Delta ^{n}$ par
$G.F(M) = G(M).F(M)$. Soit $h$ un \'el\'ement de l'ensemble $\{0,\dots,n\}$
et posons
$N = I_{_{n-h-1}}$. Les diagrammes (3) et (4) suivants :\par
\centerline{$\diagram{\Phi (N,0_{_{h+1}})&
\Fhd{F(N,0_{_{h+1}})}{}&\Psi (N,0_{_{h+1}})&\Fhd{G(N,0_{_{h+1}})}{}&
\Theta (N,0_{_{h+1}})\cr \fvb{s,b}{}&(3)&\fvb{s,b}{}&(4)&\fvb{s,b}{}\cr
\Phi (N,0_{_{h}})&\Fhd{}{F(N,1,0_{_{h}})}&\Psi
(N,0_{_{h}})&\Fhd{}{G(N,1,0_{_{h}})}&
\Theta (N,0_{_{h}})}$}\par
poss\'edent la propriet\'e de $h$-\'equivalence, donc pour tout
$w\in\Theta (N,1,0_{_{h}})$  et  tout $u,v\in\Phi (N,0_{_{h+1}})$  tels que
$s(w) =  G.F(N,0_{_{h+1}})(u)$  et  $b(w) = G.F(N,0_{_{h+1}})(v)$ on a
:\par\vskip 2mm
{\bf (i)} $\exists x^{'}\in \Psi (N,1,0_{_{h}})$ tel que
$s(x^{'})=F(N,0_{_{h+1}})(u)$,
$b(x^{'})=F(N,0_{_{h+1}})(v)$ \par\vskip 2mm\hskip 5mm
et  $t^{h}(w)=t^{h}\bigl[G(N,1,0_{_{h}})(x^{'})\bigr]=
(T^{h}G)(N,1,0_{_{h}})\bigl[t^{h}(x^{'})\bigr]$
\hskip 1cm (diagramme (4))\par\vskip 2mm
{\bf (ii)}$\exists x\in \Phi (N,1,0_{_{h}})$ tel que $s(x)=u$,  $b(x)=v$
\par\vskip 2mm\hskip 5mm
et
$t^{h}(x^{'})=t^{h}\bigl[F(N,1,0_{_{h}})(x)\bigr]=(T^{h}F)(N,1,0_{_{h}})
\bigl[t^{h}(x)\bigr]$
\hskip 1cm (diagramme (3))\par\vskip 2mm
En appliquant $T^{h}G(N,1)$ \`a la deuxi\`eme formule on obtient :\par\vskip
3mm
 \centerline{$\matrix{T^{h}G(N,1)\bigl[t^{h}(x)\bigr]&=&
T^{h}G\bigl[T^{h}F(N,1)(t^{h}(x))\bigr]\cr\cr
t^{h}(w)&=&T^{h}(F.G)\bigl[t^{h}(x)\bigr]\cr\cr
t^{h}(w)&=&T^{h}\bigl[(G.F)(N,1)(x)\bigr]}$}\par\vskip 3mm
l'unicit\'e de $x$ est \`a $h$-\'equivalence pr\`es d'o\`u $G.F$ est une
$h$-\'equivalence.
\hfill


\par\vskip 5mm\hskip 5mm{\bf Remarque (1.3.7) :} Pour tout $h\in
\{1,\dots,n-1\}$, en composant $h$ fois le foncteur $T$ on obtient $T^{h}$
qui sera un
foncteur covariant de $n$-Cat vers ($n$-$h$)-Cat.\par
\vskip 2mm\hskip 5mm {\it Conjecture :} Soit $HO$-$n$-Cat le localis\'e
 ([5] P. Gabriel-M. Zisman) de la cat\'egorie $n$-Cat par rapport aux
$n$-\'equivalences
ext\'erieures. On peut d\'efinir un ($n$+1)-nerf $n$-CAT telle que
$HO$-$n$-Cat = $T^{n}$($n$-CAT).\par\vskip 5mm\hskip 5mm


{\bf Rappel (1.3.8) :} Une {\it $n$-cat\'egorie stricte}
$\cal C$ est la donn\'ee :\par
{\bf (1)} D'un diagramme d'ensembles et d'applications
:\par\centerline{$\diagram{C_i&
\fhg{\fhd{s_i}{}}{\fhd{}{b_i}}&C_{i+1}}$\hskip 5mm avec\hskip 5mm
$\left\{\matrix{s_{i-1}
s_i&=&s_{i-1}b_i\cr b_{i-1}s_i&=&b_{i-1}b_i\cr
s_{i}e_i&=&b_{i}e_i&=&I_{C_i}}\right.$}
\par o\`u $i\in\{0,\dots,n-1\}$.\par\vskip 3mm
{\bf (2)} D'une famille d'applications : ${\star}_{_{ij}} :
C_{_{j}}{\times}_{_{C_i}}C_{_{j}}
\ \fhd{}{}\ C_j$ , pour tout $i,j$ tels que $0\leq i < j\leq n$.\par\vskip 3mm
V\'erifiant les deux axiomes suivants :\par\vskip 3mm
\hskip 2mm{\bf i)} Pour tout $i,j$ tels que  $0\leq i < j\leq n$,  ${\cal
C}_{ij} =
(C_i,C_j,{\star}_{_{ij}})$  est une cat\'egorie.\par\hskip 2mm{\bf ii)}
Pour tout
$i,j,k$ tels que $0\leq i < j < k \leq n$ les compositions
${\star}_{_{ij}}$ et ${\star}_{_{ik}}$
v\'erifient la relation de Godement.\par \vskip 2mm\hskip 5mm
Soit $\cal C$ une $n$-cat\'egorie stricte. on appelle
le nerf multiple de $\cal C$ d'ordre $n$ le $n$-pr\'e-nerf :
\par\centerline{$\diagram{{\Delta }^{n}&\Fhd{\Phi }{}&Ens\cr M&
\Fhd{}{}&{\cal N}_{_{M}}(\cal C)}$}\par o\`u ${\cal N}_{_{M}}(\cal C)$ est
d\'efini par
 : \par\vskip 3mm\centerline{${\cal N}_{_{M}}({\cal C}) =
\displaystyle\prod ^{m_1}_{_{C_0}}\dots
\displaystyle\prod ^{m_{n}}_{_{C_{n-1}}}C_n$\hskip 3mm o\`u \hskip 3mm
$M = (m_1,\dots,m_n)$.}\par\vskip 3mm\hskip 3mm
Les images des applications $d^{^{k}}_{_{i}}$ et $\varepsilon ^{^{k}}_{_{i}}$,
pour $k$ dans $\{1,\dots,n\}$ et $i$ dans $\{0,\dots,m_{_{k}}\}$, de la
m\^eme fa\c con
que dans la Remarque (1.2.9) : pour d\'efinir $d^{^{k}}_{_{i}}$ on utilise
la composition
$\star_{_{nk}}$ pour les $n$-fl\`eches de la $k$-\`eme composante, quant aux
$\varepsilon ^{^{k}}_{_{i}}$ on ins\'ere les identit\'es.\par


\par\vskip 5mm\hskip 5mm {\bf Proposition (1.3.9) :} {\it Le nerf multiple
d'une
$n$-cat\'egorie stricte est un $n$-nerf strict (c.\`a.d les $\delta
_{_{[m]}}$ sont des
isomorphismes).}\par\vskip 4mm\hskip 4mm{\bf Preuve :} Si on d\'esigne par
$C_{_{n-1}}^{\sim}$ les classes d'isomorphismes des ($n$-$1$)-fl\`eches,
alors : \par
$T{\cal C} = (C_1,\dots,C_{_{n-2}},C_{_{n-1}}^{\sim},{\star}_{_{ij}})$ est une
($n$-$1$)-cat\'egorie stricte.\par Montrons que $\Phi $ est $1$-troncable :
Soit
$(m_1,\dots,m_n)\in Ob(\Delta ^{n})$ et posons :\par\vskip 2mm
\centerline{$M = (m_1,\dots,m_{_{n-1}})$, \hskip 4mm $X = \displaystyle
\prod ^{m_1}_{_{C_0}}\dots\displaystyle\prod
^{m_{n-1}}_{_{C_{n-2}}}C_{_{n-1}}$} \par
alors l'application :\par\centerline{$\diagram{\displaystyle\prod
^{m_1}_{_{C_0}}\dots
\displaystyle\prod ^{m_{n}}_{_{C_{n-1}}}C_n&\Fhd{\delta
^{m_n}_{[M]}}{}&\displaystyle
\prod ^{m_n}_X\biggl[\displaystyle\prod ^{m_1}_{_{C_0}}\dots\displaystyle
\prod ^{m_{n-1}}_{_{C_{n-2}}}C_n\biggr]\cr (x_{_{i_1\dots
i_n}})_{_{i_1\dots i_n}}&
\Fhd{}{}&\bigl( x_{_{i_1\dots i_{n-1},1}},\dots,x_{_{i_1\dots
i_{n-1},m_n}}\bigr)}$}\par
est une bijection, donc $\Phi $ est $1$-troncable.
\par Or  $T\Phi (m_1,\dots,m_n) = \displaystyle\prod
^{m_1}_{_{C_0}}\dots\displaystyle
\prod ^{m_{n-1}}_{_{C_{n-2}}}{C_{_{n-1}}^{\sim}}$,  donc  $T(\Phi )$ n'est
autre que le
($n$-$1$)-pr\'e-nerf associ\'ee \`a $T{\cal C}$. On en d\'eduit alors que
$\Phi $ est $2$-troncable et par la suite $n$-troncable.\par\hskip 5mm
Soient maintenant $s\in \{1,\dots,n\}$ et $(M,m)$ un objet de $\Delta
^{s}\times\Delta $,
le morphisme $\delta ^{[m]}_{M}$ est un isomorphisme. En effet si
$N$ est un objet de $\Delta ^{^{n-s-2}}$ alors l'application :\par
\centerline{$\diagram{\Phi _{_{M,m}}(N)&\Fhd{\delta
^{m}_{[M]}(N)}{}&\displaystyle
\prod ^{m_n}_{\Phi _{M,0}}\Phi _{_{M,1}}(N)\cr (x_{_{i_1\dots
i_n}})_{_{i_1\dots i_n}}&
\Fhd{}{}&\bigl( (x_{_{i_1.. i_{s-1},1,i_{s+1}..i_n}}),\dots, (x_{_{i_1..
i_{s-1},m_{n},i_{s+1}..i_n}})
\bigr)}$}\par est une une bijection, d'o\`u $\delta ^{m}_{[M]}$ est une
($n$-$s$-$1$)-\'equivalence ext\'erieure et par cons\'equent $\Phi $ est un
$n$-nerf strict.\par\vskip 5mm\hskip 5mm


{\bf (1.4).--- 2-cat\'egorie usuelle et 2-nerf.}
\par\vskip 5mm\hskip 5mm{\bf Rappel ([2], [3], [4]) :} Une $2$-cat\'egorie
large ${\cal C}$
est la donn\'ee :\par {\bf (a)} D'un diagramme d'ensembles et
d'applications :\par
\centerline{$\diagram{C_0&\fhd{\fhg{s_1}{}}{\fhg{}{b_1}}&C_1
\fhd{\fhg{s_2}{}}{\fhg{}{b_2}}&C_2}$\hskip 5mm avec
\hskip 5mm $\left\{\matrix{s_1s_2&=&b_1s_2\cr b_1s_2&=&b_1b_2\cr s_1e_2&=&
b_1e_2&=&I_{C_0}}\right.$}\par
{\bf (b)} De trois applications :\par\hskip 1cm -Composition des fl\`eches
:\par
\centerline{$\diagram{C_1{\times}_{C_0} C_1&\fhd{}{}&C_1\cr
(f,g)&\fhd{}{}&gf}\hskip 2cm
\diagram{x\fhd{f}{}y\fhd{g}{}z}\hskip 3mm\Rightarrow \hskip 3mm
\diagram{x\fhd{gf}{}z}$}\par
\hskip 1cm -Composition verticale des 2-fl\`eches :\par
$\diagram{C_2{\times}_{C_1} C_2&\fhd{}{}&C_2\cr (\varphi ,\varphi ^{'}
)&\fhd{}{}&
\varphi ^{'} \cdot\varphi }\hskip 2cm \diagram{\fvb{f}{}\fhd{\varphi
}{}\fvb{}{}
\fhd{\varphi ^{'} }{}\fvb{}{h}}\hskip 5mm\Rightarrow\hskip 5mm
\diagram{\fvb{f}{}
\fhd{\varphi ^{'}\cdot\varphi }{}\fvb{}{h}}$\par
\hskip 1cm -Composition horizontale des 2-fl\`eches :\par
\centerline{$\diagram{C_2{\times}_{C_0} C_2&\fhd{}{}&C_2\cr (\alpha ,\beta )&
\fhd{}{}&\beta *\alpha }\hskip 2cm \diagram{\fhd{f}{}&\fhd{g}{}\cr
\fvb{}{\alpha }&\fvb{}{\beta }\cr \fhd{}{f^{'}}&\fhd{}{g^{'}}}\hskip 3mm
\Rightarrow \hskip 3mm \diagram{\fhd{gf}{}\cr \fvb{}{\beta *\alpha }\cr
\fhd{}{g^{'}f^{'}}}$}\par
{\bf (c)} De trois isomorphismes de coh\'erences :\par
\hskip 1cm -D'associativit\'e :\par
\centerline{$\forall (f,g,h)\in C_1{\times}_{C_0}C_1{\times}_{C_0} C_1$
\hskip 5mm
on\ a \hskip 5mm $\diagram{(hg)f&\Fhd{\Phi (f,g,h)}{\sim}&h(gf)}$}\par
\hskip 1cm -D'identit\'es : \hskip 1cm $\forall (x,f,y)\in
C_0{\times}_{C_0}C_1{\times}_
{C_0}C_0 $\hskip 1cm on a  \par\centerline{$
\diagram{fI_{x}&\Fhd{U(f)}{\sim}&f}$ \hskip 1cm et \hskip 1cm
$\diagram{I_{y}f&\Fhd{V(f)}{\sim}&f}$}
\par\hskip 5mm V\'erifiant la commutativit\'e des diagrammes
ferm\'es de $2$-fl\`eches suivants, faisant intervenir les donn\'ees
pr\'ec\'edentes. \par
\vskip 5mm\hskip 5mm {\bf (1)} (Axiome du Pentagone)
\hskip 5mm pour \hskip 4mm $x\fhd{f}{}y\fhd{g}{}z\fhd{h}{}t\fhd{k}{}s$
\hskip 4mm on a
$$\diagram{[(kh)g]f&\Fhd{\Phi (k,h,g)*I_{f}}{}&[k(hg)]f&\Fhd{\Phi
(k,hg,f)}{}&k[(hg)f]\cr
\Vert &&&&\fvb{}{I_{k}*\Phi (h,g,f)}\cr
[(kh)g]f&\Fhd{}{\Phi (kh,g,f)}&(kh)(gf)&\Fhd{}{\Phi (k,h,gf)}&k[h(gf)]}$$
\hskip 5mm{\bf (2)} (Associativite de la composition verticale des
2-fl\`eches) \par
\centerline{Pour \hskip 5mm$\diagram{\fvb{}{}\fhd{\alpha
}{}\fvb{}{}\fhd{\beta }{}
\fvb{}{}\fhd{\gamma }{}\fvb{}{}}$\hskip 5mm on a\hskip 5mm
$(\gamma \cdot\beta )\cdot\alpha  = \gamma \cdot (\beta \cdot\alpha)$}\par
\centerline{{\bf (3)}\hskip 5mm Pour \hskip 5mm$\diagram{\fhd{}{}&\fhd{}{}\cr
\fvb{}{\alpha }&\fvb{}{{\alpha }^{'}}\cr \fhd{}{}&\fhd{}{}\cr \fvb{}{\beta }&
\fvb{}{{\beta }^{'}}\cr \fhd{}{}&\fhd{}{}}$\hskip 1cm on a \hskip
1cm$({\beta }^{'}
\cdot{\alpha }^{'})*(\beta \cdot\alpha ) = ({\beta }^{'}*\beta )\cdot({\alpha
}^{'}*\alpha )$}
\centerline{{\bf (4)}\hskip 5mm Pour \hskip
5mm$\diagram{\fhd{f}{}&\fhd{I_{_{y}}}{}\cr
\fvb{}{\alpha }&\fvb{}{I^{^{2}}_{_{y}}}\cr
\fhd{}{g}&\fhd{}{I_{_{y}}}}$\hskip 1cm on a
\hskip 1cm$V(g)\cdot(I^{^{2}}_{_{y}}*\alpha ) = \alpha\cdot V(f)$}\par
\centerline{{\bf (5)}\hskip 5mm Pour \hskip
5mm$\diagram{\fhd{I_{_{x}}}{}&\fhd{f}{}\cr
\fvb{}{I^{^{2}}_{_{y}}}&\fvb{}{\alpha }\cr
\fhd{}{I_{_{x}}}&\fhd{}{g}}$\hskip 1cm on a
\hskip 1cm $U(g)\cdot(\alpha *I^{^{2}}_{_{y}}) = \alpha\cdot U(f)$}
\par\vskip 2mm{\bf (6)} \hskip 2cm$\left\{\matrix{V(I_{x})&=&U(I_{x})\cr &&\cr
I_{g}\star I_{f}&=&I_{gf} \hskip 5mm pour \hskip 5mm\fhd{f}{}\fhd{g}{}}\right.$
\par\vskip 2mm
Et la commutativit\'e des diagrammes suivants :\par
\centerline{{\bf (7)}\hskip 2cm$\diagram{(hg)f&\Fhd{\Phi (h,g,f)}{}&h(gf)\cr
\fvb{(\gamma *\beta )*\alpha }{}&&\fvb{}{\gamma *(\beta *\alpha )}\cr
(h^{'}g^{'})f^{'}&\Fhd{}{\Phi (h^{'},g^{'},f^{'})}&h^{'}(g^{'}f^{'})}$
\hskip 1cm avec\hskip 5mm
$\diagram{\fhd{f}{}&\fhd{g}{}&\fhd{h}{}\cr \fvb{}{\alpha }&\fvb{}{\beta }
&\fvb{}{\gamma }\cr \fhd{}{f^{'}}&\fhd{}{g^{'}}&\fhd{}{h^{'}}}$}\par
\centerline{{\bf (8)}\hskip 1cm$\diagram{(gI_{y})f&\fhd{\Phi
(g,I_{y},f)}{}&g(I_{y}f)\cr
\Vert &&\fvb{}{I_{g}*V(f)}\cr (gI_{y})f&\fhd{}{U(g)*I_{f}}&gf}$\hskip 1cm avec
\hskip 5mm $x\fhd{f}{}y\fhd{I_y}{}y\fhd{g}{}z$}\par
\centerline{{\bf (9)}\hskip 1cm
$\diagram{(gf)I_{x}&\fhd{\Phi (g,f,I_{x})}{}&g(fI_{x})\cr
\Vert &&\fvb{}{I_{g}*U(f)}\cr (gf)I_{x}&\fhd{}{U(gf)}&gf}$\hskip 1cm avec
\hskip 5mm
$x\fhd{I_{x}}{}x\fhd{f}{}y\fhd{g}{}z$}\par
\centerline{{\bf (10)}\hskip 1cm$\diagram{(I_{z}g)f&\fhd{\Phi
(I_{z},g,f)}{}&I_{y}(gf)\cr
\Vert &&\fvb{}{V(gf)}\cr (I_{z}g)f&\fhd{}{V(g)*I_{f}}&gf}$\hskip 1cm avec
\hskip 5mm$x\fhd{f}{}y\fhd{g}{}z\fhd{I_{z}}{}z$}\par\hskip 5mm
\par\hskip 5mm
Soit ${\cal C}$ une 2-cat\'egorie large. On appelle {\it nerf double} de
${\cal C}$ le
2-pr\'e-nerf $\Phi $ d\'efini pour tout objet $(m,n)$ de $\Delta  ^{^{2}}$
par  :\par
$\Phi (m,n)$ est l'ensemble des familles de quadruplets
$(x_{_{i}}, f_{_{ij}}^{^{\alpha }}, \lambda _{_{ij}}^{^{\alpha \beta }},
\varepsilon_{_{ijk}}^{^{\alpha }})$ \'el\'ements du produit \par\vskip 2mm
${{\cal C}_{_{0}}\times {\cal C}_{_{1}}\times {\cal C}_{_{2}}\times {\cal
C}_{_{2}}}$,
\hskip 3mm avec \hskip 3mm$0\leq i < j < k \leq m$\hskip 3mm ,\hskip 3mm
$ 0\leq \alpha  < \beta < \gamma \leq n$  \hskip 3mm et \par\vskip 3mm
$$x_{_{i}}\Fhd{f_{_{ij}}^{^{\alpha }}}{}x_{_{j}}\hskip 1cm
f_{_{ij}}^{^{\alpha }}
\Fhd{\lambda _{_{ij}}^{^{\alpha \beta }}}{}f_{_{ij}}^{^{\beta }}\hskip 1cm
f_{_{jk}}^{^{\alpha }}f_{_{ij}}^{^{\alpha }}\Fhd
{\varepsilon_{_{ijk}}^{^{\alpha }}}{\sim}f_{_{ik}}^{^{\alpha }}$$\par\vskip 3mm
tels qu'on a les relations de coh\`erences suivantes :\par\vskip 4mm
\centerline{$\matrix{\lambda _{_{ik}}^{^{\alpha \beta }}\cdot
\varepsilon_{_{ijk}}^{^{\alpha }}&=&\varepsilon_{_{ijk}}^{^{\beta }}\cdot
(\lambda _{_{ik}}^{^{\alpha \beta }}\star\lambda _{_{ij}}^{^{\alpha \beta }})
\hskip 2cm (i)\cr\lambda _{_{ij}}^{^{\beta \gamma }}\cdot
\lambda _{_{ij}}^{^{\alpha \beta }}&=&\lambda _{_{ij}}^{^{\alpha \gamma }}
\hskip 4cm (ii)}$}\par\vskip 5mm\hskip 5mm
L'image par $\Phi $ de la fl\`eche \'el\'ementaire $d^{^{1}}_{_{a}}$ (
respectivement
$d^{^{2}}_{_{a}}$ ) est l'application qui \`a la famille $\{(x_{_{i}},
f_{_{ij}}^{^{\alpha }},
\lambda _{_{ij}}^{^{\alpha \beta }}, \varepsilon_{_{ijk}}^{^{\alpha }})\}$
fait correspondre
la famille $\{(x_{_{i}}, f_{_{ij}}^{^{\alpha }}, \lambda _{_{ij}}^{^{\alpha
\beta }},
\varepsilon_{_{ijk}}^{^{\alpha }})\}$ tels que $i, j,
k\in\{0,\dots,m\}\setminus\{a\}$ (
respectivement $\alpha , \beta \in\{0,\dots,n\}\setminus\{a\}$ ). Et pour les
$\varepsilon ^{^{k}}_{_{a}}$ l'action de $\Phi (\varepsilon
^{^{k}}_{_{a}})$ sur un quadruplet
est l'insertion de l'identit\'e entre la position d'ordre $a$ et $a$+1.


\par\vskip 5mm\hskip 5mm{\bf Th\'eor\`eme (1.4.1) :} {\it Le nerf double
d'une 2-cat\'egorie large est un 2-nerf.}\par\vskip 5mm\hskip 5mm {\bf Preuve
:}
(a) Soit $\Phi $ le nerf double d'une 2-cat\'egorie large ${\cal C}$, alors
il est claire qu'on a :
$\Phi (0,n) =\Phi (0,0) ={{\cal C}_{_{0}}}$.\par
(b) Montrons d'abord que $\Phi $ est 2-troncable. Pour cela on va
v\'erifier que
$\Phi $ et $T\Phi $ sont 1-troncables. Soit $(m,n)$ un couple d'entiers
positifs non nuls, et
posant :\par\vskip 3mm\centerline{$\matrix{\delta ^{^{[n]}}_{_{m}}&=&
\Phi (\delta ^{^{01}}_{_{m}})\times\dots\times\Phi (\delta ^{^{n-1,n}}_{_{m}})
\cr\cr\delta ^{^{n}}_{_{[m]}}&=&\Phi (\delta ^{^{n}}_{_{01}})\times\dots\times
\Phi (\delta ^{^{n}}_{_{m-1,m}})}$\hskip 1cm avec}\par\vskip 3mm
\centerline{$\diagram{[m]\times[1]&\fhd
{\delta ^{^{\alpha ,\beta }}_{_{m}}}{}&[m]\times[n]\cr
(x,o)&\fhd{}{}&(x,\alpha )\cr (x,1)&
\fhd{}{}&(x,\beta )}$\hskip 1cm$\diagram{[1]\times[n]&\fhd{\delta
^{^{n}}_{_{ij}}}{}&
[m]\times[n]\cr (0,x)&\fhd{}{}&(i,x)\cr (1,x)&\fhd{}{}&(j,x)}$}\par\hskip 5mm
L'application suivante :\par
\centerline{$\diagram{\Phi _{_{m,n}}&\fhd{\delta ^{^{[n]}}_{_{m}}}{}&\Phi
_{_{m,1}}
{\times}_{_{\Phi _{_{m,0}}}}\dots{\times}_{_{\Phi _{_{m,0}}}}\Phi _{_{m,1}}\cr
(x_{_{i}}, f_{_{ij}}^{^{\alpha }}, \lambda _{_{ij}}^{^{\alpha \beta }},
\varepsilon_{_{ijk}}^{^{\alpha }})&\fhd{}{}&\Bigr((x_{_{i}},
f_{_{ij}}^{^{\alpha }},
\lambda _{_{ij}}^{^{01}}, \varepsilon_{_{ijk}}^{^{\alpha }})_{_{\alpha =0,1}},
\dots,(x_{_{i}}, f_{_{ij}}^{^{\alpha }}, \lambda _{_{ij}}^{^{n-1,n}},
\varepsilon_{_{ijk}}^{^{\alpha }})_{_{\alpha =n-1,n}}\Bigr)}$}
 est bien une bijection car si on pose $\lambda _{_{ij}}^{^{\alpha \beta }} =
\lambda _{_{ij}}^{^{\beta -1 ,\beta }}\dots\lambda _{_{ij}}^{^{\alpha
,\alpha +1}}$
on construit de fa\c con unique l'ant\'ec\'edent par $\delta
^{^{[n]}}_{_{m}}$ d'un
\'el\'ement de l'ensemble $\Phi _{_{m,1}}{\times}_{_{\Phi _{_{m,0}}}}
\dots{\times}_{_{\Phi _{_{m,0}}}}\Phi _{_{m,1}}$.\par Ce qui entraine que
$\Phi $ est 1-troncable.\par\hskip 5mm
Posons $\Psi =T\Phi $ , alors $\Psi $ est 1-troncable si et seulement si
c'est une cat\'egorie. D'abord on sait que : \par\vskip 3mm
\centerline{$\matrix{\Psi _{_{1}}{\times}_{_{\Psi _{_{0}}}},\dots,
{\times}_{_{\Psi _{_{0}}}}\Psi _{_{1}}&=&(T\Phi )_{_{1}}{\times}_{(T\Phi
)_{_{0}}}
\dots{\times}_{(T\Phi )_{_{0}}}
(T\Phi )_{_{1}}\cr\cr &=&[\Phi _{_{0,1}}]^{\sim}{\times}_{_{[\Phi
_{_{0,0}}]^{\sim}}}
\dots{\times}_{_{[\Phi _{_{0,0}}]^{\sim}}}[\Phi _{_{1,0}}]^{\sim}\cr\cr
&=&[\Phi _{_{1,0}}
{\times}_{_{\Phi _{_{0,0}}}}\dots{\times}_{_{\Phi _{_{0,0}}}}\Phi
_{_{1,0}}]^{\sim}}$}
\par\vskip 3mm alors l'application suivante :\par
\centerline{$\diagram{\Psi _{_{m}}&\Fhd{\delta _{_{[m]}}}{}&\Psi _{_{1}}
{\times}_{_{\Psi _{_{0}}}},\dots,{\times}_{_{\Psi _{_{0}}}}\Psi _{_{1}}
\cr t(f_{_{ij}},\varepsilon
_{_{ijk}})&\Fhd{}{}&t\Bigr(f_{_{01}},\dots,f_{_{m-1,m}}\Bigr)}$}
\par est de m\^eme une bijection en effet : soit
$Y =t[(x_{_{i}},f_{_{i,i+1}})_{_{0\leq i\leq m}}]$ un \'el\'ement de
$\Psi _{_{1}}{\times}_{_{\Psi _{_{0}}}},\dots,{\times}_{_{\Psi
_{_{0}}}}\Psi _{_{1}}$,
et posons :\par\vskip 2mm\centerline{$f_{_{ij}} =
\Bigr((f_{_{j-1,j}}f_{_{j-2,j-1}})
\dots\Bigr)f_{_{i,i+1}}$ \hskip 5mm pour tout \hskip 5mm$0\leq i < j\leq m$}
\par\vskip 2mm on peut alors construire gr\^ace aux isomorphismes
d'associativives
des isomorphismes de coh\`erences
$\varepsilon _{_{ijk}} : f_{_{jk}}f_{_{ij}}\fhd{}{}f_{_{ik}}$ pour tout
$i,j,k$ tels que
$0\leq i < j < k\leq m$. Alors $t(x_{_{i}},f_{_{i,j}},\varepsilon
_{_{ijk}})$ est un
l'ant\'ec\'edent de $Y$ par $\delta _{_{[m]}}$, d'o\`u la surjection. Soient
$t(x_{_{i}},f_{_{i,j}},\varepsilon _{_{ijk}})$ et
$t(x_{_{i}},g_{_{i,j}},\mu _{_{ijk}})$
deux \'el\'ement de $\Psi _{_{m}}$ qui ont m\^eme image par $\delta
_{_{[m]}}$, donc il
existe une famille d'isomorphismes $\lambda _{_{i,i+1}}:
f_{_{i,i+1}}\fhd{}{}g_{_{i,i+1}}$
pour $0\leq i\leq m-1$. Gr\^ace aux compositions horizontale et verticale
des 2-fl\`eches
on peut construire une famille d'isomorphismes $\lambda _{_{i,j}}: f_{_{i,j}}
\fhd{}{}g_{_{i,j}}$ de fa\c con recurrente par :\par\vskip 3mm
\centerline{$\matrix{\lambda _{_{i,i+2}}&=&\mu _{_{i,i+1,i+2}}\cdot(\lambda
_{_{i+1,i+2}}
\star\lambda _{_{i,i+1}})\cdot(\varepsilon _{_{i,i+1,i+2}})^{^{-1}}\cr\cr
\lambda _{_{i,i+3}} &=& \mu _{_{i,i+2,i+3}}\cdot(\lambda _{_{i+2,i+3}}\star
\lambda _{_{i,i+2}})\cdot(\varepsilon
_{_{i,i+2,i+3}})^{^{-1}}\cr\dots&&\dots\cr
\lambda _{_{i,j}} &=& \mu _{_{i,j-1,j}}\cdot(\lambda
_{_{j-1,j}}\star\lambda _{_{i,j-1}})
\cdot(\varepsilon _{_{i,j-1,j}})^{^{-1}}}$}\par \vskip 3mm
et qui v\'erifie la relation $(i)$ par construction. Ce qui montre que
$(x_{_{i}},f_{_{i,j}},\varepsilon _{_{ijk}})$ et $(x_{_{i}},g_{_{i,j}},\mu
_{_{ijk}})$
repr\'esentent la m\^eme classe, et par cons\'equent $\delta _{_{[m]}}$ est
injective.
\par\vskip 5mm
(c) V\'erifions maintenent que $\delta _{_{[m]}}:\Phi _{_{m}}\fhd{}{}\Phi
_{_{1}}
{\times}_{\Phi _{_{0}}}\dots{\times}_{\Phi _{_{0}}}\Phi _{_{1}}$ est
une 1-\'equivalence ext\'erieure. Il y a deux cas \`a v\'erifier
:\par\hskip 4mm
{\bf (h=0) :} On montre que le diagramme suivant poss\'ede la propriet\'e de
0-\'equivalence :\par\centerline{$\diagram{\Phi _{_{m,0}}&\Fhd{\delta
^{^{0}}_{_{[m]}}}{}&
\Phi _{_{1,0}}{\times}_{_{C_0}}\dots{\times}_{_{C_0}}\Phi _{_{1,0}}\cr
\fvh{s,b}{}&&\fvh{}{s,b}\cr \Phi _{_{m,1}}&\Fhd{\delta ^{^{1}}_{_{[m]}}}{}&
\Phi _{_{1,1}}{\times}_{_{C_0}}\dots{\times}_{_{C_0}}\Phi _{_{1,1}}}$}\par
\centerline{Soient $\left\{\matrix{u, v &\in &\Phi _{_{m,0}}\cr w &\in &
\Phi _{_{1,1}}{\times}_{_{C_0}}\dots{\times}_{_{C_0}}\Phi _{_{1,1}}}\right.$
\hskip 3mm tels que \hskip 3mm $\left\{\matrix{\delta
^{^{0}}_{_{[m]}}(u)&=&s(w)\cr
\delta ^{^{0}}_{_{[m]}}(v)&=&b(w)}\right.$}
\par\vskip 3mm Donc $u$ ,$v$ et $w$ sont de la forme :\par\vskip 3mm
\centerline{$w =
\Bigr(x_{_{i}},f^{^{0}}_{_{i,i+1}},f^{^{1}}_{_{i,i+1}},\lambda _{_{i,i+1}}
\Bigr)_{_{0\leq i\leq m-1}}$\hskip 1cm et \hskip 1cm $\left\{\matrix
{u&=&{(x_{_{i}},f^{^{0}}_{_{i,j}},\varepsilon ^{^{0}}_{_{ijk}})}_{_{0\leq i
<j < k\leq m}}
\cr v&=&{(x_{_{i}},f^{^{1}}_{_{i,j}},\varepsilon
^{^{1}}_{_{ijk}})}_{_{0\leq i <j < k\leq m}}}
\right.$}\par\vskip 3mm Une construction analogue \`a celle de (a) permet
de construire
des 2-fl\`eches $\lambda _{_{ij}}:f^{^{0}}_{_{ij}}\fhd{}{}f^{^{1}}_{_{ij}}$
qui v\'erifient la
relation $(i)$, alors $X ={(x_{_{i}},f^{^{\alpha }}_{_{ij}},\lambda _{_{ij}},
\varepsilon ^{^{\alpha }}_{_{ijk}})}_{_{\alpha = 0, 1}}$ est l'unique
\'el\'ement de
$\Phi _{_{m,1}}$ tel que $\delta ^{^{1}}_{_{[m]}}(X) = w $ ,  $s(X) = u $
et $b(X) = v$.
\par\vskip 4mm{\bf (h=1) :} Il s'agit dans ce cas de montrer que
l'appliction suivante :\par
\centerline{$\diagram{\Phi _{_{m,0}}&\Fhd{\delta ^{^{0}}_{_{[m]}}}{}&\Phi
_{_{1,0}}
{\times}_{_{C_0}}\dots{\times}_{_{C_0}}\Phi _{_{1,0}}\cr (x_{_{i}},f_{_{ij}},
\varepsilon _{_{ijk}})&\Fhd{}{}&\Bigr(f_{_{01}},\dots,f_{_{m-1,m}}\Bigr)}$}
poss\'ede la propriet\'e de 1-\'equivalence, qui est une cons\'equence du
fait que
$T\Phi $ est une cat\'egorie. On conclu alors que $\delta _{_{[m]}}$ est
une 1-\'equivalence
ext\'erieure et par suite $\Phi $ est un 2-nerf. \hfill


\par\vskip 5mm\hskip 5mm{\bf Th\'eor\`eme (1.4.2) :}
{\it A chaque 2-nerf correspond une 2-cat\'egorie large .}
\par\vskip 5mm\hskip 5mm{\bf Preuve :} Soit $\Phi $ un 2-nerf, on va
lui associer une une 2-cat\'egorie large $\cal C$ tel que ${\cal
C}_{_{0}}=\Phi _{_{0,0}}$ ,
${\cal C}_{_{1}}=\Phi _{_{1,0}}$ et ${\cal C}_{_{2}}=\Phi _{_{1,1}}$. Les
applications source,
but et identit\'e sont obtenues en prenant l'image par $\Phi $ du diagramme de
$\Delta ^{^{2}}$ suivant : \par
\centerline{$\diagram{[0]\times[0]&\Fhg{\Fhd{\delta _{0}\times
I}{}}{\Fhd{}{\delta _{1}
\times I}}&[1]\times[0]&\Fhg{\Fhd{I\times\delta
_{0}}{}}{\Fhd{}{I\times\delta _{1}}}&[1]
\times[1]}$}\par qui donne le diagramme :\par
\centerline{$\diagram{C_0&\Fhd{\Fhg{s_1}{}}{\Fhg{}{b_1}}&C_1&\Fhd{\Fhg{s_2}{}}
{\Fhg{}{b_2}}&C_2}$\hskip 4mm tels que \hskip 4mm
$\left\{\matrix{s_{1}s_2&=&b_{1}s_2\cr b_{1}s_2&=&b_{1}b_2\cr s_{i}e_i&=&
b_{i}e_i&=&I_{C_i}}\right.$}


\par\vskip 5mm\hskip 5mm{\bf Lemme (1.4.3) :} {\it (i) Pour toute fl\`eche
$\theta :[m]\fhd{}{}[n]$ de $\Delta $ on a un foncteur $\Theta $ de $\Phi
_{_{n}}$ vers
$\Phi _{_{m}}$ d\'efini par :\par\centerline{$\matrix{\sigma \in\Phi
_{_{n,0}}&\Fhd{}{}&
\theta ^{^{0}}(\sigma )\in \Phi _{_{m,0}}\cr\varepsilon :\sigma \fhd{}{}\tau
\in
\Phi _{_{n,1}}&\Fhd{}{}&\theta ^{^{1}}(\varepsilon ) : \theta
^{^{0}}(\sigma )\fhd{}{}
\theta ^{^{0}}(\tau )\in\Phi _{_{m,1}}\cr\cr avec\hskip 4mm \theta
^{^{0}}(\sigma ):=
\Phi (\theta \times I_{_{[0]}})(\sigma ) & et &\theta ^{^{1}}(\varepsilon
):=\Phi (\theta
\times I_{_{[1]}})(\varepsilon )}$}\par\vskip 2mm
(ii) Si $\varepsilon :\sigma \fhd{}{}\tau $ est un isomorphisme alors
pour tout $\theta $ dans fl$\Delta $ , $\theta ^{^{1}}(\varepsilon )$ est
un isomorphisme.
\par (iii) Soit $\varepsilon $ dans $\Phi _{_{m,1}}$, si $\delta
^{^{1}}_{_{[m]}}(\varepsilon )$
est un isomorphisme alors $\varepsilon $ est un isomorphisme.}
\par\vskip 5mm\hskip 5mm
{\bf Preuve :} (i) la commutativit\'e  du diagramme suivant, obtenu comme
image par
$\Phi $ d'un diagramme commutatif dans $\Delta ^{^{2}}$, \par
\centerline{$\diagram{\Phi _{_{n,0}}&\Fhd{e^{^{0}}_{_{n}}}{}&
\Phi _{_{n,1}}&\Fhd{s,b}{}&\Phi _{_{n,0}}\cr
\fvb{\theta ^{^{0}}}{}&(a)&\fvb{}{\theta ^{^{1}}}&(b)&
\fvh{}{\theta ^{^{0}}}\cr \Phi _{_{m,0}}&\Fhd{}{e^{^{0}}_{_{m}}}&\Phi
_{_{m,1}}&
\Fhd{}{s,b}&\Phi _{_{1,0}}}$\hskip 1cm o\`u \hskip 1cm
$\left\{\matrix{e^{^{0}}_{_{n}}&=&\Phi (I_{_{[n]}}\times e_{_{0}})\cr
e^{^{0}}_{_{m}}&=&\Phi (I_{_{[m]}}\times e_{_{0}})}\right.$}\par
(a) montre que $\theta ^{^{1}}(I_{_{\sigma }})=I_{_{\theta ^{^{0}}(\sigma
)}}$ et
(b) montre que $\theta ^{^{1}}(\varepsilon ) : \theta ^{^{0}}(\sigma )\fhd{}{}
\theta ^{^{0}}(\tau )$.\par\vskip 2mm\hskip 2mm
La fonctorialit\'e de $\Theta $ est due \`a la commutativit\'e du diagramme
suivant :\par
\centerline{$\diagram{\Phi _{_{n,2}}&\Fhd{\delta ^{^{[2]}}_{_{n}}}{\sim}&
\Phi _{_{n,1}}{\times}_{_{\Phi _{_{n,0}}}}\Phi _{_{n,1}}&\Fhd{\theta
^{^{1}}\times
\theta ^{^{1}}}{}&\Phi _{_{m,1}}{\times}_{_{\Phi _{_{m,0}}}}\Phi _{_{m,1}}\cr
\fvb{\delta ^{^{02}}_{_{n}}}{}&&(\star)&&\fvh{\wr}{\delta ^{^{[2]}}_{_{m}}}\cr
\Phi _{_{n,1}}&\Fhd{}{\theta ^{^{1}}}&\Phi _{_{m,1}}&\Fhd{}{\delta
^{^{02}}_{_{m}}}&
\Phi _{_{m,2}}}$}\par L'application
$\Phi (\theta \times\delta _{_{02}}) : \Phi _{_{n,2}}\fhd{}{}\Phi
_{_{m,2}}$ coupe le
diagramme pr\'ec\'edent en deux diagrammes dont la commutativit\'e est
cons\'equence de l'image par $\Phi $ du diagramme commutatif de $\Delta
^{^{2}}$ suivant :
\par \centerline{$\diagram{[n]\times [2]&\Fhg{\theta \times I_{_{[2]}}}{}&
[m]\times [2]\cr\fvh{I_{_{[n]}}\times\delta
_{_{02}}}{}&&\fvh{}{I_{_{[m]}}\times
\delta _{_{02}}}\cr [n]\times [1]&\Fhg{}{\theta \times
I_{_{[1]}}}&[m]\times [1]}$}
\par D'autre part on sait que $\delta ^{^{[2]}}_{_{n}}$ et $\delta
^{^{[2]}}_{_{m}}$ sont des
bijections, ce qui entraine la commutativit\'e du diagramme ($\star$). Par
cons\'equent pour tout $(\varepsilon ,\varepsilon ^{'})$ dans
$\Phi _{_{m,1}}{\times}_{_{\Phi _{_{m,0}}}}\Phi _{_{m,1}}$ on a
$\theta ^{^{1}}(\varepsilon ^{'}\cdot\varepsilon )=\theta ^{^{1}}
(\varepsilon ^{'})\cdot\theta ^{^{1}}(\varepsilon )$.\par\vskip 2mm	\hskip 2mm
(ii) est une cons\'equence immediate de (i).	\par\vskip 2mm\hskip 2mm
(iii) Soit $\lambda :\delta ^{^{0}}_{_{[m]}}(\tau )\fhd{}{}\delta
^{^{0}}_{_{[m]}}(\sigma )$
l'inverse de $\delta ^{^{1}}_{_{[m]}}(\varepsilon )$. Comme $\Phi $ est une
2-cat\'egorie
il existe un unique $\varepsilon ^{'}:\tau \fhd{}{}\sigma $ tel que
$\delta ^{^{1}}_{_{[m]}}(\varepsilon ^{'})=\lambda $. On v\'erifie
facilement d'apr\`es (i)
que $\varepsilon ^{'}$ est l'inverse de $\varepsilon $.\hfill
\par\vskip 3mm	\hskip 5mm
{\bf (a)} {\bf Composition verticale des $2$-fl\`eches :} D'apr\`es la
proposition (1.1.4),
 pour tout entier naturel non nul $m$,  $\Phi _m$ est un 1-nerf qui
correspond \`a la cat\'egorie $(\Phi _{_{m,0}}, \Phi _{_{m,1}}, {\mu}_m)$
o\`u ${\mu}_m$
est l'unique application qui rend commutatif le diagramme suivant :
\par\centerline{
$\diagram{\Phi _{_{m,2}}&\Fhd{{\delta }^{[2]}_{m}}{\sim}&\Phi _{_{m,1}}
{\times}_{_{\Phi _{_{m,0}}}}\Phi _{_{m,1}}\cr\fvb{{\delta }^{'}_{02}}{}&&
\fvb{}{{\mu}_{_{m}}}\cr \Phi _{_{m,1}}&\Fhd{}{I_{\Phi _{_{m,1}}}}&\Phi
_{_{m,1}}}$
\hskip 1cm $\matrix{{\mu}_{_{m}}({\delta }^{[2]}_{m})&=&{\delta
}^{'}_{02}}$}\par
On d\'esigne par $L^{^{[2]}}_{_{m}}$ l'inverse de ${\delta }^{[2]}_{m}$, et
pour tout
$(\alpha ,\beta )$ dans $\Phi _{_{m,1}}{\times}_{_{\Phi _{_{m,0}}}}\Phi
_{_{m,1}}$ on pose
$\beta \cdot\alpha =\delta ^{'}_{02}L^{^{[2]}}_{_{m}}(\alpha ,\beta )$.
Donc pour $m = 1$
on retrouve la composition verticale des $2$-fl\`eches de $\Phi $, qui est
associative car
$\Phi _{_{m}}$ est une cat\'egorie ce qui entraine l'Axiome (2) de la
d\'efinition d'une
2-cat\'egorie large.\par\vskip 5mm\hskip 2mm
{\bf (b)} {\bf Composition des fl\`eches :} L'\'equivalence ext\'erieure
$\delta _{_{[m]}}$ permet de choisir \`a l'aide de l'Axiom de choix une
application  $L_{_{m}} : C_{1}{\times}_{C_0}\dots{\times}_{C_0}C_1 \fhd{}{}
\Phi _{_{m,0}}$
avec pour tout \'el\'ement $X$ de
$C_{1}{\times}_{C_0}\dots{\times}_{C_0}C_1$ un isomorphisme :\par\vskip 2mm
\centerline{$a_{_{m}}(X) : \delta ^{0}_{[m]}(L_{_{m}}(X))
\fhd{}{}X$}	\par	\vskip 2mm
Lorsque $m=2$ on obtient la composition
des fl\`eches par $gf:=\delta ^{^{0}}_{_{02}}(L_{_{2}}(f,g))$.\par\vskip
2mm\hskip 5mm
{\bf (c)} {\bf Composition horizontale des $2$-fl\`eches :} Consid\'erons
le diagramme :
\par\centerline{$\diagram{\Phi _{2}(0)&\Fhd{\delta
^{0}_{[2]}}{}&C_{1}{\times}_{C_0}C_1
\cr \fvh{s,b}{}&&\fvh{}{s,b}\cr \Phi _{2}(1)&\Fhd{}{\delta
^{1}_{[2]}}&C_{2}{\times}_{C_0}
C_2}$}\par Pour tout couple $(\alpha ,\beta )$ dans
$C_{2}{\times}_{C_0}C_2$, on pose
$s(\alpha ,\beta )=(f,g)$ , $b(\alpha ,\beta )=(f^{'},g^{'})$ , $\sigma
=L_{_{2}}(f,g)$ et
$\sigma ^{'}=L_{_{2}}(f^{'},g^{'})$. Soit $\lambda $ la compos\'ee dans la
cat\'egorie
$\Phi _{1}{\times}_{_{\Phi _0}}\Phi _{1}$ des trois fl\`eches suivantes : \par
\centerline{$\diagram{\delta ^{0}_{[2]}(\sigma
)&\Fhd{a_{_{2}}(f,g)}{\sim}&(f,g)&
\Fhd{(\alpha ,\beta )}{}&(f^{'},g^{'})&\Fhd{a_{_{2}}(f^{'},g^{'})}{\sim}&
\delta ^{0}_{[2]}(\sigma ^{'})}$}\par
Alors il existe un unique $\varepsilon _{_{\alpha ,\beta }}: \sigma
\fhd{}{} \sigma ^{'}$
dans $\Phi _{_{2,1}}$ tel que $\delta ^{1}_{[2]}(\varepsilon _{_{\alpha
,\beta }}) =\lambda $.
D'apr\`es le Lemme (1.4.4) on a le morphisme
$\delta ^{^{1}}_{_{02}}(\varepsilon _{_{\alpha ,\beta }}) :
\delta ^{^{0}}_{_{02}}(\sigma )\fhd{}{} \delta ^{^{0}}_{_{02}}(\sigma
^{'})$, donc on d\'efinit
alors la composition horizontale par $\beta \star \alpha :=
\delta ^{^{1}}_{_{02}}(\varepsilon _{_{\alpha ,\beta }}) : gf \fhd{}{}
g^{'}f^{'}$
\par\vskip 2mm\hskip 5mm
{\bf (d)} {\bf Les isomorphismes de coh\'erences d'identit\'e :}
Soient les fl\`eches de $\Delta ^{2}$ suivantes :\par
\centerline{$\diagram{[1]&\Fhd{\delta _{_{ij}}}{}&
[2]&\Fhd{\delta _{_{011}},\delta _{_{001}}}{}&[1]}$}\par
o\`u $\delta _{_{011}}(0,1,2)=(0,1,1)$ et $\delta _{_{001}}(0,1,2)=(0,0,1)$,
alors en appliquant $\Phi $ on obtient les relations :
\par\vskip 2mm\centerline{$\left\{\matrix{\delta ^{^{0}}_{_{01}}
\delta ^{^{0}}_{_{011}}&=&\delta ^{^{0}}_{_{01}}&=&I_{_{\Phi _{_{1,0}}}}\cr
\delta ^{^{0}}_{_{12}}\delta ^{^{0}}_{_{011}}&=&\delta
^{^{0}}_{_{11}}&=&I_{_{b}}
\cr\delta ^{^{0}}_{_{02}}\delta ^{^{0}}_{_{011}}&=&\delta ^{^{0}}_{_{01}}&=&
I_{_{\Phi _{_{1,0}}}}}\right.$\hskip 5mm et \hskip 5mm
$\left\{\matrix{\delta ^{^{0}}_{_{01}}\delta ^{^{0}}_{_{001}}&=&
\delta ^{^{0}}_{_{01}}&=&I_{_{\Phi _{_{1,0}}}}\cr \delta ^{^{0}}_{_{12}}
\delta ^{^{0}}_{_{001}}&=&\delta ^{^{0}}_{_{00}}&=&I_{_{s}}\cr\delta
^{^{0}}_{_{02}}
\delta ^{^{0}}_{_{001}}&=&\delta ^{^{0}}_{_{01}}&=&
I_{_{\Phi _{_{1,0}}}}}\right.$}\par\vskip 2mm
Soient $f$ un \'el\'ement de $\Phi _{_{1,0}}$ , $\tau
=L_{_{2}}(f,I_{_{b(f)}})$ et
$\sigma =\delta ^{^{0}}_{_{011}}(f)$. Donc d'apr\`es les relations
pr\'ec\'edentes on a
$\delta^{^{0}} _{_{[2]}}(\sigma )=(f,I_{_{b(f)}})$ et $\delta
^{^{0}}_{_{02}}(\sigma )=f$ ,
par cons\'equent il existe un unique isomorphisme $\varepsilon _{_{f}}: \tau
\fhd{}{}\sigma
$ tel que $\delta^{^{1}} _{_{[2]}}(\varepsilon
_{_{f}})=a_{_{2}}(f,I_{_{b(f)}})$. On en
d\'eduit alors un isomorphisme $V(f):=\delta ^{^{1}}_{_{02}}(\varepsilon
_{_{f}}) :
I_{_{b(f)}}f\fhd{}{}f$. Si on cosid\'ere la fl\`eche $\delta
^{^{0}}_{_{001}}(f)$ on
obtient avec la m\^eme construction un isomorphisme
$U(f):=\delta ^{^{1}}_{_{02}}(\varepsilon ^{^{f}}) :
fI_{_{s(f)}}\fhd{}{}f$. Lorsque
$f=I_{_{x}}$ on trouve $V(I_{_{x}})=U(I_{_{x}})$, car $\delta
^{^{0}}_{_{011}}(I_{_{x}})=
\delta ^{^{0}}_{_{001}}(I_{_{x}})$.\par\vskip 2mm\hskip 5mm
{\bf (e)} {\bf L'isomorphisme de coh\'erence d'associativit\'e :}
Soit $(f,g,h)$ dans ${\cal C}_{_{1}}{\times}_{_{{\cal C}_{_{0}}}}{\cal
C}_{_{1}}
{\times}_{_{{\cal C}_{_{0}}}}{\cal C}_{_{1}}$, et posons $\sigma
_{_{1}}=L_{_{2}}(f,g)$ ,
$\tau _{_{1}}=L_{_{2}}(g,h)$ , $\sigma _{_{2}}=L_{_{2}}(gf,h)$ , $\tau
_{_{2}}=L_{_{2}}(f,hg)$
et $T=L_{_{3}}(f,g,h)$. On veut construire de fa\c con naturelle un
isomorphisme $A(f,g,h)$
de $(hg)f=\delta ^{^{0}}_{_{02}}(\tau _{_{2}})$ vers
$h(gf)=\delta ^{^{0}}_{_{02}}(\sigma _{_{2}})$. Pour cela on va montrer que
$\sigma _{_{2}}$ et $\tau _{_{2}}$ sont isomorphes \`a deux faces du
3-simplexe $T$ ayant
$\delta ^{^{0}}_{_{03}}(T)$ comme c\^ot\'e commun, ce qui permet en se
restreignons \`a
ce c\^ot\'e de realiser notre isomorphisme. Rappelons que $\delta
^{^{n}}_{_{ijk}}=\Phi (\delta _{_{ijk}}\times I_{_{[n]}})$ et posons
$T_{_{ijk}}=\delta ^{^{0}}_{_{ijk}}(T)$ pour tout $\scriptstyle i,j,k$ tels que
$\scriptstyle 0\leq i\leq j\leq k\leq 3$. On v\'erifie facilement qu'on a
$\delta ^{^{0}}_{_{[2]}}(T_{_{ijk}})=(\delta _{_{ij}}(T),\delta
_{_{jk}}(T))$ et
les isomorphismes :\par
\centerline{$\diagram{\delta ^{^{0}}_{_{[2]}}(\tau
_{_{1}})&\Fhd{a_{_{2}}(g,h)}{\sim}
&(g,h)&\Fhg{(a^{^{2}}_{_{3}},a^{^{3}}_{_{3}})}{\sim}&\delta
^{^{0}}_{_{[2]}}(T_{_{123}})\cr
&&&&\cr\delta ^{^{0}}_{_{[2]}}(\sigma _{_{1}})&\Fhd{\sim}{a_{_{2}}(f,g)}&(f,g)&
\Fhg{\sim}{(a^{^{1}}_{_{3}},a^{^{2}}_{_{3}})}&\delta
^{^{0}}_{_{[2]}}(T_{_{012}})}$}\par
Ce qui entraine l'existence unique de deux isomorphismes :\par\vskip 2mm
\centerline{$\left\{\matrix{\eta _{_{1}}:T_{_{123}}\fhd{}{}\tau _{_{1}}\cr
\mu _{_{1}}:T_{_{012}}\fhd{}{}\sigma _{_{1}}}\right.$\hskip 5mm tels que
\hskip 5mm
$\left\{\matrix{\delta ^{^{1}}_{_{[2]}}(\eta _{_{1}})&=&
\bigl[a_{_{2}}(g,h)\bigr]^{^{-1}}\cdot (a^{^{2}}_{_{3}},a^{^{3}}_{_{3}})\cr
\delta ^{^{1}}_{_{[2]}}(\mu _{_{1}})&=
&\bigl[a_{_{2}}(f,g)\bigr]^{^{-1}}\cdot
(a^{^{1}}_{_{3}},a^{^{2}}_{_{3}})}\right.$}
\par\vskip 2mm Consid\'erons maintenant les isomorphismes \par
\centerline{$\diagram{\delta ^{^{0}}_{_{[2]}}(\tau
_{_{2}})&\Fhd{a_{_{2}}(f,hg)}{\sim}
&(f,hg)&\Fhg{(a^{^{1}}_{_{3}},\delta ^{^{1}}_{_{02}}(\eta _{_{1}}))}{\sim}&
\delta ^{^{0}}_{_{[2]}}(T_{_{013}})\cr &&&&\cr\delta
^{^{0}}_{_{[2]}}(T_{_{023}})&
\Fhd{\sim}{(\delta ^{^{1}}_{_{02}}(\mu _{_{1}}),a^{^{3}}_{_{3}})}&(gf,h)&
\Fhg{\sim}{a_{_{2}}(gf,h)}&\delta ^{^{0}}_{_{[2]}}(\sigma _{_{2}})}$}\par
De m\^eme ils existent deux isomorphisme uniques :\par\vskip 2mm
\centerline{$\left\{\matrix{\eta _{_{2}}:\tau _{_{2}}\fhd{}{}T_{_{013}}\cr
\mu _{_{2}}:T_{_{023}}\fhd{}{}\sigma _{_{2}}}\right.$\hskip 5mm tels que
\hskip 5mm
$\left\{\matrix{\delta ^{^{1}}_{_{[2]}}(\eta _{_{2}})&=
&\bigl[(a^{^{1}}_{_{3}},\delta ^{^{1}}_{_{02}}(\eta
_{_{1}})\bigr]^{^{-1}}\cdot a_{_{2}}(f,hg)
\cr\delta ^{^{1}}_{_{[2]}}(\mu _{_{2}})&=&
\bigl[a_{_{2}}(gf,h)\bigr]^{^{-1}}\cdot (\delta ^{^{1}}_{_{02}}(\mu
_{_{1}}),a^{^{3}}_{_{3}})}
\right.$}\par\vskip 2mm
Finalement on retrouve $A(f,g,h)$ en composant les isomorphismes :
\par\centerline{$\diagram{\delta ^{^{0}}_{_{02}}(\tau _{_{2}})=(hg)f&
\Fhd{\delta ^{^{1}}_{_{02}}(\eta _{_{2}})}{\sim}&\delta
^{^{0}}_{_{02}}(T_{_{013}})=
\delta ^{^{0}}_{_{02}}(T_{_{023}})&\Fhd{\delta ^{^{1}}_{_{02}}(\mu
_{_{2}})}{\sim}&h(gf)=
\delta ^{^{0}}_{_{02}}(\sigma  _{_{2}})}$}\par\hskip 5mm
V\'erifions maintenant que les compositions ainsi que les isomorphismes de
coh\'erences
qu'on vient de construire satisfont les dix axiomes d'une 2-cat\'egorie large.


\par\vskip 5mm\hskip 5mm\underbar{\bf Axiome (1)} (Axiome du Pentagone) Soient
$f,g,k,h$ quatre fl\`eches de $\Phi $ telles que $b(f)=s(g)$ , $b(g)=s(h)$
et $b(h)=s(k)$.
D'abord on a besoin de fixer un certain nombre de donn\'ees, et pour cela on va
utiliser les notations suivantes :\par\vskip 4mm
\centerline{$\matrix{R=L_{_{4}}(f,g,h,k)&R^{'}=L_{_{3}}(f,g,h)&R^{"}=L_{_{3}
}(g,h,k)\cr\cr
T=L_{_{3}}(f,g,kh)&\tau _{_{1}}=L_{_{2}}(g,kh)&\sigma
_{_{1}}=L_{_{2}}(f,g)\cr\cr
\rho =L_{_{2}}(g,h)&\tau _{_{2}}=L_{_{2}}(f,(kh)g)&\sigma
_{_{2}}=L_{_{2}}(gf,kh)\cr\cr
T^{'}=L_{_{3}}(f,hg,k)&\tau ^{'}_{_{1}}=L_{_{2}}(hg,k)&\sigma^{'}
_{_{1}}=L_{_{2}}(f,hg)\cr\cr
&\tau ^{'}_{_{2}}=L_{_{2}}(f,k(hg))&\sigma ^{'}_{_{2}}=L_{_{2}}((hg)f,k)
\cr\cr T^{"}=L_{_{3}}(gf,h,k)&\tau ^{"}_{_{1}}=L_{_{2}}(h,k)&\sigma^{"}
_{_{1}}=L_{_{2}}(gf,h)
\cr\cr &&\sigma ^{"}_{_{2}}=L_{_{2}}(h(gf),k)}$}\par\vskip 4mm\hskip 5mm
Maintenant on va consid\'erer des diagrammes commutatifs dans la cat\'egorie
$\Phi _{_{1}}{\times}_{_{\Phi _{_{0}}}}\Phi _{_{1}}$ et puis on s'int\'eressera
aux diagrammes commutatifs dans la cat\'egorie $\Phi _{_{2}}$ qui
repr\'esentent ces
premiers (un diagramme dans $\Phi _{_{1}}{\times}_{_{\Phi _{_{0}}}}\Phi
_{_{1}}$
est l'image par $\delta _{_{[2]}} $ de son repr\'e\'esentant dans $\Phi
_{_{2}}$).
\par\hskip 5mm
Soient les diagrammes dans $\Phi _{_{1}}{\times}_{_{\Phi _{_{0}}}}\Phi _{_{1}}$
et leur repr\'esentant dans $\Phi _{_{2}}$ ( les nouvelles fl\'eches qui
apparaissent sur les
diagrammes sont d\'efinies de fa\c con unique tels que ces diagrammes soient
commutatifs ) :\par\centerline{$\diagram{\delta ^{^{0}}_{_{[2]}}(\sigma
_{_{1}})&
\Fhd{a_{_{2}}(f,g)}{}&(f,g)&\Fhg{(a^{^{1}}_{_{4}},a^{^{2}}_{_{4}})}{}&
\delta ^{^{0}}_{_{[2]}}(R_{_{012}})\cr \fvb{\delta ^{^{1}}_{_{[2]}}
(\mu _{_{1}})}{}&&\fvb{Id}{}&&\fvh{}{\delta ^{^{1}}_{_{[2]}}
(\lambda _{_{012}})}\cr \delta ^{^{0}}_{_{[2]}}(R^{'} _{_{012}})&
\Fhd{}{a^{^{12}}_{_{3}}(f,g,h)}&(f,g)&\Fhg{}{a^{^{12}}_{_{3}}(f,g,kh)}&
\delta ^{^{0}}_{_{[2]}}(T_{_{012}})}$}\par
\centerline{$\diagram{\sigma _{_{1}}&\Fhd{\eta _{_{1}}}{}&R_{_{012}}\cr
\fvb{\mu _{_{1}}}{}&(a)&\fvh{}{\lambda _{_{012}}}\cr
R^{'}_{_{012}}&\Fhd{}{\eta^{'} _{_{1}}}&T_{_{012}}}$}
\par\centerline{$\diagram{\delta ^{^{0}}_{_{[2]}}(\tau ^{"}_{_{1}})&
\Fhd{a_{_{2}}(h,k)}{}&(h,k)&\Fhg{(a^{^{3}}_{_{4}},a^{^{4}}_{_{4}})}{}&
\delta ^{^{0}}_{_{[2]}}(R_{_{234}})\cr \fvb{\delta ^{^{1}}_{_{[2]}}
(\mu _{_{2}})}{}&&\fvb{Id}{}&&\fvh{}{\delta ^{^{1}}_{_{[2]}}
(\lambda ^{"}_{_{123}})}\cr \delta ^{^{0}}_{_{[2]}}(R^{"} _{_{123}})&
\Fhd{}{a^{^{23}}_{_{3}}(g,h,k)}&(h,k)&\Fhg{}{a^{^{23}}_{_{3}}(gf,h,k)}&
\delta ^{^{0}}_{_{[2]}}(T^{"}_{_{123}})}$}\par
\centerline{$\diagram{\tau ^{"}_{_{1}}&\Fhd{\eta _{_{2}}}{}&R_{_{234}}\cr
\fvb{\mu _{_{2}}}{}&(b)&\fvh{}{\lambda ^{"}_{_{123}}}\cr
R^{"}_{_{123}}&\Fhd{}{\eta^{'} _{_{2}}}&T^{"}_{_{123}}}$}
Les diagrammes suivants vont tenir compte des pr\'ec\'edents :
\par\centerline{$\diagram{\delta ^{^{0}}_{_{[2]}}(\rho )&
\Fhd{a_{_{2}}(g,h)}{}&(g,h)&\Fhg{(a^{^{2}}_{_{4}},a^{^{3}}_{_{4}})}{}&
\delta ^{^{0}}_{_{[2]}}(R_{_{123}})\cr \fvb{\delta ^{^{1}}_{_{[2]}}
(\mu _{_{3}})}{}&&\fvb{Id}{}&&\fvh{}{\delta ^{^{1}}_{_{[2]}}
(\varepsilon ^{"}_{_{012}})}\cr \delta ^{^{0}}_{_{[2]}}(R^{'} _{_{123}})&
\Fhd{}{a^{^{23}}_{_{3}}(f,g,h)}&(g,h)&\Fhg{}{a^{^{12}}_{_{3}}(g,h,k)}&
\delta ^{^{0}}_{_{[2]}}(R^{"}_{_{012}})}$}\par
\centerline{$\diagram{\rho &\Fhd{\eta _{_{3}}}{}&R_{_{123}}\cr
\fvb{\mu _{_{3}}}{}&(c)&\fvh{}{\varepsilon ^{"}_{_{012}}}\cr
R^{'}_{_{123}}&\Fhd{}{\eta^{'} _{_{3}}}&R^{"}_{_{012}}}$}
\par\centerline{$\diagram{\delta ^{^{0}}_{_{[2]}}(\tau _{_{1}})&
\Fhd{a_{_{2}}(g,kh)}{}&(g,kh)&\Fhg{(a^{^{2}}_{_{4}},\delta ^{^{1}}_{_{02}}
(\eta _{_{2}})^{^{-1}})}{}&
\delta ^{^{0}}_{_{[2]}}(R_{_{124}})\cr \fvb{\delta ^{^{1}}_{_{[2]}}
(\mu _{_{4}})}{}&&\fvb{Id}{}&&\fvh{}{\delta ^{^{1}}_{_{[2]}}
(\lambda _{_{123}})}\cr \delta ^{^{0}}_{_{[2]}}(R^{"} _{_{013}})&
\Fhd{}{(a^{^{1}}_{_{3}},\delta ^{^{1}}_{_{02}}
(\mu _{_{2}})^{^{-1}})}&(g,kh)&\Fhg{}{a^{^{23}}_{_{3}}(f,g,kh)}&
\delta ^{^{0}}_{_{[2]}}(T_{_{123}})}$}\par
\centerline{$\diagram{\tau _{_{1}}&\Fhd{\eta _{_{4}}}{}&R_{_{124}}\cr
\fvb{\mu _{_{4}}}{}&(d)&\fvh{}{\lambda _{_{123}}}\cr
R^{"}_{_{013}}&\Fhd{}{\eta^{'} _{_{4}}}&T_{_{123}}}$}
\par\centerline{$\diagram{\delta ^{^{0}}_{_{[2]}}(\tau ^{'}_{_{1}})&
\Fhd{a_{_{2}}(hg,k)}{}&(hg,k)&\Fhg{(\delta ^{^{1}}_{_{02}}
(\eta _{_{3}})^{^{-1}},a^{^{4}}_{_{4}})}{}&
\delta ^{^{0}}_{_{[2]}}(R_{_{134}})\cr \fvb{\delta ^{^{1}}_{_{[2]}}
(\mu _{_{5}})}{}&&\fvb{Id}{}&&\fvh{}{\delta ^{^{1}}_{_{[2]}}
(\lambda ^{'}_{_{123}})}\cr \delta ^{^{0}}_{_{[2]}}(R^{"} _{_{023}})&
\Fhd{}{(\delta ^{^{1}}_{_{02}}(\eta ^{'}_{_{3}}\cdot\mu
_{_{3}})^{^{-1}},a^{^{3}}_{_{3}})}&
(hg,k)&\Fhg{}{a^{^{23}}_{_{3}}(f,hg,k)}&
\delta ^{^{0}}_{_{[2]}}(T^{'}_{_{123}})}$}\par
\centerline{$\diagram{\tau ^{'}_{_{1}}&\Fhd{\eta _{_{5}}}{}&R_{_{134}}\cr
\fvb{\mu _{_{5}}}{}&(e)&\fvh{}{\lambda ^{'}_{_{123}}}\cr
R^{"}_{_{023}}&\Fhd{}{\eta^{'} _{_{5}}}&T^{'}_{_{123}}}$}
\hskip 5mm D'autre part on a les isomorphismes dans
$\Phi _{_{1}}{\times}_{_{\Phi _{_{0}}}}\Phi _{_{1}}{\times}_{_{\Phi
_{_{0}}}}\Phi _{_{1}}$
suivants :\par\centerline{$\diagram{\delta
^{^{0}}_{_{[3]}}(T)&\Fhd{a_{_{3}}(f,g,kh)}{}&
(f,g,kh)&\Fhg{(a^{^{1}}_{_{4}},a^{^{2}}_{_{4}},\delta ^{^{0}}_{_{02}}
(\eta _{_{2}})^{^{-1}})}{}&\delta ^{^{0}}_{_{[3]}}(R_{_{0124}})
\cr &&&&	\cr\delta ^{^{0}}_{_{[3]}}(T^{'})&\Fhd{a_{_{3}}(f,hg,k)}{}&(f,hg,k)&
\Fhg{(a^{^{1}}_{_{4}},\delta ^{^{0}}_{_{02}}(\eta
_{_{3}})^{^{-1}},a^{^{3}}_{_{4}})}{}&
\delta ^{^{0}}_{_{[3]}}(R_{_{0134}})\cr &&&&\cr
\delta
^{^{0}}_{_{[3]}}(T^{"})&\Fhd{a_{_{3}}(gf,h,k)}{}&(gf,h,k)&\Fhg{(\delta
^{^{0}}_{_{02}}
(\eta _{_{1}})^{^{-1}},a^{^{3}}_{_{4}},a^{^{4}}_{_{4}})}{}&\delta
^{^{0}}_{_{[3]}}(R_{_{0234}})
\cr &&&&\cr \delta ^{^{0}}_{_{[3]}}(R^{'})&\Fhd{a_{_{3}}(f,g,h)}{}&(f,g,h)&
\Fhg{(a^{^{1}}_{_{4}},a^{^{2}}_{_{4}},a^{^{3}}_{_{4}})}{}&\delta
^{^{0}}_{_{[3]}}(R_{_{0123}})
\cr &&&&\cr \delta ^{^{0}}_{_{[3]}}(R^{"})&\Fhd{a_{_{3}}(g,h,k)}{}&(g,h,k)&
\Fhg{(a^{^{2}}_{_{4}},a^{^{3}}_{_{4}},a^{^{4}}_{_{4}})}{}&
\delta ^{^{0}}_{_{[3]}}(R_{_{1234}})}$}\par
qui sont repr\'esent\'es dans $\Phi _{_{3}}$ par les isomorphismes :\par
\centerline{$\diagram{T&\Fhd{\lambda }{}&R_{_{0124}}\cr T^{'}&\Fhd{\lambda
^{'}}{}&
R_{_{0134}}\cr T^{"}&\Fhd{\lambda ^{"}}{}&R_{_{0234}}}$\hskip 1cm
$\diagram{R^{'}&\Fhd{\varepsilon ^{'}}{}&R_{_{0123}}\cr
R^{"}&\Fhd{\varepsilon ^{"}}{}&
R_{_{1234}}}$}\par
\hskip 5mm On v\'erifie facilement qu'on a les relations :\par\vskip 3mm
\centerline{$\matrix{\lambda _{_{012}}=\delta _{_{012}}^{^{1}}(\lambda )&&
\varepsilon ^{'}_{_{012}}=
\delta _{_{012}}^{^{1}}(\varepsilon ^{'})\cr
\lambda _{_{123}}=\delta _{_{123}}^{^{1}}(\lambda )&&\varepsilon ^{"}_{_{123}}
=\lambda ^{"}_{_{123}}\cdot\eta ^{'}_{_{2}}\cr
\lambda ^{'}_{_{123}}=\delta _{_{123}}^{^{1}}(\lambda ^{'})&&\varepsilon
^{"}_{_{023}}
=\lambda ^{'}_{_{123}}\cdot\eta ^{'}_{_{5}}\cr
\lambda ^{"}_{_{123}}=\delta _{_{123}}^{^{1}}(\lambda ^{"})&&
\varepsilon ^{"}_{_{013}}=\lambda _{_{123}}\cdot\eta
^{'}_{_{4}}\cr}$}\par	\vskip 3mm
en remarqant que leur images par $\delta _{_{[2]}}^{^{1}}$ coinsident et,
en appliquant le
fait que $\delta _{_{[2]}}^{^{1}}$ est injective, on obteint les relations
$(*)$ suivantes :
\par \vskip 4mm
\centerline{$\matrix{\delta _{_{02}}^{^{1}}(\varepsilon ^{"}_{_{023}})&=&
\varepsilon ^{"}_{_{03}}=\delta _{_{02}}^{^{1}}(\lambda ^{'}_{_{123}})\cdot
\delta _{_{02}}^{^{1}}(\eta ^{'}_{_{5}})=\delta _{_{02}}^{^{1}}(\eta
_{_{5}})\cdot
\delta _{_{02}}^{^{1}}(\mu _{_{5}})^{^{-1}}\cr\cr
\delta _{_{02}}^{^{1}}(\varepsilon ^{"}_{_{013}})&=&
\varepsilon ^{"}_{_{03}}=\delta _{_{02}}^{^{1}}(\lambda _{_{123}})\cdot
\delta _{_{02}}^{^{1}}(\eta ^{'}_{_{4}})=\delta _{_{02}}^{^{1}}(\eta
_{_{4}})\cdot
\delta _{_{02}}^{^{1}}(\mu _{_{4}})^{^{-1}}}$\hskip 5mm (*)}\vskip 4mm\hskip
5mm
On d\'efinit $\gamma _{_{1}}$ ,  $\gamma _{_{2}}$ et  $\gamma _{_{3}}$
respectivement par la compos\'ee de chaqu'une des
paires des applications composables suivantes :\par
\centerline{$\diagram{T_{_{023}}&\Fhd{\lambda _{_{023}}}{}&R_{_{024}}&
\Fhd{(\lambda ^{"}_{_{013}})^{^{-1}}}{}&T_{_{013}}^{"}\cr
T_{_{013}}&\Fhd{\lambda _{_{013}}}{}&R_{_{014}}&
\Fhd{(\lambda ^{'}_{_{013}})^{^{-1}}}{}&T_{_{013}}^{'}\cr
T^{'}_{_{023}}&\Fhd{\lambda ^{'}_{_{023}}}{}&R_{_{034}}&
\Fhd{(\lambda ^{"}_{_{023}})^{^{-1}}}{}&T_{_{023}}^{"}\cr}$}
\par\hskip 5mm V\'erifions maintenant la commutativit\'e des diagrammes
(1) , (2) et (3) ci-dessous :\par
\centerline{$\diagram{\delta ^{^{0}}_{_{[2]}}(\sigma _{_{2}})&
\Fhd{a_{_{2}}(gf,kh)}{}&(gf,kh)&\Fhg{(\delta ^{^{1}}_{_{02}}
(\eta ^{'} _{_{1}}\cdot\mu _{_{1}})^{^{-1}},a^{^{3}}_{_{3}})}{}&
\delta ^{^{0}}_{_{[2]}}(T_{_{023}})\cr \fvb{Id}{}&&\fvb{Id}{}&(1)&\fvh{}{
\delta ^{^{1}}_{_{[2]}}(\gamma _{_{1}})}\cr \delta ^{^{0}}_{_{[2]}}(\sigma
_{_{2}})&
\Fhd{}{a_{_{2}}(gf,kh)}&(gf,kh)&\Fhg{}{(a^{^{1}}_{_{3}},\delta ^{^{1}}_{_{02}}
(\eta ^{'} _{_{2}}\cdot\mu _{_{2}})^{^{-1}})}&\delta
^{^{0}}_{_{[2]}}(T^{"}_{_{013}})}$}\par
En utilisant les diagrammes (a) et (b) on obtient les \'egalit\'es
suivantes :\par\vskip 4mm
\centerline{$\matrix{ \delta ^{^{1}}_{_{[2]}}(\gamma _{_{1}})&=&\delta
^{^{1}}_{_{[2]}}
((\lambda ^{"}_{_{013}})^{^{-1}}\cdot\lambda _{_{023}})=\delta ^{^{1}}_{_{[2]}}
((\lambda ^{"}_{_{01}})^{^{-1}}\cdot\lambda _{_{02}},(\lambda
^{"}_{_{13}})^{^{-1}}
\cdot\lambda _{_{23}})\cr \cr diagrammes\ (a),
(b)&=&(a^{^{1}}_{_{3}}(gf,h,k)\cdot
\delta ^{^{1}}_{_{02}}(\eta ^{'}_{_{1}}\cdot\mu _{_{1}})^{^{-1}},
\delta ^{^{1}}_{_{02}}(\eta ^{'}_{_{2}}\cdot\mu _{_{2}})^{^{-1}}
\cdot a^{^{3}}_{_{3}}(f,g,kh))}$}\par\vskip 4mm
D'o\`u la commutativit\'e du diagramme (1). \par
\centerline{$\diagram{\delta ^{^{0}}_{_{[2]}}(\tau _{_{2}})&
\Fhd{a_{_{2}}(f,(kh)g)}{}&(f,(kh)g)&\Fhg{(a^{^{1}}_{_{3}},
\delta ^{^{1}}_{_{02}}(\eta ^{'} _{_{4}}\cdot\mu _{_{4}})^{^{-1}})}{}&
\delta ^{^{0}}_{_{[2]}}(T_{_{023}})\cr \fvb{\delta ^{^{1}}_{_{02}}
(\varepsilon
_{_{I_{_{f}},A(g,h,k)}})}{}&&\fvb{(I_{_{f}},A(g,h,k))}{}&(2)&\fvh{}{
\delta ^{^{1}}_{_{[2]}}(\gamma _{_{2}})}\cr \delta ^{^{0}}_{_{[2]}}(\sigma
_{_{2}})&
\Fhd{}{a_{_{2}}(f,k(hg))}&(f,k(hg))&\Fhg{}{(a^{^{1}}_{_{3}},
\delta ^{^{1}}_{_{02}}(\eta ^{'} _{_{5}}\cdot\mu _{_{5}})^{^{-1}})}&
\delta ^{^{0}}_{_{[2]}}(T^{"}_{_{013}})}$}\par
De m\^eme en utilisant la commutativit\'e
des diagrammes (e) , (d) et les relations (*) on obtient :\par\vskip 3mm
\centerline{$\matrix{ \delta ^{^{1}}_{_{[2]}}(\gamma _{_{2}})&=&\delta
^{^{1}}_{_{[2]}}
((\lambda ^{'}_{_{013}})^{^{-1}}\cdot\lambda _{_{013}})=\delta ^{^{1}}_{_{[2]}}
((\lambda ^{'}_{_{01}})^{^{-1}}\cdot\lambda _{_{01}},
(\lambda ^{'}_{_{13}})^{^{-1}}\cdot\lambda _{_{13}})\cr \cr
diagrammes\ (e), (d)&=&(a^{^{1}}_{_{3}}(f,hg,k)^{^{-1}}\cdot a^{^{1}}_{_{4}}
\cdot (a^{^{1}}_{_{3}})^{^{-1}}\cdot a^{^{1}}_{_{3}}(f,g,kh)
, \delta ^{^{1}}_{_{02}}(\lambda ^{'}_{_{123}})^{^{-1}}\cdot
\delta ^{^{1}}_{_{02}}(\lambda _{_{123}})\cr\cr
relation\  (*) &=&(a^{^{1}}_{_{3}}(f,hg,k)^{^{-1}}
\cdot a^{^{1}}_{_{4}}\cdot (a^{^{1}}_{_{3}})^{^{-1}}\cdot
a^{^{1}}_{_{3}}(f,g,kh),
\delta ^{^{1}}_{_{02}}(\eta ^{'}_{_{5}})\cdot \delta ^{^{1}}_{_{02}}
(\eta ^{'}_{_{4}})^{^{-1}})}$}\par\vskip 4mm
On sait aussi que $A(g,h,k)= \delta ^{^{1}}_{_{02}}(\mu _{_{5}})^{^{-1}}\cdot
\delta ^{^{1}}_{_{02}}(\mu _{_{4}})$, ce qui montre la commutativit\'e du
diagramme (2).
La d\'emarche utilis\'ee pour montrer la commutativit\'e du diagramme (2)
est aussi
valable pour le diagramme (3) suivant :\par
\centerline{$\diagram{\delta ^{^{0}}_{_{[2]}}(\tau _{_{2}})&
\Fhd{a_{_{2}}((hg)f,k)}{}&((hg)f,k)&\Fhg{}{}&
\delta ^{^{0}}_{_{[2]}}(T^{'}_{_{023}})\cr
\fvb{\delta ^{^{1}}_{_{02}}(\varepsilon _{_{A(f,g,h),I_{_{f}}}})}{}&&
\fvb{(A(f,g,h),I_{_{f}})}{}&(3)&\fvh{}{}
\delta ^{^{1}}_{_{[2]}}(\gamma _{_{3}})\cr \delta
^{^{0}}_{_{[2]}}(\sigma^{"} _{_{2}})&
\Fhd{}{a_{_{2}}(h(gf),k)}&(h(gf),k)&\Fhg{}{}&\delta
^{^{0}}_{_{[2]}}(T^{"}_{_{023}})}$}\par
Les diagrammes (1) (2) et (3) sont obtenus en appliquant $\delta _{_{[2]}}$
aux diagrammes dans $\Phi _{_{2}}$ suivants :\par
\centerline{$\diagram{\tau _{_{2}}&\fhd{}{}&T _{_{013}}\cr
\fvb{\varepsilon _{_{I_{_{f}},A(g,h,k)}}}{}&(2)^{'}&
\fvb{}{\gamma _{_{2}}}\cr \tau ^{'}_{_{2}}&\fhd{}{}&T ^{'}_{_{013}}}$\hskip 1cm
$\diagram{T^{'}_{_{023}}&\fhd{\gamma _{_{3}}}{}&T^{"}
_{_{023}}\cr\fvh{}{}&(3)^{'}&
\fvh{}{}\cr \sigma ^{'}_{_{2}}&\fhd{}{\varepsilon _{_{A(f,g,h),I_{_{f}}}}}&
\sigma ^{"}_{_{2}}}$\hskip 1cm
$\diagram{T _{_{023}}&\fhg{}{}&\sigma _{_{2}}\cr\fvb{\gamma _{_{1}}}{}&(1)^{'}&
\fvb{}{Id}\cr T ^{"}_{_{013}}&\fhd{}{}&\sigma _{_{2}}}$}\par
Puisque $\delta _{_{[2]}}$ est une 1- equivalence ext\'erieure, ces diagrammes
sont commutatifs. En appliquant $\delta _{_{02}}$  aux diagrammes
$(1)^{'}$ , $(2)^{'}$ et $(3)^{'}$
on obtient trois diagrammes dans la cat\'egorie $\Phi _{_{1}}$ qui on obtient
un quatri\`eme diagramme qui est commutatif d'apr\`es la d\'efinition de
$\gamma _{_{1}}$ , $\gamma _{_{2}}$ et $\gamma _{_{3}}$ :	\par
\centerline{$\diagram{T_{_{03}}&\fhd{Id}{}&T _{_{03}}\cr
\fvh{\delta ^{^{1}}_{_{02}}(\gamma _{_{2}})}{}&&
\fvh{}{\delta ^{^{1}}_{_{02}}(\gamma _{_{1}})}\cr
T^{'}_{_{03}}&\fhd{}{\delta ^{^{1}}_{_{02}}(\gamma
_{_{3}})}&T^{"}_{_{03}}}$}\par
Finalement le diagramme globale form\'e par ce dernier et les images des
diagrammes
$(1)^{'}$ , $(2)^{'}$ et $(3)^{'}$ est le diagramme commutatif suivant :\par
\centerline{$\diagram{\delta ^{^{0}}_{_{02}}(\tau _{_{2}})&\Fhd{A(f,g,kh)}{}&
\delta ^{^{0}}_{_{02}}(\sigma _{_{2}})&\Fhd{A(gf,h,k)}{}&
\delta ^{^{0}}_{_{02}}(\sigma ^{"}_{_{2}})\cr
\fvb{A(g,h,k)\star I_{_{f}}}{}&&&&\fvb{}{Id}\cr\delta ^{^{0}}_{_{02}}(\tau
^{'}_{_{2}})&
\Fhd{}{A(f,hg,k)}&\delta ^{^{0}}_{_{02}}(\sigma ^{'}_{_{2}})&
\Fhd{}{I_{_{k}}\star A(f,g,h)}&\delta ^{^{0}}_{_{02}}(\sigma
^{"}_{_{2}})}$}\par
ce qui entraine l'axiome (1)\hfill


\par\vskip 5mm\hskip 5mm\underbar{\bf Axiome (2)} On sait que $\Phi
_{_{1}}$ est un
1-nerf d'apr\`es la proposition (1.1.4) donc la composition dans cette
cat\'egorie est
associative, ce qui montre que la composition verticale des 2-fl\`eches de
$\Phi $ est
associative. \hfill


\par\vskip 5mm\hskip 5mm\underbar{\bf Axiome (3)} (Relation de Godement):
Soient $(\alpha ,\beta )$ et $(\alpha ^{'},\beta ^{'})$ deux \'el\'ements de
${\cal C}_{_{2}}{\times}_{_{{\cal C}_{_{0}}}}{\cal C}_{_{2}}$ tels que
$bb(\alpha )=ss(\beta ^{'})$. On peut les repr\'esenter de la fa\c con
suivante :\par
\centerline{$\diagram{\Fhd{f}{}&\Fhd{f^{'}}{}\cr \fvb{}{\alpha
}&\fvb{}{\alpha ^{'}}\cr
\Fhd{}{g}&\Fhd{}{g^{'}}\cr \fvb{}{\beta }&\fvb{}{\beta ^{'}}\cr
\Fhd{}{h}&\Fhd{}{h^{'}}}$\hskip 5mm avec \hskip 5mm $\left\{\matrix{
f&=&s(\alpha )\cr g&=&s(\beta )&=&b(\alpha )\cr h&=&b(\beta )\cr
f^{'}&=&s(\alpha ^{'})\cr g^{'}&=&s(\beta ^{'})&=&b(\alpha ^{'})\cr
h^{'}&=&b(\beta ^{'})}
\right.$}\par Posons $\sigma =L_{_{2}}(f,f^{'})$ , $\sigma
^{'}=L_{_{2}}(g,g^{'})$  et
$\sigma ^{"}=L_{_{2}}(h,h^{'})$. On sait d'apr\`es le paragraphe (c) qu'ils
existent
trois \'el\'ements de $\Phi _{_{2,1}}$ uniques :\par\vskip 3mm
\centerline{$\varepsilon _{_{\alpha ,\alpha ^{'}}}:
\sigma \fhd{}{}\sigma ^{'}$\hskip 1cm
$\varepsilon _{_{\beta  ,\beta ^{'}}}:\sigma ^{'}\fhd{}{} \sigma ^{"}$\hskip
1cm
$\varepsilon _{_{\beta\cdot\alpha  ,\beta ^{'}\cdot\alpha ^{'}}}:
\sigma\fhd{}{}\sigma ^{"}$}
\par\vskip 3mm tels que les diagrammes (1) , (2) et celui form\'e par
leur r\'eunion, soient commutatifs :\par
\centerline{$\diagram{\delta ^{^{0}}_{_{[2]}}(\sigma )&\Fhd{\delta
^{^{1}}_{_{[2]}}(
\varepsilon _{_{\alpha ,\alpha ^{'}}})}{}&\delta ^{^{0}}_{_{[2]}}(\sigma
^{'})&\Fhd{
\delta ^{^{1}}_{_{[2]}}(\varepsilon _{_{\beta ,\beta ^{'}}})}{}&\delta
^{^{0}}_{_{[2]}}
(\sigma ^{"})\cr\Fvb{a_{_{2}}}{\wr}&(1)&\Fvb{a_{_{2}}}{\wr}&(2)&\Fvb{\wr}{
a_{_{2}}}\cr (f,f^{'})&\Fhd{}{(\alpha ,\alpha ^{'})}&
(g,g^{'})&\Fhd{}{(\beta ,\beta ^{'})}&(h,h^{'})}$}\par
On en d\'eduit que :  \par	\vskip 3mm\centerline{$\delta ^{^{1}}_{_{[2]}}(
\varepsilon _{_{\beta ,\beta ^{'}}})\cdot\delta ^{^{1}}_{_{[2]}}(
\varepsilon _{_{\alpha ,\alpha ^{'}}})=\delta ^{^{1}}_{_{[2]}}(
\varepsilon _{_{\beta\cdot\alpha  ,\beta ^{'}\cdot\alpha ^{'}}})$\hskip 5mm et
puis \hskip 5mm$\varepsilon _{_{\beta ,\beta ^{'}}}\cdot
\varepsilon _{_{\alpha ,\alpha ^{'}}}=\varepsilon _{_{\beta\cdot\alpha
,\beta ^{'}
\cdot\alpha ^{'}}}$}\par\vskip 3mm
D'o\`u finalement en appliquant $\delta _{_{02}}$ on retrouve
l'Axiome(3) :\par\vskip 3mm\centerline{$(\beta ^{'}\star\beta )
\cdot(\alpha ^{'}\star\alpha )=
(\beta ^{'}\cdot\alpha ^{'})\star(\beta \cdot\alpha )$}\hfill


\par\vskip 5mm\hskip 5mm\underbar{\bf Axiome (4) et (5)}: Soient $\varphi :
f\fhd{}{}f^{'}$
une 2-fl\`eche de $\Phi $ , $\sigma =L_{_{2}}(f,I_{_{b(f)}})$ et $\sigma
^{'} =L_{_{2}}
(f^{'},I_{_{b(f)}})$. \par Posons $\tau =\delta ^{^{0}}_{_{011}}(f)$ et
$\tau ^{'}=
\delta ^{^{0}}_{_{011}}(f^{'})$ et Consid\'erons les isomorphismes :\par
$$\diagram{\delta ^{^{0}}_{_{[2]}}(\sigma )&\Fhd{\delta ^{^{1}}_{_{[2]}}
(\varepsilon _{_{f}})}{\sim}&\delta ^{^{0}}_{_{[2]}}(\tau )&\Fhd{(\varphi ,
I^{^{2}}_{_{b(f)}})}{}&\delta ^{^{0}}_{_{[2]}}(\tau ^{'})&\Fhg{\delta
^{^{1}}_{_{[2]}}
(\varepsilon _{_{f^{'}}})}{\sim}&\delta ^{^{0}}_{_{[2]}}(\sigma ^{'})}$$
 L'application $\delta ^{^{0}}_{_{011}} : \Phi _{_{1,1}}\fhd{}{}\Phi
_{_{2,1}}$ fait
correspondre \`a $\varphi $ l'\'el\'ement $\delta ^{^{0}}_{_{011}}(\varphi
)$ tel que \par
\vskip 3mm\centerline{$\delta ^{^{0}}_{_{[2]}}(\delta
^{^{0}}_{_{011}}(\varphi ))=
(\varphi ,I^{^{2}}_{_{b(f)}})$\hskip 1cm et \hskip 1cm$\delta ^{^{0}}_{_{02}}
(\delta ^{^{0}}_{_{011}}(\varphi ))=\varphi $}\par	\vskip 3mm
D'autre part on sait qu'il existe un unique $\varepsilon :\sigma
\fhd{}{}\sigma ^{'}$
tel que \par
\vskip 3mm\centerline{$\delta ^{^{1}}_{_{[2]}}(\varepsilon )=\bigl[
\delta ^{^{1}}_{_{[2]}}(\varepsilon _{_{f^{'}}})\bigr]^{^{-1}}\cdot
\delta ^{^{1}}_{_{[2]}}(\delta ^{^{0}}_{_{011}}(\varphi ))\cdot\delta
^{^{1}}_{_{[2]}}
(\varepsilon _{_{f}})$}\par\vskip 3mm
D'o\`u  $\varepsilon _{_{f^{'}}}\cdot\varepsilon =\delta
^{^{0}}_{_{011}}(\varphi )\cdot
\varepsilon _{_{f}}$ , ce qui entraine en appliquant $\delta
^{^{0}}_{_{02}}$ \`a cette
\'egalit\'e, la relation :\par
\vskip 3mm\centerline{$V(f^{'})\cdot (I^{^{2}}_{_{b(f)}}\star
\varphi) =\varphi \cdot V(f)$}\par\vskip 3mm
Avec un raisonement analoque et, en consid\'erant $\delta
^{^{0}}_{_{001}}(\varphi )$,
on montre que la relation suivante est \'egalement satisfaite :\par
\vskip 3mm\centerline{$U(f^{'})\cdot (\varphi \star I^{^{2}}_{_{s(f)}}) =
\varphi \cdot U(f)$}\par\vskip 3mm


\par\vskip 5mm\hskip 5mm\underbar{\bf Axiome (6)} : (i) D'apr\`es le
paragraphe sur
les isomorphismes de coh\`erences d'identi\'e, on a
$\delta ^{^{0}}_{_{011}}(I_{_{x}})=\delta ^{^{0}}_{_{001}}(I_{_{x}})$  donc
$\varepsilon _{_{I_{_{x}}}}=\varepsilon ^{^{I_{_{x}}}}$ et par suite
$U(I_{_{x}})=V(I_{_{x}})$.\par
(ii) Soit $(f,g)$ dans $\Phi _{_{1,0}}{\times}_{_{\Phi _{_{0,0}}}}\Phi
_{_{1,0}}$, et posons
$\sigma =L_{_{2}}(f,g)$. On a le diagramme commutatif suivant :\par
\centerline{$\diagram{\delta ^{^{0}}_{_{[2]}}(\sigma
)&\Fhd{a_{_{2}}(f,g)}{}&(f,g)\cr
\fvb{\delta ^{^{1}}_{_{[2]}}(\varepsilon
_{_{I_{f},I_{g}}})}{}&&\fvb{}{(I_{f},I_{g})}\cr
\delta ^{^{0}}_{_{[2]}}(\sigma )&\Fhd{}{a_{_{2}}(f,g)}&(f,g)}$}\par
Donc \hskip 5mm $\delta ^{^{1}}_{_{[2]}}(\varepsilon _{_{I_{f},I_{g}}})=
I_{_{\delta ^{^{0}}_{_{[2]}}(\sigma )}}=\delta ^{^{1}}_{_{[2]}}(I_{\sigma
})$ \hskip 5mm
d'o\`u \hskip 5mm
$\varepsilon _{_{I_{f},I_{g}}}=I_{\sigma }$\hskip 5mm et par suite on a
:\par\vskip 2mm
\centerline{$I_{g}\star I_{f}= \delta ^{^{1}}_{_{02}}(\varepsilon
_{_{I_{f},I_{g}}})
=\delta ^{^{1}}_{_{02}}(I_{\sigma })=I_{_{\delta ^{^{0}}_{_{02}}(\sigma )}}
=I_{gf}$}\hfill


\par\vskip 5mm\hskip 5mm\underbar{\bf Axiome (7)}: Soient $(f,g,h)$ ,
$(f^{'},g^{'},h^{'})$
deux \'el\'ements de ${\cal C}_{_{1}}{\times}_{_{{\cal C}_{_{0}}}}{\cal
C}_{_{1}}
{\times}_{_{{\cal C}_{_{0}}}}{\cal C}_{_{1}}$ et $(\alpha ,\beta ,\gamma )$
un \'el\'ement
de ${\cal C}_{_{2}}{\times}_{_{{\cal C}_{_{0}}}}{\cal C}_{_{2}}
{\times}_{_{{\cal C}_{_{0}}}}{\cal C}_{_{2}}$ comme dans le diagramme
suivant :\par
\centerline{$\diagram{\fhd{f}{}&&\fhd{g}{}&&\fhd{h}{}\cr\fvb{}{\alpha }&&
\fvb{}{\beta }&&\fvb{}{\gamma
}\cr\fhd{}{f^{'}}&&\fhd{}{g^{'}}&&\fhd{}{h^{'}}}$}\par
Gardons les m\^emes notations utilis\'ees aux paragraphes (c) et (e), quant aux
\'el\'ements $f^{'}$, $g^{'}$ et $h^{'}$ on garde les m\^eme symboles mais
avec des primes
( c.\`a.d $T^{'}, \sigma _{_{1}}^{'}, \tau _{_{1}}^{'}, \mu _{_{1}}^{'},
\eta _{_{1}}^{'},
\sigma _{_{2}}^{'}, \tau _{_{2}}^{'}, \mu _{_{2}}^{'}, \eta _{_{2}}^{'}$ ).\par
\vskip 3mm\hskip 3mm{\bf (i)} On sait qu'il existe un unique morphisme
$\omega :T\fhd{}{}T^{'}$ tel que le diagramme suivant est commutative :\par
\centerline{$\diagram{\delta ^{^{0}}_{_{[3]}}(T)&\Fhd{a_{_{3}}(f,g,h)}{\sim}
&(f,g,h)\cr\fvb{\delta ^{^{1}}_{_{[3]}}(\omega )}{}&&\fvb{}{(\alpha ,\beta
,\gamma )}\cr
\delta
^{^{0}}_{_{[3]}}(T^{'})&\Fhd{\sim}{a_{_{3}}(f^{'},g^{'},h^{'})}&(f^{'},g^{'}
,h^{'})}$}
Consid\'erons les deux diagrammes commutatifs :\par
\centerline{$\diagram{\delta ^{^{0}}_{_{[2]}}(\tau
_{_{1}})&\Fhd{a_{_{2}}(g,h)}{\sim}
&(g,h)&\Fhg{(a^{^{2}}_{_{3}},a^{^{3}}_{_{3}})}{\sim}&\delta
^{^{0}}_{_{[2]}}(T_{_{123}})\cr
\fvb{\delta ^{^{1}}_{_{[2]}}(\varepsilon _{_{\beta ,\gamma
}})}{}&&\fvb{}{(\beta ,\gamma )}
&&\fvb{}{(\delta ^{^{1}}_{_{12}(\omega )},\delta ^{^{1}}_{_{23}}(\omega ))}\cr
\delta ^{^{0}}_{_{[2]}}(\tau
_{_{1}}^{'})&\Fhd{\sim}{a_{_{2}}(g^{'},h^{'})}&(g^{'},h^{'})&
\Fhg{\sim}{(a^{^{2}}_{_{3}},a^{^{3}}_{_{3}})}&\delta
^{^{0}}_{_{[2]}}(T^{'}_{_{123}})}$}
Le morphisme $\omega :T\fhd{}{}T^{'}$ induit un morphisme
$\delta ^{^{1}}_{_{123}}(\omega ) : T_{_{123}}\fhd{}{}T^{'}_{_{123}}$ tel que
$\delta ^{^{1}}_{_{[2]}}(\delta ^{^{1}}_{_{123}}(\omega ))=
(\delta ^{^{1}}_{_{12}}(\omega ),\delta ^{^{1}}_{_{23}}(\omega ))$. Donc
d'apr\`es les
paragraphes (c) et (e), le diagramme pr\'ec\'edent entraine le diagramme
commutative :\par
\centerline{$\diagram{\tau _{_{1}}&\Fhg{\eta _{_{1}}}{}
& T_{_{123}}\cr\fvb{\varepsilon _{_{\beta ,\gamma }}}{}&&
\fvb{}{\delta ^{^{1}}_{_{123}}(\omega )}\cr
\tau _{_{1}}^{'}&\Fhg{}{\eta ^{'}_{_{1}}}&T^{'}_{_{123}}}$}\par
qui, en lui appliquant $\delta _{{02}}$, donne le diagramme commutative
suivant :\par
\centerline{$\diagram{hg&\Fhg{\delta ^{^{1}}_{_{02}}(\eta _{_{1}})}{}&
\delta ^{^{0}}_{_{13}}(T)\cr\fvb{\gamma \star\beta }{}&&\fvb{}{\delta
^{^{1}}_{_{13}}
(\omega )}\cr h^{'}g^{'}&\Fhg{}{\delta ^{^{1}}_{_{02}}(\eta ^{'}_{_{1}})}&
\delta ^{^{0}}_{_{13}}(T^{'})}$}\par
\par\vskip 3mm\hskip 3mm{\bf (ii)} Consid\'erons maintenant les deux
diagrammes commutatifs suivants :\par
\centerline{$\diagram{\delta ^{^{0}}_{_{[2]}}(\tau
_{_{2}})&\Fhd{a_{_{2}}(f,hg)}{\sim}
&(f,hg)&\Fhg{(a^{^{1}}_{_{3}},\delta ^{^{1}}_{_{02}}(\eta _{_{1}}))}{\sim}&
\delta ^{^{0}}_{_{[2]}}(T_{_{013}})\cr\fvb{\delta
^{^{1}}_{_{[2]}}(\varepsilon _{_{\alpha ,
\gamma \star\beta }})}{}&&\fvb{}{(\alpha ,\gamma \star\beta )}&&\fvb{}{(
\delta ^{^{1}}_{_{01}}(\omega ),\delta ^{^{1}}_{_{13}}(\omega ))}\cr
\delta ^{^{0}}_{_{[2]}}(\tau _{_{2}}^{'})&\Fhd{\sim}{a_{_{2}}
(f^{'},h^{'}g^{'})}&(f^{'},h^{'}g^{'})&\Fhg{\sim}{(a^{^{1}}_{_{3}},
\delta ^{^{1}}_{_{02}}(\eta ^{'}_{_{1}}))}&\delta
^{^{0}}_{_{[2]}}(T^{'}_{_{013}})}$}\par
Comme $\delta ^{^{1}}_{_{[2]}}(\omega )((\delta ^{^{1}}_{_{123}}(\omega ))
=(\delta ^{^{1}}_{_{01}}(\omega ),\delta ^{^{1}}_{_{13}}(\omega ))$ alors
d'apr\`es (c) et
(e) on obtient le diagramme commutatif :\par
\centerline{$\diagram{\tau _{_{2}}&\Fhd{\eta _{_{2}}}{}
& T_{_{013}}\cr\fvb{\varepsilon _{_{\alpha ,\gamma \star\beta }}}{}&&
\fvb{}{\delta ^{^{1}}_{_{013}}(\omega )}\cr
\tau _{_{2}}^{'}&\Fhd{}{\eta ^{'}_{_{2}}}&T^{'}_{_{013}}}$}\par
En lui appliquant $\delta _{_{02}}$, on d\'eduit le diagramme commutatif :\par
\centerline{$\diagram{(hg)f&\Fhd{\delta ^{^{1}}_{_{02}}(\eta _{_{2}})}{}&
\delta ^{^{0}}_{_{03}}(T)\cr\fvb{(\gamma \star\beta )\star\alpha }{}&(a)
&\fvb{}{\delta ^{^{1}}_{_{03}}(\omega )}\cr
(h^{'}g^{'})f^{'}&\Fhd{}{\delta ^{^{1}}_{_{02}}(\eta ^{'}_{_{2}})}&
\delta ^{^{0}}_{_{03}}(T^{'})}$}\par
\vskip 3mm\hskip 3mm{\bf (iii)} Une construction analogue au pr\'ec\'edente
avec
($\sigma _{_{i}}, \mu _{_{i}}, T_{_{012}}, T_{_{023}}$) et \par
($\sigma ^{'}_{_{i}}, \mu ^{'}_{_{i}}, T^{'}_{_{012}}, T^{'}_{_{023}}$)
permet d'aboutir au
diagramme commutatif :\par\centerline{$\diagram{T_{_{023}}&\Fhg{\mu _{_{2}}}{}
&\sigma _{_{2}}\cr\fvb{\delta ^{^{1}}_{_{023}}(\omega )}{}&&
\fvb{}{\varepsilon _{_{\beta \star\alpha ,\gamma }}}\cr T_{_{023}}^{'}&
\Fhg{}{\mu _{_{2}}^{'}}&\sigma ^{'}_{_{2}}}$}\par
Et par cons\'equent le diagramme commutatif : \par
\centerline{$\diagram{\delta ^{^{0}}_{_{03}}(T)&\Fhg{\delta
^{^{1}}_{_{02}}(\mu _{_{2}})}{}&
h(gf)\cr\fvb{\delta ^{^{1}}_{_{03}}(\omega )}{}&(b)&\fvb{}{\gamma
\star(\beta \star\alpha )}\cr
\delta ^{^{0}}_{_{03}}(T^{'})&\Fhg{}{\delta ^{^{0}}_{_{02}}(\mu ^{'}_{_{2}})}&
h^{'}(g^{'}f^{'})}$}\par
La commutativit\'e des diagrammes (a) et (b) entraine alors l'Axiome (7).\par
\vskip 2mm\centerline{$\bigl[\gamma \star(\beta \star\alpha )\bigr]\cdot
A(f,g,h)=
A(f^{'},g^{'},h^{'})\cdot\bigl[(\gamma \star\beta )\star\alpha
\bigr]$}\hfill \par


\par\vskip 5mm\hskip 5mm\underbar{\bf Axiomes (8), (9) et (10)}: Les
Axiomes (8), (9)
et (10) se resemblent donc il suffit de montrer l'un de ces Axiomes et pour
les autres la
m\^eme d\'emarche est valable. V\'erifions alors l'Axiome (8). Soit
$(f,I_{_{y}},g)$ dans
${\cal C}_{_{1}}{\times}_{_{{\cal C}_{_{0}}}}{\cal
C}_{_{1}}{\times}_{_{{\cal C}_{_{0}}}}
{\cal C}_{_{1}}$, repr\'esent\'es par le diagramme :\par\vskip 3mm
\centerline{$x\Fhd{f}{}y\Fhd{I_{_{y}}}{}z\Fhd{g}{}t$} \par\vskip 2mm\hskip 5mm
Dans ce paragraphe on r\'ef\`ere de m\^eme aux notations des paragraphes
(c) et (e).
En plus on a besoin des notations : $T=(f,I_{_{y}},g)$ ,
$\rho =(f,g)$ , $F=\delta _{_{011}}(f)$ , $G=\delta _{_{001}}(g)$ et $T^{'}=
\delta ^{^{0}}_{_{0112}}(\rho )$. On sait par le fait que $\Phi _{_{3}}$
est un 1-nerf
qu'il existe un unique $\omega :T\fhd{}{}T^{'}$ tel que $\delta
^{^{1}}_{_{[3]}}(\omega )$
soit \'egale au compos\'e des isomorphismes suivants :\par
\centerline{$\diagram{\delta
^{^{1}}_{_{[3]}}(T)&\Fhd{a_{_{3}}(f,I_{_{y}},g)}{\sim}&
(f,I_{_{y}},g)&\Fhg{(a^{^{1}}_{_{2}},I^{^{2}}_{_{y}},a^{^{2}}_{_{2}})}{\sim}&
(\delta ^{^{0}}_{_{02}}(\rho ),I^{^{2}}_{_{y}},\delta ^{^{0}}_{_{12}}(\rho ))=
\delta ^{^{1}}_{_{[3]}}(T^{'})}$}\par
Soit $\eta ^{"}_{_{1}}=\delta ^{^{0}}_{_{001}}(a^{^{2}}_{_{2}}(f,g)$ , donc
$\delta ^{^{0}}_{_{[2]}}(\eta
^{"}_{_{1}})=(I^{^{2}}_{_{y}},a^{^{2}}_{_{2}}(f,g))$ et
$\delta ^{^{0}}_{_{02}}(\eta ^{"}_{_{1}})=a^{^{2}}_{_{2}}(f,g)$. On obtient
alors un
diagramme commutatif :\par\centerline{$\diagram{&&\delta
^{^{0}}_{_{[2]}}(T_{_{123}})&
\Fhd{\delta ^{^{1}}_{_{[2]}}(\omega _{_{123}})}{}&\delta
^{^{0}}_{_{[2]}}(T^{'}_{_{123}})
\cr &&\Fvb{(a^{^{2}}_{_{3}},a^{^{3}}_{_{3}})}{\wr}&&\Fvb{\wr}{\delta
^{^{1}}_{_{[2]}}
(\eta ^{"}_{_{1}})}\cr \delta ^{^{0}}_{_{[2]}}(\tau _{_{1}})&\Fhd{\sim}
{\delta ^{^{1}}_{_{[2]}}(\varepsilon ^{^{g}})}&(I_{_{y}},g)&=&\delta
^{^{0}}_{_{[2]}}(G)}$}\par
Il existe un unique $\eta ^{'}_{_{1}}:T_{_{123}}\fhd{}{}G$ tel que $\delta
^{^{0}}_{_{[2]}}
(\eta ^{'}_{_{1}})=(a^{^{2}}_{_{3}},a^{^{3}}_{_{3}})$ , ce qui implique le
diagramme
commutatif :\par
\centerline{$\diagram{T_{_{123}}&\Fhd{\omega _{_{123}}}{}&T_{_{123}}^{'}\cr
\fvb{}{\eta ^{'}_{_{1}}}&&\fvb{}{\eta ^{"}_{_{1}}}\cr
G&=&G}$\hskip 1cm puis \hskip 1cm
$\diagram{\delta ^{^{0}}_{_{13}}(T)&\Fhd{\omega _{_{13}}}{}&\delta
^{^{0}}_{_{13}}(T^{'})
\cr \fvb{}{\delta ^{^{0}}_{_{02}}(\eta
^{'}_{_{1}})}&&\fvb{}{(a^{^{2}}_{_{3}},a^{^{3}}_{_{3}})}
\cr g&=&g}$}\par
Soit maintenant les diagrammes :\par
\centerline{$\diagram{\delta ^{^{0}}_{_{[2]}}(\tau
_{_{2}})&\Fhd{a_{_{2}}(f,gI_{_{y}})}{}&
(f,gI_{_{y}})&\Fhg{(a^{^{1}}_{_{3}},\delta ^{^{1}}_{_{02}}(\eta _{_{1}})}{}&
\delta ^{^{0}}_{_{[2]}}(T_{_{013}})
\cr \Fvb{\delta ^{^{1}}_{_{[2]}}(\varepsilon ^{^{g}})}{}&(1)&
\Fvb{(a^{^{2}}_{_{3}},a^{^{3}}_{_{3}})}{\wr}&(2)&\Fvb{\wr}{\delta
^{^{1}}_{_{[2]}}
(\eta ^{"}_{_{1}})}\cr \delta ^{^{0}}_{_{[2]}}(\rho )&\Fhd{\sim}
{a_{_{[2}}(f,g)}&(f,g)&\Fhg{\sim}{a_{_{2}}(f,g)}&\delta
^{^{0}}_{_{[2]}}(\rho )}$}\par
La commutativit\'e de (1) est d'apr\`es le paragraphe (c), quant au
diagramme (2) sa
commutativit\'e est due au fait qu'on a :\par\vskip 2mm
\centerline{$\matrix{\eta _{_{1}}&=&[\varepsilon ^{^{g}}]^{^{-1}}\cdot \eta
_{_{1}}^{'}\cr\cr
T^{'}_{_{013}}&=&\delta ^{^{0}}_{_{013}}(T^{'})= \rho \cr\cr \omega _{_{13}}&=&
[a^{^{2}}_{_{2}}(f,g)]^{^{-1}}\cdot U(g)\cdot\delta ^{^{0}}_{_{02}}(\eta
_{_{1}})}$}
\par\vskip 2mm On en d\'eduit le diagramme commutatif :\par
\centerline{$\diagram{\tau _{_{2}}&\Fhd{\eta _{_{2}}}{}&T_{_{013}}\cr
\fvb{\varepsilon _{_{I_{_{f}},U(g)}}}{}&&\fvb{}{\omega _{_{013}}}\cr
\rho &=&\rho }$\hskip 1cm puis \hskip 1cm
$\diagram{(gI_{_{y}})f&\Fhd{\delta ^{^{0}}_{_{02}}(\eta _{_{2}})}{}&\delta
^{^{0}}_{_{03}}(T)
\cr \fvb{U(g)\star I_{_{f}}}{}&(a)&\fvb{}{(\omega _{_{03}})}\cr gf&=&gf}$}\par
D'autre part un raisonement analogue au pr\'ec\'edent permet d'obtenir des
diagrammes du
m\^eme genre :\par
\centerline{$\diagram{T_{_{023}}&\Fhd{\mu _{_{2}}}{}&\sigma _{_{2}}\cr
\fvb{\omega _{_{023}}}{}&&\fvb{}{\varepsilon _{_{V(f),I_{_{g}}}}}\cr
\rho &=&\rho }$\hskip 1cm puis \hskip 1cm
$\diagram{\delta ^{^{0}}_{_{03}}(T)&\Fhd{\delta ^{^{0}}_{_{02}}(\eta
_{_{2}})}{}&g(I_{_{y}}f)
\cr \fvb{\omega _{_{03}}}{}&(b)&\fvb{}{I_{_{g}}\star V(f)}\cr gf&=&gf}$}\par
On en d\'eduit alors de (a) et (b) l'Axiome  (8) : \par\vskip 3mm
\centerline{$\bigl[I{_{g}}\star V(f)\bigr]\cdot A(f,I_{_{y}},g)=U(g)\star
I_{_{f}}$}\hfill
\par\vskip 5mm\hskip 5mm


{\bf Remarques (1.4.4) :} (1) Soit ${\cal C}$ une 2-cat\'egorie large, et
consid\'erons
la 2-cat\'egorie large ${\cal C}^{'}$, qui correspond \`a son nerf double.
Il est claire qu'on a
${\cal C}_{_{0}}={\cal C}^{'}_{_{0}}$ , ${\cal C}_{_{1}}={\cal
C}^{'}_{_{1}}$ et
${\cal C}_{_{2}}={\cal C}^{'}_{_{2}}$, la composition verticale des
2-fl\`eches est
la m\^eme dans ${\cal C}$ et ${\cal C}^{'}$, mais la composition des
fl\`eches dans ${\cal C}$
n'est pas forcement la m\^eme que celle dans ${\cal C}^{'}$, ils sont
isomorphes. Cela se
traduit par le fait qu'il y a un 2-foncteur large ([4] O. Leroy) entre
${\cal C}$ et ${\cal C}^{'}$
qui admet un inverse stricte (ou exacte).\par
(2) Soit maintenant $\Phi $ un 2-nerf, et consid\'erons le nerf double
$\Psi $ de la 2-cat\'egorie large ${\cal C}$ qui correspond \`a $\Phi $.
Alors ils existent
un 2-nerf strict ${\cal S}$ et deux 2-\'equivalences ext\'erieures
$\alpha $ , $\beta $ comme suite :\par
\centerline{$\diagram{\Phi &\Fhd{\alpha }{\sim}&{\cal S}&\Fhg{\beta
}{\sim}&\Psi }$}
En effet on d\'efinit ${\cal S}$ pour tout objet $(m,n)$ de $\Delta ^{2}$
par :\par\vskip 3mm
\centerline{$\matrix{{\cal S}_{_{0,0}}={\cal C}_{_{0}}\hskip 5mm,\hskip 5mm
{\cal S}_{_{1,0}}={\cal C}_{_{1}}\hskip 5mm,\hskip 5mm{\cal
S}_{_{1,1}}={\cal C}_{_{2}}
\cr\cr {\cal S}_{_{m,1}}={\cal S}_{_{1,1}}{\times}_{_{{\cal S}_{_{0,0}}}}
\dots{\times}_{_{{\cal S}_{_{0,0}}}}{\cal S}_{_{1,1}}\hskip 1cm (m fois)
\cr\cr {\cal S}_{_{m,0}}={\cal S}_{_{1,0}}{\times}_{_{{\cal S}_{_{0,0}}}}
\dots{\times}_{_{{\cal S}_{_{0,0}}}}{\cal S}_{_{1,0}}\hskip 1cm (m fois)
\cr\cr{\cal S}_{_{m,n}}={\cal S}_{_{m,1}}{\times}_{_{{\cal S}_{_{m,0}}}}
\dots{\times}_{_{{\cal S}_{_{m,0}}}}{\cal S}_{_{m,1}}\hskip 1cm (n
fois)}$}\par\vskip 2mm
\hskip 5mm Soit $(m,n)$ un objet de $\Delta ^{^{2}}$, on d\'efinit les
images des
applications $d^{^{k}}_{_{i}}$ et $\varepsilon ^{^{k}}_{_{i}}$, pour
$(k,i)$ dans
$\{1\}\times\{0,\dots,m\}$ ou dans $\{2\}\times\{0,\dots,n\}$, par
:\par\vskip 5mm
\centerline{$\matrix{{\cal S}_{_{m,n}}&\Fhd{d^{^{1}}_{_{i}}}{}&{\cal
S}_{_{m-1,n}}\cr\cr
(f^{^{\alpha }}_{_{01}},\dots,f^{^{\alpha }}_{_{m-1,m}})_{_{1\leq \alpha
\leq n}}&\Fhd{}{}&
\matrix{(f^{^{\alpha }}_{_{12}},\dots,f^{^{\alpha }}_{_{m-1,m}})
\hskip 4mm si \hskip 4mm i=0\cr
(f^{^{\alpha }}_{_{01}},..,f^{^{\alpha }}_{_{i,i+1}}f^{^{\alpha
}}_{_{i-1,i}},..,
f^{^{\alpha }}_{_{m-1,m}})\hskip 4mm si \hskip 4mm0< i< m\cr (f^{^{\alpha
}}_{_{01}},
\dots,f^{^{\alpha }}_{_{m-2,m-1}})\hskip 4mm si \hskip 4mm i=m}}$}\par\vskip
1cm
\centerline{$\matrix{{\cal S}_{_{m,n}}&\Fhd{\varepsilon
^{^{1}}_{_{i}}}{}&{\cal S}_{_{m+1,n}}
\cr\cr(f^{^{\alpha }}_{_{01}},\dots,f^{^{\alpha }}_{_{m-1,m}})_{_{1\leq
\alpha \leq n}}&
\Fhd{}{}&\matrix{(f^{^{\alpha }}_{_{01}},..,I_{_{s(f^{^{\alpha
}}_{_{i,i+1}})}},..,
f^{^{\alpha }}_{_{m-1,m}})\hskip 4mm si \hskip 4mm 0\leq i< m\cr
(f^{^{\alpha }}_{_{01}},
\dots,f^{^{\alpha }}_{_{m-1,m}},I_{_{b(f^{^{\alpha }}_{_{m-1,m}})}})\hskip
4mm si
\hskip 4mm i=m}}$}\par\vskip 5mm
en consid\'erant la composition des fl\`eches ou la composition verticale
des 2-fl\`eches
suivant que $n$ est nul ou pas. On d\'efinit de m\^eme $d^{^{2}}_{_{i}}$ et
$\varepsilon ^{^{2}}_{_{i}}$ , en consid\'erant la composition horizontale
des 2-fl\`eches.
\par\hskip 5mm Le 2-pr\`e-nerf ${\cal S}$ est 1-troncable par construction,
en plus on a :\par\vskip 3mm
\centerline{$(T{\cal S})_{_{m}}=(T{\cal S})_{_{1}}{\times}_{_{(T{\cal
S})_{_{0}}}}
\dots{\times}_{_{(T{\cal S})_{_{0}}}}(T{\cal S})_{_{1}}$}\par\vskip 2mm
ce qui montre que ${\cal S}$ est 2-troncable. On sait que :\par\vskip 3mm
\centerline{$\matrix{{\cal S}_{_{m,n}}=\displaystyle\prod_{_{{\cal S}_{_{1,0}}
{\times}_{_{{\cal S}_{_{0,0}}}}\dots{\times}_{_{{\cal S}_{_{0,0}}}}{\cal
S}_{_{1,0}}}}^{n}
\bigr({\cal S}_{_{1,1}}{\times}_{_{{\cal
S}_{_{0,0}}}}\dots{\times}_{_{{\cal S}_{_{0,0}}}}
{\cal S}_{_{1,1}}\bigr)\cr\cr\cr \bigr({\cal S}_{_{1}}{\times}_{_{{\cal
S}_{_{0}}}}
\dots{\times}_{_{{\cal S}_{_{0}}}}{\cal S}_{_{1}}\bigr)(n)=\displaystyle
\prod_{_{{\cal S}_{_{0,0}}}}^{m}\bigr({\cal S}_{_{1,1}}{\times}_{_{{\cal
S}_{_{1,0}}}}
\dots{\times}_{_{{\cal S}_{_{1,0}}}}{\cal S}_{_{1,1}}\bigr)}$}\par\vskip 2mm
On d\'eduit, d'apr\`es la remarque (1.2.9), qu'on a un isomorphisme naturel
:\par
\centerline{$\diagram{{\cal S}_{_{m}}&\Fhd{\delta _{_{[m]}}}{}&{\cal S}_{_{1}}
{\times}_{_{{\cal S}_{_{0}}}}\dots{\times}_{_{{\cal S}_{_{0}}}}{\cal
S}_{_{1}}}$}\par
D'o\`u ${\cal S}$ est un 2-nerf strict. De m\^eme gr\^ace \`a la
remarque (1.2.9) et aux 1-\'equivalences \par\vskip 2mm
\centerline{$\delta _{_{[m]}}:\Phi _{_{m}}\fhd{}{}\Phi _{_{1}}
{\times}_{_{\Phi _{_{0}}}}\dots{\times}_{_{\Phi _{_{0}}}}\Phi _{_{1}}$ et
$\delta _{_{[m]}}:
\Psi \fhd{}{}\Psi _{_{1}}{\times}_{_{\Psi _{_{0}}}}\dots{\times}_{_{\Psi
_{_{0}}}}\Psi _{_{1}}$}
\par\vskip 2mm on peut construire les 2-\'equivalences ext\'erieures
$\alpha $ et $\beta $.\par	\vfill\eject


{\bf CHAPITRE 2 :}\par\vskip 1cm
\centerline{\bf NOTION de $n$-GROUPOIDE}\vskip 2cm\hskip 5mm
Dans ce chapitre on donne une d\'efinition de $n$-groupoide comme un
$n$-nerf ($n$-cat\'egorie large) dans lequel on peut inveser les $i$-fl\`eches
 \`a ($n$-$i$)-\'equivalence pr\`es.
D'autre part on fait associer \`a chaque espace topologique un
$n$-groupoide qui
g\'en\'eralise pour $n=1$ le groupoide de Poincar\'e.  \par\vskip 5mm\hskip 5mm


{\bf (2.1).--- Notations et d\'efinitions :}\par\vskip 5mm\hskip 5mm
Soit $\Phi $ un $n$-nerf avec $n\geq 1$, on d\'esigne par $\pi _{_{0}}(\Phi )$
l'ensemble des classes de $n$-\'equivalence d'objets de $\Phi $, il est
d\'efini comme suite :
\par\vskip 2mm\centerline{$\pi _{_{0}}(\Phi ) = T^{^{n}}\Phi = \{t^{n}(x) ,
x\in\Phi (0_{_{n}})
\} $}\par\vskip 2mm\hskip 5mm On appelle {\it $n$-groupoide} la donn\'ee d'un
$n$-nerf $\Phi $ telle que : \par pour tout $i\in\{1,\dots,n\}$ la cat\'egorie
${\cal C}_{_{i}}(\Phi ) = T^{^{n-i}}\Phi _{_{N}}$ avec $N=I_{_{(i-1)}}$ est
un groupoide.
Soit $F:\Phi \fhd{}{}\Psi $ un morphisme entre
$n$-nerfs, alors pour tout $i\in\{1,\dots,n\}$ , $F$  induit de fa\c con
naturelle un
morphisme :\par \vskip 2mm\centerline{${\cal C}_{_{i}}(F) : {\cal C}_{_{i}}
(\Phi )\Fhd{}{}{\cal C}_{_{i}}(\Psi )$}\par	\vskip 2mm et pour tout $i,j$
tels que
$1\leq i< j\leq n$ et tout objet $f$ de ${\cal C}_{_{i}}(\Phi )$ on a la
relation suivante :\par
\centerline{${\cal C}_{_{i}}(F)(0)(I^{^{i-j}}_{f})=I^{^{i-j}}_{{\cal
C}_{_{j}}(F)(0)(f)}$}\par
\vskip 2mm qui se d\'eduit facilement des d\'efinitions de $F$ et ${\cal
C}_{_{i}}(F)$.
Soient $\Phi $ un $n$-groupoide et $f$ un objet de ${\cal C}_{i}(\Phi )$.
On appelle
{\it $i$-\`eme groupe d'homotopie} de $\Phi $ de base $f$, le groupe
$Aut_{_{{\cal C}_{i}(\Phi )}}(f)$ qu'on notera par ${\pi }_{_{i}}(\Phi
,f)$. Pour tout objet
$f$ de $\Phi $ on d\'esignera par ${\pi }_{_{i}}(\Phi ,f)$ le groupe
${\pi }_{_{i}}(\Phi ,I^{^{i-j}}_{_{f}})$.\par\vskip 5mm


{\bf (2.2).--- Propriet\'es des $n$-groupoides :}
\par\vskip 5mm\hskip 5mm{\bf Th\'eor\`eme (2.2.1)} {\it Soit $\Phi $ un
$n$-groupoide,
alors pour tout $i,j$ tels que $2\leq j\leq i\leq n$, et tout objet $f$ de
${\cal C}_{_{j}}(\Phi )$, le groupe ${\pi }_{_{i}}(\Phi ,f)$ est
ab\'elien.}\par\vskip 5mm
\hskip 5mm{\bf Preuve :} Soient $i,j$ tels que $2\leq j\leq i\leq n$, et
posons :
$S = I_{_{i-2}}$ , $N = I_{_{i-1}}$. Le 2-nerf ${\cal A}_{_{i}} =
T^{^{n-i}}\Phi _{_{S}}$ est
un 2-groupoide tel que : \par\vskip 2mm
\centerline{2-fl${\cal A}_{_{i}}$ =$C_{_{i}}$= fl${\cal C}_{_{i}}(\Phi
)$\hskip 1cm et
\hskip 1cm fl${\cal A}_{_{i}}$ =$C_{_{i-1}}$= Ob${\cal C}_{_{i}}(\Phi
)$}\par\vskip 2mm
Donc la composition v\'erticale des 2-fl\`eches de ${\cal A}_{_{i}}$
coinside avec celle des
fl\`eches de ${\cal C}_{_{i}}$. Comme ${\cal A}_{_{i}}$ est une
2-cat\'egorie, alors elle
poss\'ede en plus une composition horisontale $\ast$ pour les 2-fl\`eches
et qui v\'erifie
avec la premi\`ere la relation de Godement suivante :\par\vskip
2mm\centerline{$
(\beta ^{'}\bar{\bullet}_{_{i}}\alpha ^{'})\ast(\beta
\bar{\bullet}_{_{i}}\alpha ) =
(\beta ^{'}\ast\beta )\bar{\bullet}_{_{i}}(\alpha ^{'}\ast\alpha )$}\vskip 2mm
o\`u  $\alpha $ , $\beta $ , $\alpha ^{'}$ et $\beta ^{'}$ sont des
2-fl\`eches de
${\cal A}_{_{i}}$ tels que $b(\alpha ) = s(\beta )$ , $b(\alpha ^{'}) =
s(\beta ^{'})$ et
$ss(\beta ^{'}) = bb(\alpha )$ , suivant le diagramme :\par
\centerline{$\diagram{&\Fhd{}{}&\Fhd{}{}\cr &\fvb{}{\alpha }&\fvb{}{{\alpha
}^{'}}\cr
&\Fhd{}{}&\Fhd{}{}\cr &\fvb{}{\beta }&\fvb{}{{\beta }^{'}}\cr
&\Fhd{}{}&\Fhd{}{}}$}
Lorsqu'on prend $\alpha $ et $\beta $ dans $\pi _{_{i}}(\Phi ,f)$ o\`u $f$
est un objet de
${\cal C}_{_{j}}(\Phi )$, la relation de Godement pour les diagrammes :\par
\centerline{$\diagram{&\Fhd{}{}&&\Fhd{}{}\cr &\fvb{}{\alpha
}&&\fvb{}{I^{^{i-j+1}}_{_{f}}}
\cr &\Fhd{}{}&&\Fhd{}{}&\cr &\fvb{}{I^{^{i-j+1}}_{_{f}}}&&\fvb{}{\beta }\cr
&\Fhd{}{}&&
\Fhd{}{}}$\hskip 2cm $\diagram{&\Fhd{}{}&&\Fhd{}{}\cr
&\fvb{}{I^{^{i-j+1}}_{_{f}}}&&
\fvb{}{\beta }\cr &\Fhd{}{}&&\Fhd{}{}&\cr &\fvb{}{\alpha
}&&\fvb{}{I^{^{i-j+1}}_{_{f}}}
\cr &\Fhd{}{}&&\Fhd{}{}}$} donne les relations : \par\vskip 2mm
\centerline{$\beta \ast\alpha  = (\beta \ast I)\bar
{\bullet}_{_{i}}(I\ast\alpha )=
(I\ast \alpha )\bar {\bullet}_{_{i}}(\beta \ast I)$\hskip 4mm o\`u\hskip 4mm
$I=I^{^{i-j+1}}_{_{f}}$}\par\vskip 2mm et d'apr\`es les axiomes (4), (5) et
(6) on obtient
les relations :\par\vskip 2mm\centerline{$\beta \ast\alpha =U^{-1}\bar
{\bullet}_{_{i}}
(\beta \bar{\bullet}_{_{i}}\alpha )\bar{\bullet}_{_{i}}U=U^{-1}\bar
{\bullet}_{_{i}}
(\alpha \bar{\bullet}_{_{i}}\beta )\bar{\bullet}_{_{i}}U$\hskip 4mm
o\`u\hskip 4mm
$U=U(I^{^{i-j}}_{_{f}})$}\par\vskip 2mm On en d\'eduit que :
\hskip 2mm $\alpha \bar{\bullet}_{_{i}}\beta  = \beta
\bar{\bullet}_{_{i}}\alpha $ ,
$\beta \ast\alpha =\alpha \ast\beta $ et par suite $\pi _{_{i}}(\Phi ,f)$
est ab\'elien.
\hfill


\par\vskip 5mm\hskip 5mm{\bf Proposition (2.2.2) :} {\it Soient $F : \Phi
\fhd{}{}\Psi $
une morphisme entre deux $n$-groupoides, $f$ un objet de
${\cal C}_{i}(\Phi )$ et $f^{'} ={\cal C}_{i}(F)(0)(f)$.
Alors $F$ induit un morphisme de groupes
${\pi }_{_{i}}(F,f) : {\pi }_{_{i}}(\Phi ,f)\ \Fhd{}{}\ {\pi }_{_{i}}(\Psi
,f^{'})$ pour tout
$i\in\{1,\dots,n\}$, et une application ${\pi }_{_{0}}(F) : {\pi
}_{_{0}}(\Phi ) \fhd{}{}
{\pi }_{_{0}}(\Psi )$ telle que ${\pi }_{_{0}}(F) = T^{^{n}}(F)$. }
\par\vskip 5mm\hskip 5mm{\bf Lemme (2.2.3) :} {\it Si $\lambda  : {\cal
A}\fhd{}{}{\cal B}$
est une morphisme entre deux cat\'egories, alors pour tout
$(x,y)\in {\cal A}(1){\times}_{_{{\cal A}(0)}}{\cal A}(1)$ on a :
$\lambda (1)(y.x) = \lambda (1)(y).\lambda (1)(x)$. }
\par\vskip 5mm\hskip 5mm{\bf Preuve du Lemme :} Dans la cat\'egorie $\Delta
$ les
fl\`eches $\delta _{ij} : [1]\fhd{}{}[2]$ nous donnent les diagrammes
commutatifs (1) et (2) :
\par\centerline{$\diagram{{\cal A}(1)&\Fhg{\delta ^{'}_{02}}{}&{\cal A}(2)&
\Fhd{\delta ^{'}_{01}\times\delta ^{'}_{12}}{}&
{\cal A}(1){\times}_{_{{\cal A}(0)}}{\cal A}(1)\cr \fvb{\lambda
(1)}{}&(1)&\fvb{\lambda (2)}{}
&(2)&\fvb{}{\lambda (1)\times \lambda (1)}\cr{\cal B}(1)&\Fhg{}{\delta
^{'}_{02}}&
{\cal B}(2)&\Fhd{}{\delta ^{'}_{01}\times\delta ^{'}_{12}}&{\cal
B}(1){\times}_{_{{\cal B}(0)}}
{\cal B}(1)}$}\par
Or, on sait que $\delta ^{'}_{01}\times\delta ^{'}_{12}$ est une bijection
donc le grand
diagramme est commutatifs ce qui montre que $\lambda $ respecte la
composition des
fl\`eches.\hfill \par\vskip 5mm\hskip 5mm{\bf Preuve de la proposition :}
On sait que pour tout $i\in\{1,\dots,n\}$, on a une morphisme
${\cal C}_{i}(F) :{\cal C}_{i}(\Phi )\fhd{}{}{\cal C}_{i}(\Psi )$, et par
suite un
 diagramme commutatif :\par
\centerline{$\diagram{{\cal C}_{i}(\Phi )(1)&\Fhd{{\cal
C}_{i}(F)(1)}{}&{\cal C}_{i}
(\Psi )(1)\cr \fvb{s,b}{}&(1)&\fvb{}{s,b}\cr
{\cal C}_{i}(\Phi )(0)&\Fhd{}{{\cal C}_{i}(F)(0)}&{\cal C}_{i}(\Psi )(0)}$}\par
Soit $\alpha $ est un \'el\'ement de ${\cal C}_{i}(\Phi )(1)$ tel que
$s(\alpha ) = b(\alpha )$ ;
alors d'apr\`es le diagramme (1) on obtient :\par
\centerline{$s\bigl[{\cal C}_{i}(F)(1)(\alpha )\bigr] = b\bigl[{\cal
C}_{i}(F)(1)(\alpha )
\bigr] ={\cal C}_{i}(F)(0)\bigr(s(\alpha )\bigr)$}\par\vskip 3mm
On en d\'eduit que pour tout objet $f$ de ${\cal C}_{i}(\Phi )$ la
restriction de
${\cal C}_{i}(F)(1)$ \`a ${\pi }_{_{i}}(\Phi ,f)$ est une application ${\pi
}_{_{i}}(F,f) :
{\pi }_{_{i}}(\Phi ,f)\ \fhd{}{}\ {\pi }_{_{i}}(\Psi ,f^{'})$. En
appliquant le Lemme (2.2.3)
au morphisme ${\cal C}_{i}(F)$, on montre que ${\pi }_{_{i}}(F,f)$
respecte la composition des fl\`eches de ${\cal C}_{i}(\Phi )$. Par
cons\'equent
${\pi }_{_{i}}(F,f)$ est un morphisme de groupes.\hfill


\par\vskip 5mm\hskip 5mm{\bf Proposition (2.2.4) :} {\it Soit $F : \Phi
\fhd{}{}\Psi $
un morphisme entre deux $n$-groupoides.
$F$  est une $n$-\'equivalence ext\'erieure si et seulement si pour tout
$i\in\{1,\dots,n\}$ et tout objet $f$ de ${\cal C}_{i}(\Phi )$; ${\pi
}_{_{i}}(F,f)$ est
un isomorphisme, et ${\pi }_{_{0}}(F)$ est une bijection.}\par
\vskip 5mm\hskip 5mm{\bf Preuve :} On sait que $F$ est une $n$-\'equivalence
ext\'erieure si et seulement si pour tout $i\in\{1,\dots,n\}$ et tout $u,
v$ objets de
${\cal C}_{i}(\Phi )$, l'application :	\par\vskip
2mm\par\centerline{$G^{u,v}_{i} :
Hom_{_{{\cal C}_{i}(\Phi )}}(u,v)\Fhd{}{}Hom_{_{{\cal C}_{i}(\Phi )}}(u,u) $}
\vskip 2mm induite par ${\cal C}_{i}(F)$ ainsi que l'application
$T^{^{N}}F:T^{^{N}}\Phi \fhd{}{}T^{^{N}}\Psi$ sont
 bijectives, o\`u $u^{'} = {\cal C}_{i}(F)(1)(u)$ (Proposition (1.3.1)).
Lorsqu'on prend $i\in\{1,\dots,n\}$ et $u=v$, on obtient
${\pi }_{_{i}}(F,u)=G^{u,u}_{i}$ et qui est donc un isomorphisme.\par
Inversement lorsque $ Hom_{_{{\cal C}_{i}(\Phi )}}(u,v)$ est non vide, la
donn\'ee d'un
\'el\'ement $\alpha $ de $ Hom_{_{{\cal C}_{i}(\Phi )}}(v,u)$ (il existe
toujours car
${\cal C}_{i}(\Phi )$ est un groupoide) permet de construire une bijection
:	\par
\centerline{$\diagram{Hom_{_{{\cal C}_{i}(\Phi )}}(u,v)
&\Fhd{t_{_{\alpha }}}{}&Hom_{_{{\cal C}_{i}(\Phi )}}(u,u) \cr \varphi &\Fhd{}{}
&\alpha \cdot\varphi }$}\par Si on pose
$\alpha ^{'}= G^{v,u}_{i}(\alpha )$ on obtient le diagramme commutative
suivant :
\par\centerline{$\diagram{Hom_{_{{\cal C}_{i}(\Phi )}}(u,v)
&\Fhd{t_{_{\alpha }}}{}&Hom_{_{{\cal C}_{i}(\Phi )}}(u,u)\cr
\fvb{G^{u,v}}{}&&\fvb{}{G^{u,u}_{i}}\cr Hom_{_{{\cal C}_{i}(\Phi
)}}(u^{'},v^{'})
&\Fhd{}{t_{_{\alpha ^{'}}}}&Hom_{_{{\cal C}_{i}(\Phi )}}(u^{'},u^{'})}$}\par
Et comme $G^{u,u}$ est un isomorphisme, alors $G^{u,v}$ est une bijection.
Par cons\'equent $F$ est une $n$-\'equivalence ext\'erieure.\hfill


\par\vskip 5mm\hskip 5mm{\bf Lemme (2.2.5) :} {\it Soit $\Phi $ un
($n$+1)-groupoide,
avec $n\geq 1$, alors :	\par (1) $\Phi _{1}{\times}_{_{\Phi _{0}}}\Phi
_{1}$ est un
$n$-groupoide;\par (2) Pour tout $i,j$ tels que $2\leq j\leq i\leq n$ et
tout objet $(f,g)$ de
${\cal C}_{j}(\Phi _{1}{\times}_{_{\Phi _{0}}}\Phi _{1})$ on a :\par\vskip 2mm
	\centerline{$\pi _{i}\bigr(\Phi _{1}{\times}_{_{\Phi _{0}}}\Phi
_{1},(f,g)\bigr)=
\pi _{i}(\Phi _{1},f){\times}_{_{\Phi (0_{{n}})}}\pi _{i}(\Phi
_{1},g)$}}\par\vskip 5mm
\hskip 5mm{\bf Preuve :} (1) Comme cons\'equence de la proposition (1.3.5),
$\Phi _{1}{\times}_{_{\Phi _{0}}}\Phi _{1}$ est un ($n$-1)-pr\'e-nerf
($n$-1)-troncable.\par Soient $s\in\{1,\dots,n-1\}$ , $(M,m)$ un objet de
$\Delta ^{^{n-s-1}}\times\Delta $ et posons $\Psi =\Phi
_{1}{\times}_{_{\Phi _{0}}}\Phi _{1}$.
\par	\vskip 2mm\centerline{On sait que :\hskip 4mm
$\Psi _{_{M,0}}=(\Phi _{1})_{_{M,0}}{\times}_{_{(\Phi _{0})_{_{M,0}}}}
(\Phi _{1})_{_{M,0}}=\Phi _{_{1,M,0}}{\times}_{_{\Phi _{_{0,M,0}}}}\Phi
_{_{1,M,0}}$}
\par	\vskip 2mm Or, $\Phi _{_{1,M,0}}$ et $\Phi _{_{0,M,0}}$ sont des
foncteurs constants,
alors il en est de m\^eme pour $\Psi _{_{M,0}}$.\par\vskip 2mm\hskip 5mm
Montrons maintenant que  $\delta ^{^{[m]}}_{_{M}}:\Psi _{_{M,m}}\fhd{}{}
\Psi _{_{M,1}}{\times}_{_{\Psi _{_{M,0}}}}\dots {\times}_{_{\Psi
_{_{M,0}}}}\Psi _{_{M,1}}$
\par est une $s$-\'equivalence ext\'erieure. Posons $N=(1,M)$:\par\vskip 2mm
\centerline{$\matrix{{\cal A}&=&\Psi _{_{M,1}}{\times}_{_{\Psi _{_{M,0}}}}\dots
{\times}_{_{\Psi _{_{M,0}}}}\Psi _{_{M,1}}\cr\cr {\cal B}&=&\bigr(\Phi
_{_{N,1}}
{\times}_{_{\Phi _{_{N,0}}}}\dots {\times}_{_{\Phi _{_{N,0}}}}\Phi
_{_{N,1}}\bigr)
{\times}_{_{\Phi _{_{0_{n-s}}}}}\bigr(\Phi _{_{N,1}}
{\times}_{_{\Phi _{_{N,0}}}}\dots {\times}_{_{\Phi _{_{N,0}}}}\Phi
_{_{N,1}}\bigr)}$}\par
\vskip 2mm On sait que :\par\vskip 2mm
	\centerline{$\matrix{\Psi _{_{M,m}}&=&\Phi _{_{N,m}}{\times}_{_{\Phi
_{_{0_{n-s}}}}}
\Phi _{_{N,m}}\cr\cr \delta ^{^{[m]}}_{_{N}}:\Phi _{_{N,m}}&\fhd{}{}&
\Phi _{_{N,1}}{\times}_{_{\Phi _{_{N,0}}}}\dots {\times}_{_{\Phi
_{_{N,0}}}}\Phi _{_{N,1}}
\hskip 3mm(s-equivalence)}$}\par\vskip 2mm D'autre part, d'apr\`es la
remarque (1.2.8),
pour tout objet $R$ de $\Delta ^{s}$ on a une bijection naturelle
$F(R):{\cal A}(R)\fhd{}{}{\cal B}(R)$ qui rend commutatif le diagramme
suivant :	\par
\centerline{$\diagram{\Psi _{_{M,m}}(R)&\Fhd{\delta
^{^{[m]}}_{_{M}}}{}&{\cal B}(R)\cr
\fvb{\delta ^{^{[m]}}_{_{N}}\times \delta ^{^{[m]}}_{_{N}}}{}&&\parallel\cr
{\cal A}(R)&
\Fhd{}{F(R)}&{\cal B}(R)}$}\par
L'isomorphisme naturel $F$ est alors une $s$-\'equivalence ext\'erieure
entre les
$s$-pr\'e-nerfs ${\cal A}$ et ${\cal B}$. D'o\`u $\delta ^{^{[m]}}_{_{M}}$
est une
$s$-\'equivalence ext\'erieure, et par cons\'equent $\Psi $ est un
$n$-nerf.\par
\hskip 5mm Soient $i\in\{1,\dots,n\}$ et $N=I_{_{i-1}}$, alors :\par\vskip 3mm
\centerline{${\cal C}_{i}(\Psi )=T{^{n-i}}\Psi _{_{N}}=
\bigr(T{^{n-i}}\Phi _{_{1,N}}\bigr){\times}_{_{\Phi
_{_{0_{n-1}}}}}\bigr(T{^{n-i}}
\Phi _{_{1,N}}\bigr)={\cal C}_{i}(\Phi _{1}){\times}_{_{\Phi _{_{0_{n-1}}}}}
{\cal C}_{i}(\Phi _{1})\hskip 4mm (\star)$}\par	\vskip 3mm
Comme ${\cal C}_{i}(\Phi _{1})$ est un groupoide, il s'en suit que ${\cal
C}_{i}(\Psi )$ est
aussi un groupoide, donc $\Psi $ est un $n$-groupoide.\par
(2) Soient $i,j$ tels que $2\leq j\leq j\leq n$ et $(f,g)$ un objet de
${\cal C}_{j}(\Psi )$,
alors d'apr\`es la relation $(\star)$ ona :\par\vskip 3mm
\centerline{$\matrix{\pi _{i}\bigr(\Psi ,(f,g)\bigr)&=&Aut_{_{{\cal
C}_{i}(\Psi )}}
\bigr((f,g)\bigr)\cr \cr&=&Aut_{_{{\cal C}_{i}(\Phi _{1}){\times}_{_{\Phi
_{_{0_{n-1}}}}}
{\cal C}_{i}(\Phi _{1})}}\bigr((f,g)\bigr)\cr \cr&=&Aut_{_{{\cal
C}_{i}(\Phi _{1})}}(f)
{\times}_{_{\Phi _{_{0_{n}}}}}Aut_{_{{\cal C}_{i}(\Phi _{1})}}(g)\cr\cr &=&
\pi _{i}(\Phi _{1},f){\times}_{_{\Phi (0_{{n}})}}\pi _{i}(\Phi
_{1},g)}$}\hfill


\par\vskip 5mm\hskip 5mm{\bf Th\'eor\`eme (2.2.6) :} {\it Soient $\Phi $ un
2-groupoide,
$\tau $ un \'el\'ement de $\Phi _{_{2,0}}$ tel que $\delta
^{^{0}}_{_{[2]}}(\tau )=(f,g)$ et
 $\delta ^{^{0}}_{_{02}}(\tau )=h$. Alors il existe un isomorphisme  \par
\vskip 2mm
\centerline{${\cal L} :\pi _{_{2}}(\Phi ,f)\ \Fhd{}{}\ \pi _{_{2}}(\Phi ,h)$}}
\par\vskip 5mm\hskip 5mm{\bf Preuve :} D'apr\`es les propositions (2.2.2) ,
(2.2,4) et le
Lemme (2.2.5), le diagramme : \par \vskip 4mm
\centerline{${\Phi _{_{1}}\Fhg{\delta _{_{02}}}{}\Phi _{_{2}}\Fhd{\delta
_{_{[2]}}}{}
\Phi _{_{1}}{\times}_{_{\Phi _{_{0}}}}\Phi _{_{1}}}$}\par\vskip 3mm o\`u
$\delta _{_{[2]}}$ est une 1-\'equivalence ext\'erieure, induit le morphisme et
l'isomorphisme suivants  :\par\centerline{$\diagram{\pi _{_{1}}(\Phi
_{_{1}},h)&
\Fhg{}{}&\pi _{_{1}}(\Phi _{_{2}},\tau )&\Fhd{}{\sim}&\pi _{_{1}}(\Phi
_{_{1}},f)
{\times}_{_{\Phi (0,0)}}\pi _{_{1}}(\Phi _{_{1}},g)}$}\par Or, on a ${\cal
C}_{_{2}}(\Phi )=
\Phi _{_{1}}={\cal C}_{_{1}}(\Phi _{_{1}})$ et ${\cal C}_{_{1}}(\Phi
_{_{2}})=\Phi _{_{2}}$
alors les morphismes pr\'ec\'edents coinsident avec :\par
\centerline{$\diagram{\pi _{_{2}}(\Phi ,h)&\Fhg{\delta ^{^{1}}_{_{02}}}{}&
Aut_{_{\Phi _{_{2}}}}(\tau )&\Fhd{\delta^{^{1}} _{_{[2]}}}{\sim}&
\pi _{_{2}}(\Phi ,f){\times}_{_{\Phi (0,0)}}\pi _{_{2}}(\Phi ,g)}$}\par
Consid\'erons maintenant le morphisme ${\cal L}$ compos\'e des morphismes :\par
\centerline{$\diagram{\pi _{_{2}}(\Phi ,f)&\Fhd{{\cal I}}{}&
\pi _{_{2}}(\Phi ,f){\times}_{_{\Phi (0,0)}}\pi _{_{2}}(\Phi ,g)&
\Fhg{\delta ^{^{1}}_{_{[2]}}}{\sim}&Aut_{_{\Phi _{_{2}}}}(\tau )&
\Fhd{\delta^{^{1}} _{_{02}}}{}&\pi _{_{2}}(\Phi ,h)}$}\par
o\`u ${\cal I}$ est le morphisme qui \`a $\alpha $ fait correspondre le couple
$(\alpha ,I_{_{g}})$. Pour montrer que ${\cal L}$ est un isomorphisme on va
l'expliciter.
Soient $\alpha $ un \'el\'ement de $\pi _{_{2}}(\Phi ,f)$ et $\sigma
=L_{_{2}}(f,g)$;
on a alors le diagramme commutatif suivant :\par
\centerline{$\diagram{\delta^{^{0}} _{_{[2]}}(\sigma
)&\Fhd{a_{_{2}}(f,g)}{\sim}&
\delta ^{^{0}}_{_{[2]}}(\tau )=(f,g)\cr \fvb{\delta ^{^{1}}_{_{[2]}}
(\varepsilon _{_{\alpha ,I_{_{g}}}})}{}&&\fvb{}{(\alpha ,I_{_{g}})}\cr
\delta ^{^{0}}_{_{[2]}}(\sigma )&\Fhd{\sim}{a_{_{2}}(f,g)}&
\delta ^{^{0}}_{_{[2]}}(\tau )=(f,g)}$}\par
Donc ils existent un isomorphisme $\mu :\sigma \fhd{}{}\tau $ et un morphisme
$\eta :\tau \fhd{}{}\tau $ uniques tels que $\delta ^{^{1}}_{_{[2]}}(\mu
)=a_{_{2}}(f,g)$ et
$\delta ^{^{1}}_{_{[2]}}(\eta )=(\alpha ,I_{_{g}})$. On en d\'eduit que :
$\eta \cdot\mu =\mu \cdot\varepsilon _{_{\alpha ,I_{_{g}}}}$, et par
cons\'equent :\par
\centerline{${\cal L}(\alpha )=\delta ^{^{1}}_{_{02}}(\eta )=\delta
^{^{1}}_{_{02}}(\mu )
\cdot (I_{_{g}}\star\alpha )\cdot[\delta ^{^{1}}_{_{02}}(\mu
)]^{^{-1}}$}\par\vskip 2mm
Montrons maintenant que ${\cal L}$ est bijective. \par
\underbar{(i) Injection :} Soit $(\alpha ,\beta )$ dans $\pi _{_{2}}(\Phi
,f)$ tel que
${\cal L}(\alpha )={\cal L}(\beta )$. Comme $\Phi $ est un 2-groupoide, on
suppose
la donn\'ee pour toute fl\`eche $k$ de $\Phi $ d'un isomorphisme de
coh\'erences
$\omega (k):k^{^{-1}}k\fhd{}{}I_{_{s(k)}}$ compatible avec les autres
isomorphismes de
coh\'erences de $\Phi $. Pour simplifier les notaions posons :
$\Omega =\omega (g)\star I_{_{f}}$ . On a les \'equivalences suivantes
:\par\vskip 4mm
\centerline{$\matrix{{\cal L}(\alpha )={\cal L}(\beta
)&\Leftrightarrow&I_{_{g}}\star
\alpha =I_{_{g}}\star\beta \cr &&\cr &\Leftrightarrow&
I_{_{g^{^{-1}}}}\star(I_{_{g}}\star\alpha
)=I_{_{g^{^{-1}}}}\star(I_{_{g}}\star\beta )
\cr &&\cr Axiomes\ (7) , (6) &\Leftrightarrow&
I_{_{g^{^{-1}}g}}\star\alpha =I_{_{g^{^{-1}}g}}\star\beta
\cr &&\cr &\Leftrightarrow&
\Omega \cdot(I_{_{g^{^{-1}}g}}\star\alpha)
\cdot \Omega ^{^{-1}} =\Omega \cdot(I_{_{g^{^{-1}}g}}
\star \beta )\cdot \Omega ^{^{-1}}
\cr &&\cr Axiome\ (3)&\Leftrightarrow&
I^{^{2}}_{_{s(g)}}\star\alpha =I^{^{2}}_{_{s(g)}}\star\beta
\cr &&\cr s(g)=b(f)&\Leftrightarrow&V(f)\cdot (I^{^{2}}_{_{b(f)}}\star\alpha)
(V(f))^{^{-1}}=V(f)\cdot (I^{^{2}}_{_{b(f)}}\star\alpha)(V(f))^{^{-1}}
\cr &&\cr Axiome\ (4) &\Leftrightarrow& \alpha =\beta }$}\vskip 4mm\par
\underbar{(ii) Surjection :} Soit $\lambda $ un \'el\'ement de $\pi
_{_{2}}(\Phi ,h)$ , on
va construire un l'\'el\'ement $\alpha $ de $\pi _{_{2}}(\Phi ,f)$ tel que
${\cal L}(\alpha )=\lambda $. Pour simplifier un peu les notations posons
:\par\vskip 2mm
\centerline{$\beta =[\delta ^{^{1}}_{_{02}}(\mu )]^{^{-1}}\cdot\lambda \cdot
\delta ^{^{1}}_{_{02}}(\mu )$  et  $A=A(f,g,g^{^{-1}})$  l'isomorphisme
d'associativit\'e.}\par
\vskip 2mm On alors les \'equivalences suivantes :\par\vskip 4mm
\centerline{$\matrix{{\cal L}(\alpha )=\lambda &\Leftrightarrow&I_{_{g}}\star
\alpha =\beta \cr &&\cr &\Leftrightarrow&
I_{_{g^{^{-1}}}}\star(I_{_{g}}\star\alpha )=I_{_{g^{^{-1}}}}\star\beta
\cr &&\cr Axiomes\ (7) , (6) &\Leftrightarrow&
I_{_{g^{^{-1}}g}}\star\alpha =A^{^{-1}}\cdot (I_{_{g^{^{-1}}}}
\star\beta )\cdot A\cr &&\cr  Axiome\ (3)&\Leftrightarrow&
I^{^{2}}_{_{s(g)}}\star\alpha=\Omega \cdot
A^{^{-1}}\cdot (I_{_{g^{^{-1}}}}\star\beta )\cdot A\cdot \Omega ^{^{-1}}
\cr &&\cr Axiome\ (4) &\Leftrightarrow& \alpha =V(f)\cdot \Omega \cdot
A^{^{-1}}\cdot (I_{_{g^{^{-1}}}}\star\beta )\cdot A\cdot \Omega ^{^{-1}}
[V(f)]^{^{-1}}}$}\vskip 4mm\par
On en d\'eduit que ${\cal L}$ est surjective, donc bijecive. Comme ${\cal
L}$ est compos\'e
de morphismes de groupes alors c'est un isomorphisme.\hfill
\vskip 5mm\hskip 5mm


\par\vskip 5mm\hskip 5mm{\bf Corollaire (2.2.7) :} {\it Soit $\Phi $ un
2-groupoide et
$f$ une fl\`eche de $\Phi $. Alors il existe un isomorphisme
${\cal L} :\pi _{_{2}}(\Phi ,I_{_{s(f)}}) \fhd{}{}\pi _{_{2}}(\Phi ,f)$}
\par\vskip 5mm\hskip 5mm{\bf Preuve :} Il suffit d'appliquer le
Th\'eor\`eme (2.2.5) \`a
l'\'el\'ement $\tau =\delta ^{^{0}}_{_{001}}(f)$ de $\Phi _{_{1,0}}$.\hfill


\vskip 5mm\hskip 5mm{\bf Proposition (2.2.8) :} {\it Soit $\Phi $ un
$n$-groupoide ,
$n\geq 2$. Alors pour tout entier $i$ tel que $2\leq i\leq n$ et
tout objet $f$ de ${\cal C}_{_{i}}(\Phi )$, on a un isomorphisme naturel :\par
$$\pi _{_{i}}(\Phi ,I_{_{s(f)}})\ \Fhd{}{\sim}\ \pi _{_{i}}(\Phi ,f)$$} \par
\hskip 5mm{\bf Preuve :} Consid\'erons le 2-nerf
${\cal A}_{_{i}}(\Phi )=T^{^{n-i}}\Phi _{_{I_{_{i-2}}}}$, c'est un
2-groupoide pour
lequel l'ensemble des fl\`eches coinside avec l'ensemble d'objets du 1-nerf
${\cal C}_{_{i}}(\Phi )$. Donc d'apr\`es le Corollaire (2.2.7) on a
l'isomorphisme :\par
$$\pi _{_{2}}({\cal A}_{_{i}}(\Phi ),I_{_{s(f)}})\ \Fhd{}{\sim}\
\pi _{_{2}}({\cal A}_{_{i}}(\Phi ),f)$$\vskip 2mm\par
Or, ${\cal C}_{_{2}}({\cal A}_{_{i}}(\Phi ))={\cal C}_{_{i}}(\Phi )$ alors
$\pi _{_{2}}({\cal A}_{_{i}}(\Phi ),I_{_{s(f)}})=\pi _{_{i}}(\Phi
,I_{_{s(f)}})$ et
$\pi _{_{2}}({\cal A}_{_{i}}(\Phi ),f)=\pi _{_{i}}(\Phi ,f)$, ce qui montre
la proposition.
\hfill \vskip 5mm\hskip 5mm
{\bf Remarque :} Si on consid\'ere $I_{_{b(f)}}$ dans le corollaire (2.2.7)
on obtient de
m\^eme un isomorphisme $\pi _{_{2}}(\Phi ,I_{_{b(f)}})\fhd{}{}\pi
_{_{2}}(\Phi ,f)$
et par suite un isomorphisme \par	\vskip 2mm\centerline{
$\varphi _{_{f}} : \pi _{_{2}}(\Phi ,I_{_{s(f)}})\Fhd{}{}\pi _{_{2}}(\Phi
,I_{_{b(f)}})$}
\par\vskip 2mm Lorsque $s(f)=b(f)=x$ l'automorphisme $\varphi _{_{f}}$ ne
depend que
de la classe de $f$ dans $\pi _{_{1}}(\Phi ,x)$, donc $\varphi $ est une
action de
$\pi _{_{1}}(\Phi ,x)$ sur $\pi _{_{2}}(\Phi ,x)$.
On est en train de g\'en\'eraliser la proposition (2.2.8) \`a $f$ un objet de
${\cal C}_{_{j}}(\Phi )$ avec $2\leq  j< i \leq n$. Ce qui entrainera
l'action de
$\pi _{_{1}}(\Phi ,x)$ sur les $\pi _{_{i}}(\Phi ,x)$.
\par\vskip 5mm\hskip 5mm


{\bf (2.3).--- $n$-Groupoide de Poincar\'e associ\'e \`a un espace topologique
:}
\par\vskip 2mm\hskip 5mm
{\bf Rappel (2.3.1) :} Soit $\Delta $ la cat\'egorie simpliciale, parmi ses
fl\`eches on
distingue la famille suivante qui jouera un r\^ole important dans la suite
: \par\hskip 3cm
$\diagram{[0]&\fhd{\delta _{i}}{}&[m]\cr 0&\fhd{}{}&i}$ \hskip 1cm
$\diagram{[1]&\fhd{\delta _{ij}}{}&[m]\cr 0, 1&\fhd{}{}&i, j}$
\par\vskip 4mm {\it Le foncteur $\cal R$ :} Pour tout $m$ entier positif on
d\'esigne par le
$m$-simplexe fondamental l'ensemble : \par
$R^{m}=\{ (t_0,\dots,t_n)$ tel que $0\leq t_i\leq 1$ et $\sum t_i =1\}$. Et
pour tout
$i\in\{0,\dots,m\}$, on d\'efinit les applications :\par
$\diagram{R^{m-1}\fhd{d^{"}_{i}}{}R^{m}\cr
(t_0,..,t_{m-1})\fhd{}{}(t_0,..,0,..,t_{_{m-1}})}$
\hskip 1cm $\diagram{R^{m+1}\fhd{\varepsilon ^{"}_{i}}{}
R^{m}\cr
(t_{0},..,t_{_{m+1}})\fhd{}{}(t_{0},..,t_{i}+t_{_{i+1}},..,t_{_{m+1}})}$
\par On vient de construire un foncteur covariant $\cal R$ : $\Delta
\fhd{}{}Ens$,
qui envoie les applications $\delta _{i}$ et $\delta _{ij}$ vers les
op\'erateurs c\^ot\'es
et sommets d'un $m$-simplexe suivants :\par
\hskip 5mm $\diagram{\{1\}\fhd{\delta _{i}^{"}}{}R^{m}\cr
1\fhd{}{}(0,..,0,1,0,..,0)}$
\hskip 2cm $\diagram{R\fhd{\delta _{ij}^{"}}{}R^{m}\cr
(t,s)\fhd{}{}(0,..,0,t,0,..,0,s,0,..,0)}$


\par\vskip 5mm {\bf Le multi-complexe singulier d'un espace topologique :}
\par\vskip 2mm\hskip 5mm
{\bf D\'efinition (2.3.2) :} On appelle {\it multi-complexe singulier}
(o\`u {\it $\infty$-complexe simgulier}) d'un espace topologique $X$ la famille
${\cal X} =(X^{n})_{n\geq 1}$ de $n$-pr\'e-nerfs $X^{n}$ avec $n\geq 1$,
d\'efinie de fa\c con r\'ecurrente pour tout $n\geq 1$ et tout $(M,m)$ objet de
$\Delta ^{n}\times \Delta $ par : \par\vskip 3mm
\centerline{$\matrix{X^{1}(m)=Hom(R^{m},X),\hskip 5mm X^{1}(0)=X
\cr\cr X^{n+1}(M,m) =\{f\in Hom(R^{m},X^{n}(M))\mid \forall x\in R^{m} ,
\forall i\in \{0,\dots,m_n\}\hskip 4mm  \delta ^{'}_i (f(x)) =
f_i\}}$}\par	\vskip 3mm
O\`u $\delta ^{'}_i : X^{n+1}_{_{M,m}}\fhd{}{}X^{n+1}_{_{M,0}}=X^{n}_{_{M}}$
est l' image de l'application $\delta _i$ par le foncteur $X^{n+1}_{_{M}}$
\par\vskip 2mm
et $f_i$ est un \'el\'ement de $X^{n}_{_{M}}$ ind\'ependant de $x$.
Dans la suite on notera les ensembles $X^{n}_{_{M}}$ par $X_{_{M}}$.\par
\vskip 5mm\hskip 5mm On d\'esigne par ${\cal R}_n : \Delta ^{n}\fhd{}{}Ens$ le
foncteur covariant d\'efini pour tout objet $(m_{_{1}},\dots,m_{_{n}})$ de
$\Delta ^{n}$,
pour tout $k\in \{1,\dots,n\}$ et tout $i\in\{0,\dots,m_k\}$ par :\par\vskip
2mm
\centerline{$\matrix{{\cal R}_{n}(m_{_{1}},\dots,m_{_{n}}) = R^{m_n}\times..
\times R^{m_1}\cr\cr{\cal R}_{n}(d^{k}_i) = I\times..d^{"}_i..\times I
	\hskip 5mm et
\hskip 5mm{\cal R}_{n}(\varepsilon ^{k}_i) = I\times..\varepsilon
^{"}_i..\times I}$}
\par\vskip 2mm o\`u les applications $d^{"}_i$ et $\varepsilon ^{"}_i$ sont
plac\'ees
\`a la $k$-\`eme position. On d\'efinit ${\cal H}_n$ par le foncteur
contravariant
$Hom({\cal R}_n, X)$, compos\'e du foncteur $Hom(\ , X)$ et ${\cal
R}_n$.\par\vskip 5mm
\hskip 5mm Soit $X$ un espace topologique. On d\'esigne par ${\cal L}$ la
famille des
$n$-pr\'e-nerfs ${\cal L}^{n}_{_{X}}$ avec $n\geq 1$ d\'efinie pour tout
$(M,m)$ objet de $\Delta ^{n}\times \Delta $ par :\par\vskip 3mm
\centerline{${\cal L}^{m}_{_{X}}=Hom(R^{m},X)\hskip 5mm et \hskip 5mm
{\cal L}^{n+1}_{_{X}}(M,m)=Hom(R^{m},{\cal L}^{n}_{_{X}}(M))$}\par
\vskip 3mm
{\bf Remarque :} $X^{n+1}_{_{X}}(M,m)$ est un sous ensemble de
${\cal L}^{n+1}_{_{X}}(M,m)$, plus p\'ecis\'ement $X^{n+1}_{_{X}}$ est un
sous foncteur
 de ${\cal L}^{n+1}_{_{X}}$.\vskip 3mm\hskip 5mm


{\bf Lemme (2.3.3) :} {\it Soit $\alpha :{\cal H}_{n}\fhd{}{}{\cal
L}^{n}_{_{X}}$
la morphisme d\'efinie pour tout $M=(m_{1},\dots,m_{n})$
objet de $\Delta ^{n}$ par :\par\centerline{$\diagram{{\cal H}_{n}(M)&
\Fhd{\alpha _{_{M}}}{}&{\cal L}^{n}_{_{X}}(M)\cr f&\Fhd{}{}&\alpha
_{_{M}}(f)}$}
tel que $\forall(x_n,\dots,x_1)\in {\cal R}_{n}(M)$ on a $\alpha
_{_{M}}(f)(x_n)
\dots(x_1)=f(x_n,\dots,x_1)$. Alors $\alpha $ est un isomorphisme naturelle.}
\par\vskip 5mm\hskip 5mm{\bf Preuve :} $\alpha $ admet un inverse naturel
$\beta $
d\'efini par : $\beta _{_{M}}(g)(x_n,\dots,x_1)=g(x_n)\dots(x_1)$.\hfill


\par\vskip 5mm\hskip 5mm{\bf Proposition (2.3.4) :} {\it Soit $M =
(m_{_{1}},\dots,m_{_{n}})
$ un objet de $\Delta ^{n}$. Alors tout \'el\'ement de $X_{_M}$ correspend
de fa\c con
biunivoque \`a un \'el\'ement $f$ de ${\cal H}_{n}(M)$ tel que :\par\vskip 2mm
\centerline{$\matrix{\forall k\in\{0,\dots,n-1\},\hskip 1cm
\forall i\in\{0,\dots,m_{_{n-k}}\}\cr\cr
 \hskip 5mm et \hskip 5mm \forall (x_n,..,{\hat x_{_{n-k}}},..,x_1)\in
{\cal R}_{n-1}(m_1,..,{\hat m_{_{n-k}}},..,m_n)\hskip 5mm on \hskip 2mm a\cr
\cr
f(x_n,..,x_{_{n-k+1}},\delta ^{"}_i, x_{_{n-k-1}},..,x_1) =
f_{i}(x_{_{n-k-1}},..,x_1)}$}\par
\vskip 3mm o\`u $f_i$ est un \'el\'ement de $X_{m_{_{1}},\dots,m_{_{n-k-1}}}$
ind\'ependant des variables $x_n,..,x_{_{n-k+1}}$.\par
Si  $m_{_{n-k}} = 0$  alors  $X_{_{M}} = X_{m_{_{1}},\dots,m_{_{n-k-1}}}$.}
\par\vskip 5mm\hskip 5mm
{\bf Preuve :} D'apr\`es la remarque p\'ec\'edente un \'el\'ement $f$ de
$X^{n+1}_{_{X}}(M,m)$ est en particulier un \'el\'ement de ${\cal
L}^{n+1}_{_{X}}(M,m)$
qui s'identifie par Lemme (2.3.3) \`a un
\'el\'ement $g$ de ${\cal H}_{n}(M)$ tel que
$g(x_n,\dots,x_1)=f(x_n)\dots(x_1)$ donc
$\forall k\in\{0,\dots,n-1\}$ et $\forall i\in\{0,\dots,m_{_{n-k}}\}$
:\par\vskip 2mm
\centerline{$\matrix{ f(x_n)..(x_{_{n-k+1}})(\delta ^{"}_i)(
x_{_{n-k-1}})..(x_1)&=&
f(x_n)..(\delta ^{"}_0)(\delta ^{"}_i)( x_{_{n-k-1}})..(x_1)\cr&\ldots&\cr &=&
f(\delta ^{"}_0)..(\delta ^{"}_0)(\delta ^{"}_i)(
x_{_{n-k-1}})..(x_1)}$}\vskip 3mm
ce qui montre que $g$ ne d\'epend que des variables $x_n,..,x_{_{n-k+1}}$.
Lorsque $m_{_{n-k}}=0$, la seule composante possible \`a la ($n$-$k$)-\`eme
coordonn\'ee  est $\delta ^{"}_{0}$ donc $f=f_{0}$.\hfill


\par\vskip 5mm {\bf Equivalence d'homotopie dans $X_{_{M}}$ :} \par\vskip 2mm
\hskip 5mm Soient $f$, $g$ deux \'el\'ements de $X_{_M}$ avec $M =
(m_1,\dots,m_s)$,
on dit que $f$ et $g$ sont homotopes et note $\overline f =\overline g$ si
et seulement si
il existe $\gamma \in X_{_{M,1}}$ tel que $\delta _{0}^{'}(\gamma) = f$ et
$\delta _{1}^{'}(\gamma) = g$. L'homotopie dans $X_{_M}$ est une relation
d'\'equivalence , et on d\'esignera par $\overline{(X_{_M})}$ l'ensemble de ses
classes d'\'equivalences. \par\vskip 5mm
\hskip 5mm{\bf Remarques :} D'apr\`es la proposition (2.3.4) on peut dire que
l'ensemble les classes d'homotopies $\overline{(X_{_M})}$ de $X_{_M}$
s'identifie \`a
$X_{m_{_{1}},\dots,m_{_{s}}}$ lorsqu'il existe un
$s\in\{1,\dots,m_{_{n}}\}$ tel que
$m_{_s} = 0$.\par\hskip 5mm Les applications
${\cal H}(d^{k}_i)$ et ${\cal H}(\varepsilon ^{k}_i)$
correspondent \`a des changements du domaine de la $k$-\`eme variable, donc
d'apr\`es la proposition, leur restrictions envoient les \'el\'ement de
$X_{_{M}}$
dans $X_{_{M^{'}}}$ o\`u $M$ et $M^{'}$ sont respectivement but et sourse des
$d^{k}_i$ et $\varepsilon ^{k}_i$.\par\vskip 5mm\hskip 5mm
Soit $X$ un espace topologique, on d\'esigne par $\Phi ^{n}$ : $\Delta
^{n}\fhd{}{}
Ens$ le $n$-pr\'e-nerf d\'efini pour tout objet $M$ de $\Delta ^{n}$ par :
\par\vskip 3mm
\centerline{$\Phi ^{n}(M) = \overline{(X_{_{M}})},\hskip 5mm\Phi
^{n}(d^{k}_i) =
\overline{X^{n}(d^{k}_i)}$\hskip 3mm et \hskip 3mm
$\Phi ^{n}(\varepsilon ^{k}_i) = \overline{X^{n}(\varepsilon ^{k}_i)}$}


\par\vskip 5mm\hskip 5mm{\bf Th\'eor\`eme (2.3.5) :} {\it Pour tout $n\geq 1$ ,
$\Phi ^{n}$ est un $n$-nerf et pour tout $k\in\{1,\dots,n-1\}$,
$T^{k}(\Phi ^{n}) = \Phi ^{n-k}$.}\par\vskip 5mm\hskip 5mm{\bf Preuve :}
Montrons d'abord que $\Phi ^{n}$ est $n$-tronquable. Il suffit pour cela de
montrer que
$\Phi $ est 1-tronquable et $T(\Phi ^{n}) = \Phi ^{n-1}$. Soit $(M,m)$ un objet
de
$\Delta ^{n-1}\times\Delta$, montrons que l'application suivante est bijective
:
\par\vskip 2mm$$\matrix{\overline{X_{_{M,m}}}&\Fhd{\delta ^{M}_{[m]}}{}&
\overline{X_{_{M,1}}}{\times}_{_{X_{M}}}\dots{\times}_{_{X_M}}
\overline{X_{_{M,1}}}\cr\cr \overline{f}&\Fhd{}{}&(\overline{\delta
^{'}_{_{01}}(f)},
\dots,\overline{\delta ^{'}_{_{m-1,m}}(f)})}$$\hskip 3mm
{\bf (1)} {\bf Injection :} Soient $f$ et $g$ dans $X_{_{M,m}}$ tels que
$\delta ^{^{M}}_{_{[m]}}(\overline{f})= \delta ^{^{M}}_{_{[m]}}(\overline{f})$,
alors pour tout $i\in\{0,\dots,n-1\}$ il existe une homotopie $\lambda
_{i,i+1}$ entre
$\delta ^{'}_{_{01}}(f)$ et $\delta ^{'}_{_{01}}(g)$. les applications
:\par\vskip 2mm
 \centerline{$f : \delta ^{"}_{0}(\{1\})\fhd{}{}X_{_{M}}\ $, $g : \delta
^{"}_{0}(\{1\})
\fhd{}{}X_{_{M}}\ $ et $\ \lambda _{i,i+1} : R\times\delta ^{"}_{i,i+1}(R)$}
\par\vskip 2mm v\'erifient la condition de recollement des applications
continues,
donc on obtient une apllication $F : (\delta ^{"}_{0}(\{1\})\times R^{m})\cup
(\delta ^{"}_{1}(\{1\})\times R^{m})\cup (\displaystyle\bigcup ^{m-1}_{i=0}
R\times
\delta _{i,i+1}(R))\fhd{}{}X_{_{M}}$.\par Le domaine de $F$ est une
d\'eformation de
$R\times R^{m}$, donc ce d\'ernier admet un retracte $r$ sur le premier.
L'application
$\lambda  = Fr$ est alors une homotopie entre $f$ et $g$. En effet pour
tout $x\in R^{m}$
on a :\par \vskip 3mm
\centerline{$\matrix{\delta ^{'}_{_{0}}(\lambda )(x) =
\lambda (\delta ^{"}_{0}(\{1\}),x) &=&F(\delta ^{"}_{0}(\{1\}),x) = f(x)\cr\cr
\delta ^{'}_{_{1}}(\lambda )(x) = \lambda (\delta ^{"}_{1}(\{1\}),x) &=&
F(\delta ^{"}_{1}(\{1\}),x) = g(x)}$}\par\vskip 3mm
En plus $\lambda $ est bien dans $X_{_{M,m,1}}$, car pour tout $t\in R$ et tout
$i\in\{0,\dots,m-1\}$ \par
$\delta ^{'}_{_{i}}(\lambda (t)) = \lambda (t,\delta ^{"}_{1}(\{1\}) =
F(t,\delta ^{"}_{1}(\{1\})
 = \lambda _{i,i+1}(t,\delta ^{"}_{1}(\{1\}) = C_i$ ne d\'epend pas de $t$.
\vskip 3mm\hskip 3mm
{\bf (2)} {\bf Surjection :} Soit $(\overline{\lambda _0},\dots,\overline
{\lambda _{m-1}})\in
{\overline{X_{_{M,1}}}{\times}_{X_{_{M}}},\dots,{\times}_{X_{_{M}}}
\overline{X_{_{M,1}}}}$. Les applications $\lambda _i : {\delta
}^{"}_{i,i+1}(R^{m})
\fhd{}{}X_{_{M}}$ o\`u $i\in\{0,\dots,m-1\}$ v\'erifient la condition de
recollement en
une application globale $\lambda _{_{[m]}} : \displaystyle\bigcup
_{i=0}^{m-1}{\delta }^{"}
_{i,i+1}(R)\fhd{}{}X_{_{M}}$. \par Le domaine de $\lambda _{_{[m]}}$ est la
r\'eunion de $m$
c\^ot\'es du simplexe $R^{m}$, alors il existe un retracte $r$ de $R^{m}$
sur cette
r\'eunion, et par cons\'equent le compos\'e $\lambda = \lambda _{_{[m]}}r$
est un
$m$-simplexes simguliers tel que : $\delta _{_{[m]}}(\overline{\lambda}) =
(\overline{\lambda _0},\dots,\overline{\lambda _{m-1}})$. Pour tout $t\in R$ et
tout $i\in\{0,\dots,m-1\}$, on a : \par\vskip 2mm
\centerline{$\delta ^{"}_{_{i,i+1}}(\lambda )(t) = \lambda
(\delta ^{"}_{i,i+1}(t)) = \delta _{_{[m]}}(\delta ^{"}_{i,i+1}(t)) =
\lambda _{i,i+1}(t)$}
\par	\vskip 3mm
De m\^eme on peut voir facilement que $\lambda \in X_{_{M,m}}$. \par\vskip 2mm
\hskip 5mm Montrons maintenant que $T(\Phi ^{n}) = \Phi ^{n-1}$.
Soient $M = (m_{_{1}},\dots,m_{_{n-1}})$ un objet de $\Delta ^{n-1}$ et $f$
, $g$ deux
\'el\'ement de $X_{_{M}}$. Alors $f$ et $g$ sont
$1$-\'equivalents si et seulement si il existe $\overline{\lambda }$ dans
$\overline{X_{_{M,1}}}$ tel que $\delta ^{'}_{_{0}}(\lambda ) = f$ et
$\delta ^{'}_{_{1}}
(\lambda ) = g$, ce qui est \'equivalent de dire que $f$ et $g$ sont
homotopes, par
cons\'equent on a : \par\vskip 2mm
\centerline{$T(\Phi ^{n})(M) = [\Phi ^{n}(M,0)]^{\sim} =
[\overline{X_{_{M,0}}}]^{\sim} =
[X_{_{M}}]^{\sim} = \overline{X_{_{M}}} = \Phi ^{n-1}(M)$}\par\vskip 2mm
Ce qu'on vient de montrer pour $\Phi ^{n}$ est aussi valable pour $\Phi
^{n-1}$. On en
d\'eduit alors que $\Phi ^{n}$ est $n$-tronquable.\par\vskip 5mm\hskip 5mm
{\bf Remarque :} La notion de $k$-\'equivalence int\'erieure dans le
$n$-pr\'e-nerf
$\Phi ^{n}$ est identique \`a celle de l'homotopie dans $\Phi ^{n-k}$, pour
tout
$k\in\{0,\dots,n-1\}$.\par\vskip 5mm\hskip 5mm
Il nous reste \`a montrer que pour tout $s\in\{0,\dots,n-1\}$ et tout objet
$(M,m)$
de $\Delta ^{n-s}\times\Delta$ le morphisme :\par
$$\delta ^{^{M}}_{_{[m]}} :
\Phi ^{n}_{_{M,m}}\Fhd{}{}\Phi ^{n}_{_{M,1}}{\times}_{_{\Phi
^{n}_{_{M,0}}}}\dots
{\times}_{_{\Phi ^{n}_{_{M,0}}}}\Phi ^{n}_{_{M,1}}$$ \par
est une (n-s-1)-\'equivalence ext\'erieure. Soit $h$ tel que $0\leq h <
n-s-1$, on pose :
\par\vskip 3mm
\centerline{$\matrix{n_s&=&n-s-1\cr\cr I&=&I_{_{n_{s}-h-1}}\cr\cr
 \delta ^{^{M}}_{_{[m]}}(n_{s})&=&\delta ^{^{M}}_{_{[m]}}(I)}$\hskip 1cm
$\matrix{R^{I}&=&R\times\dots\times R\cr\cr X_{_{M,m,I}}&=&\Phi ^{n}(M,m,I)
\cr\cr X_{_{M,m,I,1}}&=&\Phi ^{n}(M,m,I,1)}$}
\par \vskip 4mm Soient $u$ , $v$ , $w$ comme dans le diagramme suivant :
$$\diagram{{ u, v \in X_{_{M,m,I}}}&\Fhd{\delta ^{^{M}}_{_{[m]}}(n_{s})}{}&
\Pi_{_{X_{_{M}}}}X_{_{M,1,I}}\cr \fvh{\delta ^{'}_{0},\delta ^{'}_{1}}{}&&
\fvh{}{\delta ^{'}_{0},\delta ^{'}_{1}}\cr X_{_{M,m,I,1}}&\Fhd{}{\delta
^{^{M}}_{_{[m]}}
(n_{s}+1)}&{\Pi_{_{X_{_{M}}}}X_{_{M,1,I,1}}}\ni w}$$ Les applications
suivantes :\par
$$\diagram{u : \delta ^{"}_{0}(\{1\})\times R^{I}\times
R^{m}&\Fhd{}{}&X_{_{M}}\cr
v : \delta ^{"}_{1}(\{1\})\times R^{I}\times R^{m}&\Fhd{}{}&X_{_{M}}\cr
w_{i,i+1} : R\times R^{I}\times\delta ^{"}_{i,i+1}(R)&\Fhd{}{}&X_{_{M}}\cr }$$
v\'erifient la condition de recollement, donc on obtient une application
continue :\par
$$F : (\delta ^{"}_{0}(R)\times R^{m})\cup (\delta ^{"}_{1}(R)\times R^{m})\cup
(\displaystyle\bigcup ^{m-1}_{i=0} R\times \delta ^{"}_{i,i+1}(R))\
\Fhd{}{}\ X_{_{M}}$$\par
\hskip 5mm Soit $r$ est un retracte de $R\times R^{I}\times R^{m}$ sur le
domaine de $F$,
le compos\'e $x_r = Fr$ est un \'el\'ement de $X_{_{M,m,I,1}}$ dont l'image par
$\delta ^{^{M}}_{_{[m]}}(n_{_{s}}+1)$ est homotope \`a $w$. L'appartenance
de $x_r$ \`a
$X_{_{M,m,I,1}}$ est due \`a sa d\'ependance de $u$, $v$ et $w$. D'autre
part si
$(t,z,t^{'})$ est un \'el\'ement de $R\times R^{I}\times\delta
^{"}_{i,i+1}(R)$, alors : \par
$$\delta ^{'}_{i,i+1}(x_r)(t,z,t^{'}) = x_{r}(t,z,\delta ^{'}_{i,i+1}(t^{'})) =
F(t,z,\delta ^{'}_{i,i+1}(t^{'})) = w_{i,i+1}(t,z,t^{'})$$\par
ce qui montre que $\delta ^{^{M}}_{_{[m]}}(n_{s}+1)(x_{r}) = w$.
Soit $y$ dans $X_{_{M,m,I,1}}$ tel que : $\delta ^{'}_{0}(y) = u$ , $\delta
^{'}_{1}(y) = v$ et
pour tout $i\in\{0,\dots,m-1\}$, il existe une homotopie $\lambda _{i,i+1}$
entre
$\delta ^{'}_{i,i+1}(y)$ et $\delta ^{'}_{i,i+1}(x_r)$. De m\^eme en
recollant les applications :
\par $$\lambda _{i,i+1} :  R\times R\times R^{I}\times\delta
^{"}_{i,i+1}(R)\Fhd{}{}X_{_{M}}$$
\par on obtient une application globale \par
$$L : \cup (\displaystyle\bigcup ^{m-1}_{i=0}R\times
R\times R^{I}\times\delta ^{"}_{i,i+1}(R))\Fhd{}{}X_{_{M}}$$\par
tel que le compos\'e $\lambda  = Lr$ de $L$ avec un retracte $r$ de
$R\times R\times R^{I}\times R^{m}$ sur le domaine de $L$
soit une homotopie entre $y$ et $x_r$. Lorsque $h = n-s-1$, la construction
de $x_r$
se fait de la m\^eme fa\c con en consid\'erant un repr\'esentant de la
classe $w$,
et on se ram\`ene aux cas pr\'ec\'edants en utilisant la transitivit\'e de
l'homotopie.
\hfill \par\vskip 5mm\hskip 5mm{\bf Conclusion :} Pour chaque entier
naturel non nul $n$ et
chaque espace topologique $X$  nous venant de lui faire associer un $n$-nerf.


\par\vskip 5mm\hskip 5mm{\bf  Th\'eor\`eme (2.3.6) :} {\it Soit $X$ un
espace topologique,
alors pour tout entier naturel non nul $n$, le $n$-nerf $\Phi ^{n}_{X}$ est un
$n$-groupoide, qu'on appellera $n$-groupoide de Poincare de $X$ et qu'on
notera par
$\Pi _{n}(X)$.}\par\vskip 5mm\hskip 5mm{\bf Preuve :} On d\'esigne par
$\Phi $ le $n$-nerf $\Phi ^{n}_{_{X}}$, et consid\'erons la cat\'egorie
${\cal C}_{h}(\Phi )$. Il est claire qu'on a : ${\cal C}_{_{h-1}} = X_{N}$
et ${\cal C}_{_{h}}= \overline{X_{_{N,1}}}$. Soit l'application :\par
\centerline{$\diagram{\overline{X_{_{N,1}}}&\Fhd{Inv}{}&\overline{X_{_{N,1}}}\cr
\overline{f}&\Fhd{}{}&\overline{f^{'}}}$\hskip 1cm o\`u \hskip 1cm
$\diagram{R&\Fhd{f^{'}}{}&X_{_{N}}\cr (s,t)&\Fhd{}{}&f(t,s)}$}\par\hskip 5mm
Montorns que l'application $Inv$ est celle qui fait associer \`a chaque
fl\`eche de
${\cal C}_{h}(\Phi )$ son inverse par la loi de composition des fl\`eches,
pour cela
il faut montrer qu'on a la commutativit\'e des deux diagrammes suivants :\par

\centerline{$\diagram{\overline{X_{_{N,1}}}{\times}_{X_{N}}\overline{X_{_{N,
1}}}&
\Fhd{\mu }{}&\overline{X_{_{N,1}}}\cr \Fvb{I_{\overline{X_{_{N,1}}}}\times
Inv}{}&&
\Fvb{}{I}\cr \overline{X_{_{N,1}}}&\Fhd{}{s}&X_{_{N}}}$\hskip 2cm
$\diagram{\overline{X_{_{N,1}}}{\times}_{X_{N}}\overline{X_{_{N,1}}}&\Fhd{\mu
}{}&
\overline{X_{_{N,1}}}\cr \Fvb{ Inv\times I_{\overline{X_{_{N,1}}}}}{}&&
\Fvb{}{I}\cr
\overline{X_{_{N,1}}}&\Fhd{}{b}&X_{_{N}}}$} \par\hskip 5mm

On sait, d'apr\`es (2) de la d\'emonstration du Th\'eor\`eme (2.3.5),
que la donn\'ee d'un retracte $r$ de $R^{2}$ sur ses deux c\^ot\'es
$\delta ^{"}_{01}(R)\cup\delta ^{"}_{12}(R)$ permet la construction pour
tout $(f,g)$
dans $X_{_{N,1}}\times X_{_{N,1}}$ d'un repr\'esentant $l$ de la classe qui
correspond
au compos\'e de $\overline{f}$ et $\overline{g}$. Si on d\'esigne par $F$
l'application de
$\delta ^{"}_{01}(R)\cup\delta ^{"}_{12}(R)$ dans $X$ qui est d\'efinie
comme le
recollement de $f$ et $g$, alors on peut prendre $l = \delta
^{"}_{02}(F\circ r)$ suivant
le diagramme :\par\centerline{$\diagram{R&\Fhd{l}{}&X_{_{N}}\cr
\fvb{\delta ^{"}_{02}}{}&&\fvh{}{F}\cr R^{2}&\Fhd{}{r}&\delta ^{"}_{01}(R)
\cup\delta ^{"}_{12}(R)}$}\par
\hskip 5mm{\bf Constuction d'un retracte :} L'ensemble $R^{2}$ est une
surface plane
de l'espace euclidien r\'eel ayant $\vec n = (1,1,1)$ comme vecteur normale.
Soient les points $B = (0,1,0)$ et $D = ({1\over 2},0,{1\over 2})$ de
$R^{2}$. Pour tout
point $M = (x,y,z)$ de $R^{2}$ on consid\`ere le plan ${\cal P}_{_{M}}$ de
l'espace euclidien
passant par $M$ et de base $(\vec n ,\overrightarrow{BD})$, son \'equation est
$X-Z = -z+x$. On d\'efini le retracte $r$ par :
\par\centerline{$\diagram{R^{2}&\Fhd{r}{}&\delta ^{"}_{01}(R)\cup\delta
^{"}_{12}(R)\cr
 M&\Fhd{}{}&{\delta ^{"}_{12}(R)\cap {\cal P}_{_{M}}\  \  si \ 0\leq x\leq
{1-y\over 2}}\cr
&&{\delta ^{"}_{01}(R)\cap {\cal P}_{_{M}}\  \  si \ {1-y\over 2}\leq x\leq
1-y}}$}\par
Qui est donn\'e explicitement par :\par
\centerline{$\diagram{R^{2}&\Fhd{r}{}&\delta ^{"}_{01}(R)\cup\delta
^{"}_{12}(R)\cr
 M&\Fhd{}{}&{(0,1+x-z,-x+z)\  \  si \ 0\leq x\leq {1-y\over 2}}\cr
&&{(x-z,1-x+z,0)\  \  si \ {1-y\over 2}\leq x\leq 1-y}}$}
Montrons maintenant qu'on a :\par\vskip 2mm
\centerline{(1) $\mu (\overline{f},Inv(\overline{f})) =
\overline{Id_{_{s(f)}}}$\hskip 5mm
et \hskip 5mm (2) $\mu (\overline{g},Inv(\overline{g})) =
\overline{Id_{_{b(g)}}}$}
\par\vskip 2mm \hskip 5mm(1) Soit $l_{1}$ le repr\'esentant correspondant
\`a la classe
du compos\'e de $\overline{f}$ et $Inv(\overline{f})$ d\'efini par :\par
\centerline{$\diagram{R&\Fhd{l_{1}}{}&X_{_{N}}\cr
 (x,z)&\Fhd{}{}&f(1-2x,2x)\  \  si\  0\leq x\leq {1\over 2}\cr
&&f(2x-1,2-2x)\  \  si\  {1\over 2}\leq x\leq 1}$}\par
L'homotopie $\Gamma _{1} : R\times R\fhd{}{}X_{_{N}}$ d\'efinie par\par
\vskip 4mm\centerline{$\Gamma _{1}\bigl[(a,b),(x,y)\bigr] = \left\{\matrix{
f(1,0)&&si\ 0\leq x\leq {b\over 2}\cr\cr
f(1-2x+b,2x-b)&&si\  {b\over 2}\leq x\leq {1\over 2}\cr\cr
f(2x+b-1,2-2x-b)&&si\ {1\over 2}\leq x\leq 1-{b\over 2}\cr\cr
f(1,0)&&si\  1-{b\over 2}\leq x\leq 1}\right.$}\par\vskip 4mm
est tels que \par\vskip 3mm
\centerline{$\matrix{\delta ^{'}_{0}(\Gamma _{1})(x,y)&=&
\Gamma _{1}\bigl[(1,0),(x,y)\bigr] = l_{1(x,y)}\cr\cr
\delta ^{'}_{1}(\Gamma _{1})(x,y)&=&\Gamma _{1}\bigl[(0,1),(x,y)\bigr] =
f(1,0) =
Id_{_{s(f)}}(x,y)}$}\par\vskip 3mm
Il est claire de voir que $\Gamma _{1}$ est dans $X_{_{N,1}}$, et par
cons\'equent $\overline{l_{1}} = \overline{Id_{_{s(f)}}}$.\par
\hskip 5mm (2) Une homotopie du m\^eme genre montrera qu'on a aussi
$\mu (\overline{g},Inv(\overline{g})) = \overline{Id_{_{b(g)}}}$. On en
d\'eduit alors que
$\Phi ^{n}_{X}$ est un $n$-groupoide.\par


\vskip 5mm{\bf  Fonctorialit\'e de $\Pi _{_{n}}$ :}\vskip 2mm\hskip 5mm
Soit $n$ un entier naturelle non nul fix\'e et soit $f : X\fhd{}{}Y$ une
fl\`eche dans
la cat\'egorie ${\cal T}op$. L'application $f$ induit un morphisme \par
$\Pi _{_{n}}(f) : \Pi _{_{n}}(X)\fhd{}{}\Pi _{_{n}}(Y)$ d\'efinie comme
suite :\par
\centerline{$\diagram{\Pi _{_{n}}(X)(M)&\Fhd{\Pi _{_{n}}(f)(M)}{}&\Pi
_{_{n}}(Y)(M)\cr
\overline{\lambda }&\Fhd{}{}&\overline{f\circ\lambda }}\hskip 2cm\diagram{
R^{^{M}}&\fhd{\lambda }{}&X\cr \fvb{Id}{}&&\fvb{}{f}\cr
R^{^{M}}&\fhd{}{f\circ\lambda }&
Y}$}\par\centerline{avec \hskip 1cm$M = (m_1,\dots,m_n)$ \hskip 5mm et
\hskip 5mm
$R^{^{M}} = R^{m_{n}}\times\dots\times R^{m_{1}}$}\par \vskip 3mm
La d\'efinition de $\Pi _{_{n}}(f)$ est bien justifi\'ee car pour tout
$(x_1,\dots,x_n)$ dans
$R^{^{M}}$ on a : \par\vskip 3mm
\centerline{$\matrix{f\circ\lambda (x_{_{n}},..,x_{_{k+1}},\delta
^{"}_{i},x_{_{k}},..,x_1)
&=&f\bigl[\lambda (x_{_{n}},..,x_{_{k+1}},\delta
^{"}_{i},x_{_{k}},..,x_1)\bigr]\cr \cr
&=&f\bigl[\lambda _{i}(x_{_{k}},..,x_1)\bigr]\hskip 5mm 0\leq i\leq
m_{_{k}}}$}\par
\vskip 3mm ce qui montre que $f\circ\lambda $ appartient \`a $Y_{_{M}}$. Soit
$\sigma  : M\fhd{}{}N$ une fl\`eche de $\Delta ^{n}$, on a \par\vskip 4mm
\centerline{$R^{^{M}}\Fhd{{\cal R}(\sigma )}{}R^{^{N}}\Fhd{\lambda
}{}X\Fhd{f}{}Y$}\par
\vskip 3mm Montrons qu'on a la commutativit\'e du diagramme suivant
:\par\vskip 2mm
\centerline{$\diagram{\Pi _{_{n}}(X)(M)&\Fhd{\Pi _{_{n}}(f)(M)}{}&\Pi
_{_{n}}(Y)(M)\cr
\fvh{\Pi _{_{n}}(X)(\sigma )}{}&&\fvh{}{\Pi _{_{n}}(Y)(\sigma )}\cr
\Pi _{_{n}}(X)(N)&\Fhd{}{\Pi _{_{n}}(f)(N)}&\Pi _{_{n}}(Y)(N)}$}
\par\vskip 2mm Soit $\overline{\lambda }\in \Pi _{_{n}}(X)(N)$ \par\vskip 2mm
\centerline{$\matrix{\Pi _{_{n}}(f)(M).\Pi _{_{n}}(X)(\sigma
)(\overline{\lambda })&=&
\Pi _{_{n}}(f)(M)\bigl(\overline{\lambda \circ {\cal R}(\sigma )}\bigr)\cr\cr
&=&\overline{f\circ\bigl(\lambda \circ {\cal R}(\sigma )\bigr)}\cr\cr
(associtivite\  de\  \circ)&=&
\overline{\bigl(f\circ\lambda \bigr)\circ {\cal R}(\sigma )}\cr\cr
&=&\Pi _{_{n}}(Y)(\sigma )(\overline{f\circ\lambda })\cr\cr
&=&\Pi _{_{n}}(Y)(\sigma ).\Pi _{_{n}}(f)(N)(\overline{\lambda
})}$}\par\vskip 3mm
On en d\'eduit qu'on a un foncteur covariant $\Pi _{_{n}}(\ )$ de la
cat\'egorie ${\cal T}op$
vers la cat\'egorie $n$-Gr , des $n$-groupoides et transformations
naturelles, qui \`a
un espace topologique $X$ fait correspondre son $n$-groupoide de Poincar\'e
$\Pi _{_{n}}(X)$.\par \vskip 5mm


{\bf (2.4).--- Groupes d'homotopie de $\Pi _{_{n}}(X)$}
\vskip 5mm\hskip 5mm
{\bf Proposition (2.4.1) :} {\it Soient $(Z,z_{_{0}})$ un espace
topologique, $n\geq 1$ et
$N = I_{_{n}}$. Alors ils existent des bijections entre les ensembles suivants
:
\par\vskip 2mm(a) $\bigl[(R^{^{N}},\delta R^{^{N}}) ; (Z,z_{_{0}})\bigr]\
\fhd{\alpha }{\sim}
\ \bigl[(R^{^{N}}/_{\delta R^{^{N}}},w_{_{0}}) ; (Z,z_{_{0}})\bigr]$\hskip 5mm
$w_{_{0}} = p(x)$ pour  $x\in\delta R^{^{N}}$\par\vskip 2mm
(b) $\bigl[(R^{^{N}}/_{\delta R^{^{N}}},w_{_{0}}) ; (Z,z_{_{0}})\bigr]\
\fhd{\beta }{\sim}\
\bigl[(S^{^{n}}, s_{_{0}}) ; (Z,z_{_{0}})\bigr]$\hskip 5mm
$s_{_{0}} = (1,0,\dots,0)$}\par
\vskip 5mm\hskip 5mm{\bf Preuve :}
(a) On d\'efinit l'application $\alpha $ comme par : $\alpha (\bar u)  =
\bar {u^{'}}$ o\`u
$u^{'} : R^{^{N}}/_{\delta R^{^{N}}}\fhd{}{}Z$ est l'unique application qui
rend
commutatif le diagramme suivant
:\par\centerline{$\diagram{R^{^{N}}&\Fhd{u}{}&Z\cr
\fvb{p}{}&&\parallel\cr R^{^{N}}/_{\delta R^{^{N}}}&\Fhd{}{u^{'}}&Z}$}\par
La continuit\'e de $u$ et de $p$ entrainent celle de $u^{'}$. Soit $\Gamma
$ une homotopie
relative entre deux \'el\'ements $u$ et $v$ de
$\bigl[(R^{^{N}},\delta R^{^{N}}) ; (Z,z_{_{0}})\bigr]$ alors elle induit
une homotopie
${\Gamma }^{'}$ qui rend commutatif le diagramme :\par
\centerline{$\diagram{I\times R^{^{N}}&\Fhd{\Gamma }{}&Z\cr
\fvb{Id\times p}{}&&\parallel\cr I\times{R^{^{N}}/_{\delta
R^{^{N}}}}&\Fhd{}{{\Gamma }^{'}}
&Z}$}\par
De m\^eme la continuit\'e de $\Gamma $ , $Id\times p$ et la relation
$I\times{R^{^{N}}/_{\delta R^{^{N}}}} = I\times R^{^{N}}/_{I\times \delta
R^{^{N}}}$
entrainent celle de ${\Gamma }^{'}$. Montrons maintenant que $\alpha $ est
une bijection.
\par{\bf Surjection :} Soit $\bar {u^{'}}$ dans $\bigl[(R^{^{N}}/
_{\delta R^{^{N}}},w_{_{0}}) ; (Z,z_{_{0}})\bigr]$. Le compos\'e $u =
u^{'}\circ p$ est une
application continue telle que $\alpha (\bar u) = \bar {u^{'}}$, donc
$\alpha $ est surjective.
\par{\bf Injection :} Soient $\bar u$ et $\bar v$ deux \'el\'ement de
$\bigl[(R^{^{N}},\delta R^{^{N}}) ; (Z,z_{_{0}})\bigr]$ tels que
$\alpha (\bar u) = \alpha (\bar v)$.  Soient $\bar {u^{'}} = \alpha (\bar u)$ ,
$\bar {v^{'}} = \alpha (\bar v)$ et $\Gamma ^{'}$ une homotopie relative
entre $u^{'}$
et $v^{'}$ alors $\Gamma  = {\Gamma }^{'} \circ (I\times p)$ est une
homotopie relative
entre $u$ et $v$, ce qui montre que $\bar u = \bar v$.
Finalement $\alpha $ est une bijection. \par\vskip 5mm\hskip 5mm
(b) La donn\'ee d'un hom\'eomorphisme
$\varphi : (S^{^{n}}, s_{_{0}})\fhd{}{}\ (R^{^{N}}/_{\delta
R^{^{N}}},w_{_{0}})$ nous
 permet de construire une bijection $\beta $ d\'efinie de la fa\c con
suivante :\par
\vskip 2mm\centerline{$\diagram{\bigl[(R^{^{N}}/_{\delta
R^{^{N}}},w_{_{0}}) ; (Z,z_{_{0}})
\bigr]&\Fhd{\beta }{}&\bigl[(S^{^{n}}, s_{_{0}}) ; (Z,z_{_{0}})\bigr]\cr
\bar f &\Fhd{}{}&\overline{f\circ \varphi }}$}\vskip 2mm\par
Si $\Gamma $ est une homotopie entre deux \'elements $f$ et $g$ de
$\bigl[(R^{^{N}}/_{\delta R^{^{N}}},w_{_{0}}) ; (Z,z_{_{0}})\bigr]$ alors
$\Gamma \circ (Id\times \varphi )$ est une homotopie entre $f\circ \varphi $ et
$g\circ \varphi $. Ce qui montre que $\beta (\bar f) = \beta (\bar g)$.
D'autre part
l'application $\beta $ admet un inverse $\beta ^{'}$ d\'efini par :
\par\vskip 2mm
$\beta ^{'}(\bar u) = \overline{u\circ {\varphi }^{^{-1}}}$, donc $\beta $
est une bijection.\hfill


\par\vskip 5mm\hskip 5mm
{\bf Proposition (2.4.2) :} {\it Soient $(Z,z_{_{0}})$ un espace
topologique, $n\geq 1$ et
$N = I_{_{(i-1)}}$ avec $1\leq i\leq n$. Alors on a :\par\vskip 2mm
(a)   $\bigl[(R,\delta R) ; (Z_{_{N}},Id_{_{N}}(z_{_{0}}))\bigr]$ est en
bijection avec $
\bigl[(R^{^{N,1}},\delta R^{^{N,1}}) ; (Z,z_{_{0}})\bigr]$ \par\vskip 2mm
(b)   $\pi _{_{i}}\bigl(\Pi _{_{n}}(Z), z_{_{0}}\bigr) =
\bigl[(R^{^{N,1}},\delta R^{^{N,1}}) ; (Z,z_{_{0}})\bigr]$}\par
\vskip 5mm\hskip 5mm{\bf Preuve :}
(a) On v\'erifiera sans peine que l'application $\sigma $ d\'efinie comme
suite :\par
\centerline{$\diagram{\bigl[(R,\delta R) ;
(Z_{_{N}},Id_{_{N}}(z_{_{0}}))\bigr]&
\fhd{\sigma }{}&\bigl[(R^{^{N,1}},\delta R^{^{N,1}}) ; (Z,z_{_{0}})\bigr]\cr
\bar f&\fhd{}{}&\bar F}$}
 o\`u $F$ est d\'efinie par $F(t_{_{1}},..,t_{_{i}}) =
f(t_{_{1}})...(t_{_{i}})$ est une bijection.
 \par\vskip 5mm
(b) On sait que : \par\centerline{
$\pi _{_{i}}\bigl(\Pi _{_{n}}(Z), z_{_{0}}\bigr) = \Bigl\{\bar f
\in\overline{Z_{_{N,1}}}\mid
\forall (t,x) \in (\delta R)\times R^{^{N}}\  f(t,x) =
Id_{_{N}}(z_{_{0}})\in Z_{_{N}}\Bigr\}$}
\par\vskip 2mm
\centerline{ Or,\hskip 1cm $\delta (R^{^{N,1}}) = (\delta R)\times
R^{^{N}}\bigcup R
\times(\delta R^{^{N}})$} \par \vskip 2mm
\centerline{et\hskip 1cm $\forall (t,x)\in R\times(\delta R^{^{N}})$ \hskip
5mm on a
\hskip 5mm $f(t,x) = f(\delta ^{'}_{_{0}},x) = f(\delta ^{'}_{_{1}},x) =
Id_{_{N}}(z_{_{0}})$}
\par\vskip 2mm\centerline{d'o\`u \hskip 1cm $\forall (t,x)\in \delta
R^{^{N,1}}$
\hskip 5mm on a \hskip 5mm $f(t,x) = Id_{_{N}}(z_{_{0}})$.}\par \vskip 2mm
On en d\'eduit alors que $\pi _{_{i}}\bigl(\Pi _{_{n}}(Z), z_{_{0}}\bigr) =
\bigl[(R^{^{N,1}},\delta R^{^{N,1}}) ; (Z,z_{_{0}})\bigr]$\hfill \par
\vskip 5mm\hskip 5mm Consid\'erant l'application suivante :\par
\centerline{$\diagram{\pi _{_{i}}\bigl(\Pi
_{_{n}}(X),x_{_{0}}\bigr)&\Fhd{\gamma }{}&
\pi _{_{1}}\bigl(X_{_{N}},Id_{_{N}}(x_{_{0}})\bigr)\cr \bar f&\Fhd{}{}&\bar
g}$}\par
\centerline{o\`u \hskip 2mm $g\ :\ \bigl(I,\delta I\bigr)\ \fhd{}{}\
\bigl(X_{_{N}},
Id_{_{N}}(x_{_{0}})\bigr)$\hskip 2mm d\'efinie par \hskip 2mm $g(t) =
f(1-t,t)$}\par
$\gamma $ est bien d\'efinie. En effet soient $\bar f$ et $\bar {f^{'}}$
sont deux
\'el\'e\'ements de $\pi _{_{i}}\bigl(\Pi _{_{n}}(X),x_{_{0}}\bigr)$ tels
que $\bar f =
\bar {f^{'}}$. Soit $\Gamma $ une homotopie entre $f$ et $f^{'}$ telle que
$\Gamma (\delta ^{"}_{_{0}},x) = f(x)$ , $\Gamma (\delta ^{"}_{_{1}},x) =
f^{'}(x)$ \par et
$\Gamma (x,\delta ^{"}_{_{i}}) = f(\delta ^{"}_{_{0}}) = f^{'}(\delta
^{"}_{_{0}}) =
Id_{_{N}}(x_{_{0}})$.\par Alors on peut d\'efinir de fa\c con naturelle une
homotopie entre
$g$ et $g^{'}$ par\par\centerline{$\diagram{I\times I&\Fhd{\Gamma
^{'}}{}&X_{_{N}}
\cr (s,t)&\Fhd{}{}&\Gamma \bigl((1-s,s), (1-t,t)\bigr)}$}\par
La continuit\'e de $\Gamma $ entraine celle de $\Gamma ^{'}$ ce qui montre que
$\bar g = \bar {g^{'}}$ et par cons\'equent $\gamma $ est bien d\'efinie.	\par


\vskip 5mm\hskip 5mm{\bf Proposition (2.4.3) :} {\it L'application $\gamma
$ est un
isomorphisme de groupes.}\par\vskip 5mm\hskip 5mm{\bf Preuve :} Montrons que
$\gamma $ est une bijection. \par (a) Soient $\bar f$ et $\bar {f^{'}}$
deux \'el\'ements de
$\pi _{_{i}}\bigl(\Pi _{_{n}}(X),x_{_{0}}\bigr)$ tels que $\gamma (f) =
\gamma ( f^{'})$.
Soit $\Gamma ^{'}$ une homotopie entre $g$ et $g^{'}$ o\`u $\bar g =\gamma
(f)$ et
$\bar {g^{'}} =\gamma ^{'}(f^{'})$ telle que $\Gamma ^{'}(0,x) = g(x)$ ,
$\Gamma ^{'}(1,x) =
g^{'}(x)$ et $\Gamma ^{'}(x,i) = g(i) = g^{'}(i) = Id_{_{N}}(x_{_{0}})$
pour $i = 0, 1$.\par
Alors on peut d\'efinir de fa\c con naturelle une homotopie $\Gamma $ entre
$f$ et
$f^{'}$ par\par\centerline{$\diagram{R\times R&\Fhd{\Gamma }{}&X_{_{N}}\cr
\bigl((s,s^{'}), (t,t^{'})\bigr)&\Fhd{}{}&\Gamma ^{'}(s^{'},t^{'})}$}\par
La continuit\'e de $\Gamma ^{'}$ entraine celle de $\Gamma $ ce qui montre que
$\bar g = \bar {g^{'}}$ et par cons\'equent $\gamma $ est injective.\par
(b) Soit $g$ un \'el\'ement de $\pi
_{_{i}}\bigl(X_{_{N}},Id_{_{N}}(x_{_{0}})\bigr)$, donc
l'application continue suivante :\par
\centerline{$\diagram{(R,\delta R)&\Fhd{f}{}&(X_{_{N}},Id_{_{N}}(x_{_{0}})\cr
(s,t)&\Fhd{}{}&g(t)}$}\par
et telle que $\gamma (\bar f) = \bar g$, ce qui entraine que $\gamma $ est
surjective.\par
Montrons maintenant que $\gamma $ est un morphisme de groupes, ce qui revient
\`a montrer que pour tout $\bar f$ , $\bar {f^{'}}$ dans
$\pi _{_{i}}\bigl(\Pi _{_{n}}(X),x_{_{0}}\bigr)$ on a : $\gamma (\bar
{f^{'}}\bullet\bar f) =
\bar {g^{'}}\cdot\bar g$.\par On sait que les compositions $\bullet$ et
$	\cdot$ sont
d\'efinies par $\bar {f^{'}}\bullet\bar f = \bar F$ et $\bar
{g^{'}}\cdot\bar g = \bar G$
o\`u \par\centerline{$\diagram{(R,\delta
R)&\Fhd{F}{}&(X_{_{N}},Id_{_{N}}(x_{_{0}}))}$
\hskip 5mm et\hskip 5mm$\diagram{(I,\delta
I)&\Fhd{G}{}&(X_{_{N}},Id_{_{N}}(x_{_{0}}))}$}
\par sont d\'efinies par :\par \centerline{$F(x,y) =
\left\{\matrix{f^{'}(2x,1-2x)&&si\
0\leq x\leq {1\over 2}\cr\cr f(2x-1,2-2x)&&si\  {1\over 2}\leq x\leq
1}\right.$}
\par\vskip 3mm\centerline{$G(t) = \left\{\matrix{g(2t)&&si\
0\leq t\leq {1\over 2}\cr\cr g^{'}(2t-1)&&si\  {1\over 2}\leq t\leq
1}\right.$}\par
\vskip 2mm Or, pour tout $(x,y)$ dans $R$ on a $x + y = 1$ donc \par\vskip 2mm
 \centerline{$F(x,y) = \left\{\matrix{f^{'}(2-2y,2y-1)&&si\  {1\over 2}\leq
y\leq 1\cr\cr
f(1-2y,2y)&&si\  0\leq y\leq {1\over 2}}\right.$}\par\vskip 2mm
On en d\'eduit que $\gamma (\bar F) = \bar G$ , donc $\gamma $ est un
morphisme de groupes et par suite un isomorphisme.\hfill


\par\vskip 5mm\hskip 5mm
{\bf Th\'eor\`eme (2.4.4) :} {\it Soient $(X,x_{_{0}})$ un espace
topologique point\'e, $n$ un
entier naturel non nul. Alors pour tout $i\in \{1,..,n\}$, il existe un
isomorphisme de
groupes entre $\pi _{_{i}}(\Pi _{_{n}}(X),x_{_{0}})$ et $\pi
_{_{i}}(X,x_{_{0}})$.}
\par\vskip 5mm\hskip 5mm{\bf Preuve :} D'apr\`es les propositions (2.4.1) et
(2.4.2),
et en rempla\c cons $N$ par $I_{_{i-1}}$ on obtient les bijections
suivantes :\par\vskip 3mm
\centerline{$\matrix{\pi _{_{1}}\bigl(X_{_{N}},Id_{_{N}}(x_{_{0}})\bigr)
&\simeq&\bigl[(S^{1},s_{_{0}}); (X_{_{N}},Id_{_{N}}(x_{_{0}})\bigr]\cr\cr
&\simeq&\bigl[(R,\delta R); (X_{_{N}},Id_{_{N}}(x_{_{0}})\bigr]\cr\cr
&\simeq&\bigl[(R^{^{N,1}},\delta R^{^{N,1}}); (X,x_{_{0}})\bigr]\cr \cr
&\simeq&\bigl[(S^{i},s_{_{0}}); (X,x_{_{0}})\bigr]\cr \cr
&\simeq&\pi _{_{i}}(X,x_{_{0}})}$}\par	\vskip 3mm
qui sont des isomorphismes de groupes par transf\`ere de la strucrure de
$\pi _{_{i}}(\Pi _{_{n}}(X),x_{_{0}})$. On d\'eduit alors d'apr\`es la
proposition (2.6.3) qu'il
existe un isomorphisme entre $\pi _{_{i}}(\Pi _{_{n}}(X),x_{_{0}})$ et $\pi
_{_{i}}(X,x_{_{0}})$.
\hfill \par\vskip 5mm


{\bf (2.5).---R\'ealisation g\'eom\`etrique d'un $n$-groupoide :}\par\vskip
5mm\hskip 5mm
Soit $\Phi $ un $n$-groupoide. Pour toute fl\`eche $\sigma :M^{'}\fhd{}{}M$ de
$\Delta ^{n}$ on a les deux applications :\par
\centerline{$\diagram{\Phi (M)&\Fhd{\sigma ^{'}=\Phi (\sigma )}{}&\Phi
(M^{'})\hskip 2cm
R^{^{M^{'}}}&\Fhd{\sigma ^{"}=R^{\sigma }}{}&R^{^{M}}}$}
Consid\'erons l'ensemble : $\bar \Phi :=\displaystyle\bigcup _{M\in
ob(\Delta ^{n})}
R^{^{M}}\times \Phi (M)$. De mani\`ere naturelle on d\'efinit sur cet
ensemble la relation
qui, pour tout $(x,y)$ dans $R^{^{M^{'}}}\times \Phi (M)$, identifie
l'\'el\'ement
$(x,\sigma ^{'}(y))$ de $ R^{^{M^{'}}}\times \Phi (M^{'})$ \`a l'\'el\'ement
$(\sigma ^{"}(x),y)$ de $ R^{^{M}}\times \Phi (M)$. Cette relation est une
relation
d'\'equivalence, et on note par $\mid \Phi \mid$ l'ensemble de ses classes
d'\'equivalences.
$\mid \Phi \mid$ est munit de fa\c con naturelle d'une structure d'espace
topologique,
o\`u $\Phi (M)$ est munit de la topologie discr\`ete (si le $n$-groupoide
$\Phi $ est \`a
valeurs dans la cat\'egorie  ${\cal T}op$ des espaces topologiques on prend
dans ce cas la
topologie de $\Phi (M)$). dans la suite on \'ecrira une classe dans $\mid
\Phi \mid$
par $\mid (x,a)_{_{M}}\mid$ o\`u  $(x,a)$ est un repr\'esentant dans
$ R^{^{M}}\times \Phi (M)$.
\par\vskip 5mm\hskip 5mm Soit $F:\Phi \fhd{}{}\Psi $ un morphisme entre
deux $n$-groupoides, alors $F$ induit une application continue :\par
\centerline{$\diagram{\mid \Phi \mid&	\Fhd{\mid F\mid}{}&\mid \Psi \mid\cr
\mid (x,a)_{_{M}}\mid&	\Fhd{}{}&\mid (x,F_{_{M}}(a))_{_{M}}\mid}$}
Cette application ne d\'epend pas du repr\'esentant $(x,a)_{_{M}}$, en
effet :\par
si $(x,y)$ est un \'el\'ement de $R^{^{M}}\times \Phi (M^{'})$ et
$\sigma :M\fhd{}{}M^{'}$ une fl\`eche de $\Delta ^{n}$ on a les relations :\par
\centerline{$\mid (x,\sigma ^{'}(y))_{_{M}}\mid =\mid (\sigma
^{"}(x),y)_{_{M^{'}}}\mid
\hskip 1cm et \hskip 1cm\sigma ^{'}F_{_{M^{'}}}=F_{_{M}}\sigma ^{'}$}\par
Alors :\par\vskip 2mm
\centerline{$\matrix{\mid F\mid\Big[\mid (x,\sigma ^{'}(y))_{_{M}} \mid\Big]
&=&\mid (x,F_{_{M}}\sigma ^{'}(y))_{_{M}}\mid\cr\cr
&=&\mid (x,\sigma ^{'}F_{_{M^{'}}}(y))_{_{M}}\mid\cr	\cr
&=&\mid (\sigma ^{"}(x),F_{_{M^{'}}}(y))_{_{M^{'}}}\mid\cr\cr
&=&\mid F \mid\Big[\mid (\sigma ^{"}(x),y)_{_{M^{'}}}\mid\Big]}$}\par\vskip 2mm
D'autre part la continuit\'e des applications $F_{_{M}}$ entraine celle de
$\mid F\mid$.
\par	\vskip 2mm\hskip 5mm
On vient de construire un foncteur $\mid\ \mid$ : $n$-Gr $\fhd{}{}$ ${\cal
T}op$
de la cat\'egorie des $n$-groupoides  et transformations naturelles vers
celle des espaces
topologiques. Au paragraphe 3 chapitre 2, nous avons construit un foncteur
dans le sense
inverse $\Pi_{_{n}}(\ ) : {\cal T}op\ \fhd{}{}$ $n$-Gr qui \`a un espace
topologique fait
associer son $n$-groupoide de Poincar\'e. \vskip 5mm\hskip 5mm
{\bf Proposition (2.5.1) :} {\it Soient $\Phi $ un $n$-groupoide et $X$ un
espace topologique.
Alors il existe une application naturelle :\par\vskip 2mm
\centerline{${\cal F}\ :\ Hom\bigl(\mid \Phi \mid,X\bigr)\ \Fhd{}{}\
Hom\bigl(\Phi ,\Pi _{_{n}}(X)\bigr)$}}\par
\vskip 5mm\hskip 5mm{\bf Preuve :} Soient $f : \mid \Phi \mid \fhd{}{}X$
une application
continue et $M$ un \'el\'ement de $\Delta ^{n}$. Pour tout $a$ dans $\Phi
(M)$ on d\'esigne
par ${\cal L}_{a}$ l'application de $R^{^{M}}$ vers $R^{^{M}}\times \Phi
(M)$ qui \`a $x$
fait correspondre le couple $(x,a)$. On d\'efinit ${\cal F}(f)$ par :\par
\centerline{$\diagram{\Phi (M)&\Fhd{{\cal F}(f)(M)}{}&\overline{X_{_{M}}}\cr
a&\Fhd{}{}&\overline{\Gamma _{a}}}$}
o\`u $\Gamma _{a}$ est la compos\'e des trois applications suivantes :	\par
\centerline{$\diagram{R^{^{M}}&\Fhd{{\cal L}_{a}}{}&R^{^{M}}\times \Phi (M)&
\Fhd{S_{_{M}}}{}&\mid\Phi \mid&\Fhd{f}{}X}$}
et $S_{_{M}}$ l'application canonique qui envoie un \'el\'ement vers sa
classe dans
$\mid\Phi \mid$. Montrons maintenant que ${\cal F}(f)$ est un morphisme.
\par\hskip 5mm Soit $\sigma : M^{'}\fhd{}{}M$ une fl\`eche de $\Delta
^{n}$, on veut
montrer la commutativit\'e du diagramme suivant :\par\centerline{$
\diagram{\Phi (M)&\Fhd{{\cal F}(f)(M)}{}&\overline{X_{_{M}}}\cr \fvb{\sigma
^{'}}{}&&
\fvb{}{\bar \sigma }\cr\Phi (M^{'})&\Fhd{}{{\cal
F}(f)(M^{'})}&\overline{X_{_{M^{'}}}}}$
\hskip 1cm o\`u \hskip 1cm$\sigma ^{'}=\Phi (\sigma )$ et $\bar \sigma =
\Pi _{_{n}}(X)(\sigma )$}
Pour tout $(x,a)$ dans $R^{^{M^{'}}}\times \Phi (M)$ on a les relations
suivantes :	\par
\vskip 2mm\centerline{$\matrix{\Big[\bar{\sigma }
\circ {\cal F}(f)(M)\Big](a)=\bar{\sigma }\Big[\overline{\Gamma _{a}}\Big]=
\overline{\Gamma _{a}\circ \sigma ^{"}}\cr\cr
\Big[{\cal F}(f)(M)\circ \sigma ^{'}\Big](a)=\overline{\Gamma _{ \sigma
^{'}(a)}}\cr\cr
\mid (x,\sigma ^{'}(a))_{_{M^{'}}}\mid =\mid (\sigma ^{"}(x),a)_{_{M}}\mid}$}
\par\vskip 2mm On en d\'eduit :\par\vskip 2mm\centerline{$\matrix{
\Gamma _{ \sigma ^{'}(a)}(x)&=&f\Big(\mid (x,\sigma
^{'}(a))_{_{M^{'}}}\mid\Big)\cr\cr
&=&f\Big(\mid (\sigma ^{"}(x),a)_{_{M}}\mid\Big)\cr	\cr
&=&\Gamma _{a}(\sigma ^{"}(x))\cr	\cr
&=&[\Gamma _{a}\circ\sigma ^{"}](x)}$}\par\vskip 4mm
ce qui montre que ${\cal F}(f)$ est un morphisme.\hfill \par\hskip 5mm
Si on prend dans la proposition (2.5.1),  $X=\mid \Phi \mid$ et
$f=Id_{_{\mid \Phi \mid}}$,
on obtient une morphisme de $\Phi $ vers $\Pi _{_{n}}(\mid \Phi \mid)$.
Nous conjecturons que le foncteur $\mid\ \mid$ est un inverse \`a
\'equivalence pr\`es
du foncteur $\Pi _{_{n}}(\ )$ de la cat\'egorie des espaces topologiques
$n$-tronqu\'es
vers la cat\'egorie $n$-Gr des $n$-groupoides.\par\vfill\eject


\centerline{\bf SOMMAIRE}\par\vskip 2cm
{\bf ITRODUCTION}\hfill 1\par\vskip 6mm
{\bf Chapitre 1 :\hskip 3mm NOTION DE $n$-NERF}
\par\vskip 6mm\hskip 2cm
{\bf (1.1) Notations et d\'efinitions }\hfill 1\par\vskip 3mm\hskip 2cm
{\bf (1.2) Troncation d'une $n$-pr\'e-Nerf}\hfill 5\par\vskip 3mm\hskip 2cm
{\bf (1.3) $n$-\'equivalence et $n$-Nerf }\hfill 11\par\vskip 3mm\hskip 2cm
{\bf (1.4) 2-cat\'egorie usuelle et 2-Nerf}\hfill 17\par\vskip 3mm
{\bf Chapitre 2 :\hskip 3mm NOTION DE $n$-GROUPOIDE}\par\vskip 3mm\hskip 2cm
{\bf (2.1) Notations et d\'efinitions}\hfill 42\par\vskip 3mm\hskip 2cm
{\bf (2.2) Propriet\'es des $n$-groupoides}\hfill 44\par\vskip 3mm\hskip 2cm
{\bf (2.3) $n$-Groupoide de Poincar\'e}\hfill 51
\par\vskip 3mm\hskip 2cm
{\bf (2.4) Groupes d'homotopies de $\Pi _{_{n}}(X)$}\hfill 60
\par\vskip 3mm\hskip 2cm
{\bf (2.5) R\'ealisation g\'eom\`etrique d'un $n$-groupoide}\hfill
64\par	\vfill\eject


 \vskip 1cm{\bf  REFERENCES BIBLIOGRAPHIQUES}\vskip 1cm
[1] A. Grothendieck, A la poursuite des champs, unpublished manuscript.
\par\vskip 3mm
[2] J. Giraud, Cohomologie non ab\'elienne, Grundelehren der mathematischen
Wissenschaften in Einzeldarstellungen 179, Springer-Verlag (1971).\par\vskip
3mm
[3] L. Breen, On the classification of 2-gerbes and 2-stacks, Ast\'erisque 225,
Soc. Math. de France (1994). \par\vskip 3mm
[4] O. Leroy, Sur une notion de 3-cat\'egorie adapt\'ee \`a l'homotopie,
preprint
Universit\'e de Montpellier 2 (1994).\par\vskip 3mm
[5] P. Gabriel-M. Zisman, Calculus of fractions and homotopy theory. Ergebnisse
der mathematik, Vol. 35. Berlin-Heidelberg-New York : Springer-Verlag (1967).
\par\vskip 3mm
[6] S. Mac Lane, Categories for the working mathematicians, Spinger-Verlag,
New York (1971).\par \vskip 3mm
[7] S. Mac Lane, Sheaves in geometry and logic, Spinger-Verlag,
New York (1992).\par \vskip 3mm
[8] P. May, Simplicial objects in algebraic topology. Princeton : Van Nostrand
mathematical studies 11, (1967).\par\vskip 3mm

\end

\end